\documentclass[12pt]{article}
\usepackage{amsmath}    
\usepackage{amssymb}    
\usepackage{amsfonts}   
\usepackage{amsthm}     
\usepackage{geometry}   
\usepackage{graphicx}   
\usepackage{enumitem}   
\usepackage{indentfirst}
\usepackage{tikz-cd}    
\usepackage{caption}    
\usepackage{mathrsfs}   
\usepackage[normalem]{ulem}  
\usepackage{color}      
\usepackage{mathtools}  
\usepackage{tcolorbox}  
\usepackage[hidelinks]{hyperref}
\usepackage{doi}

\usepackage[all]{xy}

\usepackage[numbers]{natbib}
\providecommand{\keywords}[1]{\textbf{Keywords:} #1}
\providecommand{\subjclass}[1]{\textbf{2020 Mathematics Subject Classification:} #1}
\newtheorem{definition}{Definition}[section]
\newtheorem{theorem}[definition]{Theorem}
\newtheorem{proposition}[definition]{Proposition}
\newtheorem{lemma}[definition]{Lemma}
\newtheorem{example}[definition]{Example}
\newtheorem{remark}[definition]{Remark}

\newtheorem{corollary}[definition]{Corollary}
\numberwithin{equation}{section}

\begin{document}
	
	\title{Topological Lifting of  Mundici Equivalence \\
 and Di Nola--Lettieri Equivalence
}
	\author{Yanqiao Jia
		, Yichuan Yang\thanks{Corresponding author. E-mail: ycyang@buaa.edu.cn}}
	\date{}
	
	\maketitle

\vspace{-0.7em}
\begin{center}
\textit{Dedicated to the memory of Prof.~K.~Tanaka (1955-2026)}
\end{center}
\vspace{0.1em}

\begin{abstract}
We prove that \emph{summability}
characterises precisely when a topological MV-algebra can be realized as
the unit interval of a lattice-ordered topological abelian group, thereby
providing a topological lifting of both  Mundici equivalence and Di Nola--Lettieri equivalence. Applications
include a Birkhoff--Kakutani-type metrisation theorem, a
second-countability criterion, and a 
topological extension of a key
continuity lemma from the proof theory of first-order \L{}ukasiewicz
logic---following Baaz--Metcalfe's hypersequent approach to approximate
completenesss---and provides an $A$-valued topological semantics for the approximate strong completeness of~$\mathbf{G}\text{\L}\forall$, where $A$ ranges over 
 summable topological MV-algebras.
\end{abstract}


	\medskip

\section*{Introduction}\label{sec:introduction}

\smallskip
\noindent\textbf{Topological lifting
problem.}
MV-algebras were introduced by Chang~\cite{chang1958} and proved complete for
{\L}ukasiewicz logic in~\cite{chang1959}. Two foundational categorical
equivalences shape their modern theory. Mundici~\cite{mundici1986,mundici1986joa}
established the categorical equivalence between MV-algebras and abelian
lattice-ordered groups ($\ell$-groups for short) with a distinguished strong
unit. 
Di Nola and Lettieri~\cite{dinola1994} established a categorical
equivalence between perfect MV-algebras and abelian $\ell$-groups;
Belluce and Di Nola~\cite{belluce1996} subsequently developed Yosida-type
functional representations for perfect MV-algebras using these equivalences.
These categorical bridges have enabled the transfer of structural results between
MV-algebras and abelian $\ell$-groups for nearly four decades
\cite{caramello2015,mundici2000,wehrung2019,yang2006}%
\footnote{Caramello and Russo~\cite{caramello2015} gave a topos-theoretic
recasting of Mundici's equivalence, elevating it to a Morita equivalence
of classifying topoi.}. 

These two equivalences are cornerstones of many-valued logic and abelian $\ell$-groups. Yet the question of whether—and under what intrinsic condition—they lift to the topological setting has remained unsettled for decades:\smallskip  \\ 
\emph{Can  Mundici
equivalence and  Di Nola--Lettieri equivalence  be
lifted to the topological setting?}\\ 
\noindent More concretely, given a topological
MV-algebra $(A,\tau)$, does there exist a group topology on the Mundici
group $G_A$ whose restriction to $[0,u_{G_A}]$ recovers $\tau$? And similarly
on  Di Nola--Lettieri side? 

\medskip
\noindent\textbf{Duality theories for MV-algebras and $\ell$-groups.}
Cignoli, Dubuc, and Mundici~\cite{cignoli2004} extended Stone duality to
multisets and locally finite MV-algebras, providing a dual equivalence that
mirrors the classical Stone representation of Boolean algebras.
On the $\ell$-group side, Rump and Yang~\cite{rumpyang2008} studied the
Jaffard--Ohm correspondence between abelian $\ell$-groups and B\'ezout
domains and extended it to sheaves on spectral spaces, exhibiting the
spectrum of a B\'ezout domain and that of its divisibility $\ell$-group as
a concrete instance of Hochster duality. Later~\cite{rumpyang2009}
constructed the structure sheaf of an archimedean $\ell$-group on its prime
spectrum and recovered the lateral completion from it. Independently,
Dubuc and Poveda~\cite{dubucpoveda2010} represented arbitrary MV-algebras as
global sections of a sheaf of MV-chains on a compact space, generalising
McNaughton's theorem; part of the sheaf-theoretic viewpoint had been
prefigured in Yang~\cite{yang2006}.
These dualities live on the MV/$\ell$-group side of the quantum-B-algebra
spectrum introduced by Rump and Yang~\cite{rump2014}; for the broader
class of unital lattice-ordered quantum B-algebras, Zariski spectra of prime
deductive systems have recently been analysed by Mu, Li, Jia and
Yang~\cite{mu2026qb}, under a distributivity assumption on the deductive-system
lattice.
This line of work highlights the central role of compatibility conditions
between algebraic operations and topology in yielding dual equivalences and
sheaf representations, a theme to which the present paper contributes a
topological lifting criterion for Mundici and Di Nola--Lettieri equivalences.

\medskip
\noindent\textbf{Topological ordered algebraic systems.}
The metrization of topological groups was settled independently by
Birkhoff~\cite{birkhoff1936} and Kakutani~\cite{kakutani1936}, who proved
that a topological group is metrizable if and only if the identity admits a
countable neighbourhood basis; moreover, the metric can be chosen
left-invariant. Kakutani's representation theorems for abstract $L$- and
$M$-spaces in~\cite{kakutani1941A,kakutani1941A2} built the bridge between
order-theoretic structures and spaces of continuous functions, a line of
thought that culminated in Yosida's representation theorem for unital
$\ell$-groups~\cite{yosida1941unit}. These works, together with the theory
of topological groups and the Birkhoff--Kakutani metrization theorem,
provide a foundation for understanding the interaction between algebraic
order and topology.

\medskip
\noindent\textbf{The topological study of MV-algebras.}
The topological study of MV-algebras dates back to Hoo's pioneering
work~\cite{hoo1997}, who showed that $T_0,T_1,T_2$ coincide for
topological MV-algebras. Gan, Luan, Deng and Yang~\cite{gan2024}
later proved that every topological MV-algebra is regular ($T_3$),
hence $T_0,T_1,T_2,T_3$ are equivalent; Li and Yang~\cite{li2026}
obtained $T_{3\frac12}$ (complete regularity) under local convexity.
On the side of neighbourhood characterisations, Weber~\cite{weber2012}
proved that the $0$-neighbourhood system determines the whole
topology. Luan, Weber and Yang
analysed when a filter topology makes an MV-algebra into a topological
MV-algebra (\cite[Theorem~3.5]{luanweberyang2021}), showing in particular that every MV-chain
with a filter topology is a topological MV-algebra (\cite[Theorem~4.7]{luanweberyang2021}). 
  Luan, Zhao and Yang~\cite{luan2021}
characterised Hausdorff filter-topological MV-algebras. 
All these filter-topological results rely on the fact that every filter
topology is locally convex (Proposition~3.9 of \cite{luanweberyang2021}). For a general topological MV-algebra, however, local convexity
is not automatic; a concrete counterexample is provided by
Example~\ref{ex:non-locally-convex} below.  This distinction motivates seeking an intrinsic neighbourhood condition---
one that arises naturally from the algebraic structure---that is weaker
than local convexity yet strong enough to force an $\ell$-group lifting.

\medskip
\noindent\textbf{Topological lifting within the duality framework.}
The programme of lifting algebraic dualities and categorical equivalences
of logical algebras to the topological setting has attracted sustained
attention.  For MV-algebras themselves, Dubuc and Poveda~\cite{dubucpoveda2010}
gave a sheaf-theoretic representation via classifying topoi, while
Rump and Yang~\cite{rump2014} placed them inside the quantale duality for
non-commutative logical algebras.  More recently, for semimodules over an incline,
Bai and Yang~\cite{bai2020} established an algebraic--topological duality
together with a sheaf representation, thereby lifting the structure theory
of such semimodules to the topological category.  The present paper
resolves this question by introducing \emph{summability} as the intrinsic
neighbourhood condition that captures precisely when a topological
MV-algebra arises as the unit interval of a lattice-ordered topological
abelian group.  We prove that summability characterises precisely when a
topological MV-algebra can be realized as the unit interval of a
lattice-ordered topological group, thereby resolving the topological
lifting problem for both Mundici equivalence and Di Nola--Lettieri
equivalence.  Applications include a Birkhoff--Kakutani-type metrisation
theorem, a second-countability criterion, and a topological extension of a
key continuity lemma from the proof theory of first-order {\L}ukasiewicz
logic---following Baaz--Metcalfe's hypersequent approach to approximate
completeness~\cite{baaz2010}---a completeness result recently revisited by
Duo~\cite{duo2026} via labelled tableaux, though still confined to the
standard $[0,1]$-semantics, whereas here the semantic carrier is shifted
to arbitrary summable topological MV-algebras.

The main results of this paper are as follows.
\begin{itemize}
    \item \textbf{Neighbourhood basis characterisation
    (Proposition~\ref{thm:MV-topology}).} We give a full characterisation of which filter bases can serve as a
  neighbourhood basis of $0$ in a topological MV-algebra.
  Restricting our condition to open filter bases recovers, as a corollary,
  the neighbourhood criterion of Gan et al.~\cite[Proposition~3.5]{gan2024}, which treats arbitrary (not necessarily filter-generated)   $0$-neighbourhood families.

    \item \textbf{Summability and the topological lifting theorem (Theorem~\ref{thm:iffMV}, Proposition~\ref{prop:MV-lattice-group-condition}).}
We introduce the notion of a \emph{summable} topological MV-algebra and prove
that a topological MV-algebra is homeomorphically isomorphic to the unit
interval of a lattice-ordered topological group if and only if it is summable.
Conversely, we characterize when the unit
interval of a lattice-ordered topological group is a topological MV-algebra.

    \item \textbf{Categorical adjunction
    (Theorem~\ref{thm:adjunction}).}
We establish an adjunction between the category
of summable topological MV-algebras and the category of MV lattice-ordered
topological groups, providing a categorical solution to the topological lifting
problem for Mundici equivalence. Furthermore, we identify the subcategories in which the
		adjunction restricts to strongly summable topological MV-algebras and
		topological \(\ell\)-groups (Corollary \ref{cor:strong-adjunction}), and show
		that in the totally ordered case it becomes an equivalence of categories
		(Corollary \ref{cor:totally-ordered-equivalence}).

    \item \textbf{Topological Di Nola--Lettieri equivalence (Theorem~\ref{thm:perfect-topological-equivalence}).}
We prove that Di Nola--Lettieri equivalence lifts to an equivalence of
topological categories. 
\end{itemize}

\medskip
\noindent\textbf{Topological consequences.}
As consequences of summability, we prove the following:
\begin{itemize}
\item \textbf{Separation axioms:} summability implies that
  $T_0$ entails complete regularity for totally ordered
  topological MV-algebras (Proposition~\ref{prop:separation}).
  This gives, in the totally ordered case, an independent route to the
  complete regularity established by Li and Yang~\cite{li2026} for
  arbitrary locally convex topological MV-algebras---and, notably,
  without assuming local convexity, which even totally ordered
  topological MV-algebras need not possess
  (Example~\ref{ex:non-locally-convex}).

    \item \textbf{Metrisability:} we establish a Birkhoff--Kakutani-type
    metrisation theorem for totally ordered topological
    MV-algebras (Proposition~\ref{prop:metrization}).

    \item \textbf{Second countability:} a totally ordered Hausdorff
    topological MV-algebra is second-countable if and only if
    it is separable and metrisable
    (Corollary~\ref{cor:second-countable}).

    \item \textbf{Approximate completeness:} we give a new proof
    of a key lemma in the approximate strong completeness of
    the hypersequent calculus $\mathbf{G}\text{\L}\forall$ \cite{baaz2010}
    and extend its continuity ingredient to arbitrary summable
    topological MV-algebras
    (Proposition~\ref{prop:ev-continuous}). This provides
    the main topological ingredient for an $A$-valued semantic
    version of approximate strong completeness for
    $\mathbf{G}\text{\L}\forall$ , where $A$ ranges over all summable
    topological MV-algebras, thereby going beyond the standard
    $[0,1]$-valued setting.
\end{itemize}

\medskip
The paper is organised as follows. 
Section~\ref{sec:preliminaries} recalls the necessary background.
Section~\ref{sec:0-neighbor} establishes the neighbourhood basis
characterisation for $0$ in a topological MV-algebra.
Section~\ref{sec:summable} introduces summability and proves the
topological lifting theorem.
Section~\ref{sec:adjunction} develops the categorical framework,
including the adjunction, the restriction to strong summability, and
the topological Di Nola--Lettieri equivalence.
Section~\ref{sec:application} presents topological consequences of
summability, including separation axioms, metrisability, second
countability, and approximate completeness.
Section~\ref{sec:conclusion} concludes the paper with a summary of the
main contributions and a discussion of open problems and future
directions.

\section{Preliminaries}\label{sec:preliminaries}
	In this section we recall the necessary background and fix notation that will be used throughout the paper.
	
	\begin{definition}[{\cite{chang1958,
				mundici2000}}]
		An MV-algebra is an algebra \((A,\oplus,\neg,0)\) of type
		\((2,1,0)\) such that, for all \(x,y\in A\),
		\begin{itemize}
			\item[(M1)] \((A,\oplus,0)\) is a commutative monoid;
			\item[(M2)] \(\neg\neg x=x\);
			\item[(M3)] \(x\oplus\neg0=\neg0\);
			\item[(M4)] \(\neg(\neg x\oplus y)\oplus y
			=\neg(\neg y\oplus x)\oplus x\).
		\end{itemize}
	\end{definition}
	
	Set \(1:=\neg0\), \(x\odot y:=\neg(\neg x\oplus\neg y)\), and
	\(x\ominus y:=x\odot\neg y\). Then the following identities hold:
	\begin{itemize}
		\item \(\neg1=0\);
		\item \(x\oplus y=\neg(\neg x\odot\neg y)\);
		\item \(x\oplus1=1\);
		\item \(x\oplus\neg x=1\);
		\item \(x\odot\neg x=0\).
	\end{itemize}
	
	The \emph{natural order} of an MV-algebra is defined by
	\(x\leq y\) if and only if \(\neg x\oplus y=1\). Then
	\((A,\leq)\) is a bounded distributive lattice, with bottom element
	\(0\), top element \(1\), and
	\[
	x\vee y=y\oplus(x\odot\neg y),
	\qquad
	x\wedge y=x\odot(\neg x\oplus y).
	\]
	
	A subset \(I\subseteq A\) is called an \emph{ideal} of \(A\) if:
	\begin{enumerate}[label=(\arabic*)]
		\item \(0\in I\);
		\item if \(x,y\in I\), then \(x\oplus y\in I\);
		\item if \(x\in I\) and \(y\leq x\), then \(y\in I\).
	\end{enumerate}
	A proper ideal \(M\) of \(A\) is called \emph{maximal} if it is
	maximal, with respect to inclusion, among the proper ideals of \(A\).
	
	An MV-algebra \(A\) is called \emph{totally ordered} if its natural
	order is total, that is, for every \(x,y\in A\), either \(x\leq y\)
	or \(y\leq x\). If \(A\) is totally ordered, then
	\cite[Theorem~1.6.1]{mundici2000}
	\begin{equation}\label{prop:chain-property}
		x\oplus y<1 \quad\Longrightarrow\quad x\odot y=0.
	\end{equation}

	\begin{example}[\cite{mundici2000}]
		The following are examples of MV-algebras:
		\begin{enumerate}[label=(\arabic*)]
			\item The unit interval $[0,1]$ with $x \oplus y = \min(x+y,1)$ and $\neg x = 1-x$ is called the \emph{standard MV-algebra}. In this case, $x \odot y = \max(x+y-1,0)$ and $x \ominus y = \max(x-y,0)$.
			\item For any integer $n \ge 2$, the set $L_n = \{0, \frac{1}{n-1}, \frac{2}{n-1}, \cdots, 1\}$ with the same operations is a finite subalgebra of the standard MV-algebra.
		\end{enumerate}
	\end{example}

	\begin{theorem}[Chang's Completeness Theorem] \label{thm:chang}
		An equation holds in all MV-algebras if and only if it holds in the standard MV-algebra $[0,1]$.
	\end{theorem}
For a proof, see \cite[Theorem 2.5.3]{mundici2000}.
	
	\begin{lemma}\label{lem:MV-identity}
		For any $x, y, c$ in an MV-algebra, the following identity holds:
		\[
		x \ominus y = ((x \ominus c) \oplus (c \ominus y)) \ominus ((c \ominus x) \oplus (y \ominus c)).
		\]
	\end{lemma}
	\begin{proof}
		By Chang's completeness theorem (Theorem~\ref{thm:chang}), it
		suffices to verify the identity in $[0,1]$, where
		$a\oplus b=\min\{1,a+b\}$ and
		$a\ominus b=\max\{0,a-b\}$. Since the two inner sums are at most
		$1$, the right-hand side equals
		\[
		\begin{aligned}
			\max\{0,\,&\max\{0,x-c\}+\max\{0,c-y\}\\
			&-\max\{0,c-x\}-\max\{0,y-c\}\}
			=\max\{0,x-y\},
		\end{aligned}
		\]
		where we used
		$\max\{0,t\}-\max\{0,-t\}=t$. Hence it equals $x\ominus y$.
	\end{proof}
	
	A McNaughton function on $[0,1]^n$ is a continuous piecewise linear function
	$f:[0,1]^n\to[0,1]$ whose linear pieces have integer coefficients.
	
	\begin{theorem}[McNaughton's theorem]\label{thm:mcnaughton}
		A function \(f:[0,1]^n\to[0,1]\) is induced by an MV-term if and only if it is a
		McNaughton function.
	\end{theorem}
	For a proof, see \cite[Theorem 9.1.5]{mundici2000}.

	We now recall some basic notions from the theory of lattice-ordered groups. For a comprehensive treatment, see \cite{birkhoff1967, darnel1994}.
	
	A \emph{lattice-ordered group} ($\ell$-group for short) is a group $(G, +, 0)$ equipped with a partial order $\le$ that is a lattice and such that the group operation is translation-invariant: for all $x, y, z \in G$,
	\[
	x \le y \implies x + z \le y + z \quad \text{and} \quad z + x \le z + y.
	\]
	For an $\ell$-group $(G, +, 0, \leq)$, we denote the lattice operations by $x \vee y$ and $x \wedge y$. We write
	\[
	x^+ = x \vee 0, \qquad x^- = (-x) \vee 0, \qquad |x| = x^+ + x^-.
	\]
	
	An \emph{$\ell$-group homomorphism} is a group homomorphism $f: G \to H$ that also preserves the lattice operations, i.e., $f(x \vee y) = f(x) \vee f(y)$ and $f(x \wedge y) = f(x) \wedge f(y)$. 
	
	A subset $I \subseteq G$ is called an \emph{$\ell$-ideal} if it is a subgroup of $G$ and for any $x \in I$ and $y \in G$, $|y| \le |x|$ implies $y \in I$. Equivalently, $I$ is a convex subgroup closed under the lattice operations.
	
	A positive element $u \in G$ is called a \emph{strong unit} if for every $x \in G$ there exists $n \in \mathbb{N}$ such that $x \le n u$. We denote by $(G, u)$ an $\ell$-group with a distinguished strong unit $u$.
	
	For example, the additive group $\mathbb{R}$ with the usual order is an abelian $\ell$-group, and $1$ is a strong unit.
	
	Throughout this paper,  we only consider abelian $\ell$-groups.
	\begin{proposition}[\cite{birkhoff1967, darnel1994}]\label{prop:lgroup-basic}
		Let $(G, +, 0, \leq)$ be an abelian $\ell$-group and $x, y, z \in G$. The following hold:
		\begin{enumerate}[label=(\arabic*)]
			\item $x + (y \vee z) = (x + y) \vee (x + z)$ and $x + (y \wedge z) = (x + y) \wedge (x + z)$;
			\item $-(x \vee y) = (-x) \wedge (-y)$ and $-(x \wedge y) = (-x) \vee (-y)$;
			\item $x = x^+ - x^-, |x| = x \vee (-x) = x^+ \vee x^-$;
			\item If $n x \ge 0$ for some positive integer $n$, then $x \ge 0$.
			\item $(G, +, 0, \leq)$ is a subdirect product of totally ordered
			abelian \(\ell\)-groups.
		\end{enumerate}
	\end{proposition}
	
	\begin{lemma}\label{lem:lgroup-identity}
		Let $(G, +, 0, \leq)$ be an abelian $\ell$-group and let $u \in G$ be a positive element. Then for any $a \in G$, the following holds:
		\[
		a = a \wedge u + (a - u) \vee 0.
		\]
	\end{lemma}
	
	\begin{proof}
		By Proposition \ref{prop:lgroup-basic}, we have $a - u = (a - u)^+ - (a - u)^- = (a - u) \vee 0 - (u - a) \vee 0 =  (a - u) \vee 0 + (a - u) \wedge 0 $. Hence $a =(a - u) + u = (a - u) \vee 0 + (a - u) \wedge 0 + u = a \wedge u + (a - u) \vee 0$.
	\end{proof}
	
	\begin{lemma}\label{lem:lgroup-positive-part-estimate}
		Let $(G,+,0,\leq)$ be an abelian $\ell$-group. If \(c\leq d\), then
		\[
		d^+-c^+=(d-c)\wedge d^+.
		\]
		In particular, \(0\leq d^+-c^+\leq d-c\). Consequently, for any \(y\in G\)
		and \(a,b\in G^+\), we have
		\[
		\bigl((y+a-b)^+-y^+\bigr)^+\leq a,
		\qquad
		\bigl((y+a-b)^+-y^+\bigr)^-\leq b.
		\]
	\end{lemma}
	
	\begin{proof}
		The identity \(d^+-c^+=(d-c)\wedge d^+\) is an \(\ell\)-group
		identity. By Proposition \ref{prop:lgroup-basic}(5), it is enough to verify it
		in the totally ordered case, where it follows immediately by considering the
		three cases \(d\leq0\), \(c\leq0\leq d\), and \(0\leq c\leq d\).
		The inequality \(0\leq d^+-c^+\leq d-c\) follows at once.
		
		Applying this inequality to \((c,d)=(y-b,y+a-b)\) and \((c,d)=(y-b,y)\), we
		get \(0\leq (y+a-b)^+-(y-b)^+\leq a\) and
		\(0\leq y^+-(y-b)^+\leq b\). Put
		\[
		\alpha=(y+a-b)^+-(y-b)^+,\qquad
		\beta=y^+-(y-b)^+.
		\]
		Then \(0\leq\alpha\leq a\), \(0\leq\beta\leq b\), and
		\((y+a-b)^+-y^+=\alpha-\beta\). Since
		\((\alpha-\beta)^+\leq\alpha\) and \((\alpha-\beta)^-\leq\beta\), the
		conclusion follows.
	\end{proof}
	
	\begin{example}[\cite{mundici2000}]\label{ex:lgroup-to-MV}
		Let $(G, +, 0, \leq)$ be an abelian $\ell$-group and $u \in G$ a positive element. Then the interval $[0,u] = \{ x \in G \mid 0 \le x \le u \}$ becomes an MV-algebra with the operations
		\[
		x \oplus y = (x + y) \wedge u, \qquad \neg x = u - x.
		\]
		In this case, $x \odot y = (x + y - u) \vee 0$ and $x \ominus y = (x - y) \vee 0$, and the natural order of the MV-algebra coincides with the order inherited from $G$.
		
		In particular, taking $G = \mathbb{R}$ and $u = 1$ yields the standard MV-algebra $[0,1]$.
	\end{example}
	
	Let $\ell\mathcal{G}_u$ denote the category whose objects are abelian $\ell$-groups with a distinguished strong unit $(G, u)$ and whose morphisms are $\ell$-group homomorphisms preserving the strong unit. 
Let $\mathcal{MV}$ denote the category of MV-algebras with MV-algebra homomorphisms preserving the  addition, the negation,  and  the zero.
	
	The construction in Example \ref{ex:lgroup-to-MV} actually defines a functor $\Gamma: \ell\mathcal{G}_u \to \mathcal{MV}$ as follows: for an object $(G, u)$, set $\Gamma(G, u) = ([0,u], \oplus, \neg, 0)$ with operations as in Example \ref{ex:lgroup-to-MV}; for a morphism $f: (G, u) \to (H, v)$, define $\Gamma(f) = f|_{[0,u]}$. 
	
	More generally, if $u$ is merely a positive element of an abelian
	$\ell$-group $G$, we continue to write $\Gamma(G,u)$ for the
	MV-algebra $[0,u]$. In this case, let
	$G_u:=\{g\in G:|g|\leq nu\text{ for some }n\in\mathbb N\}$. Then
	$G_u$ is an $\ell$-ideal of $G$, $u$ is a strong unit of $G_u$, and
	$\Gamma(G,u)=\Gamma(G_u,u)$.
	
	\begin{theorem}[Mundici equivalence,  \cite{mundici1986}]\label{thm:mundici}
		The functor $\Gamma: \ell\mathcal{G}_u \to \mathcal{MV}$ defined above is an equivalence of categories.
	\end{theorem}
	
	In particular, there exists a quasi-inverse functor $\Xi: \mathcal{MV} \to \ell\mathcal{G}_u$, and natural isomorphisms $\varphi: \mathrm{id}_{\mathcal{MV}} \to \Gamma\Xi$ and $\psi: \mathrm{id}_{\ell\mathcal{G}_u} \to \Xi\Gamma$.
	
	Their naturality is expressed by the following commutative diagrams:
	\begin{equation}\label{diag:naturality-varphi-psi}
		\begin{tikzcd}[row sep=large, column sep=large]
			A \arrow[r, "\varphi_A"] \arrow[d, "f"']
			& \Gamma\Xi(A) \arrow[d, "\Gamma\Xi(f)"] \\
			B \arrow[r, "\varphi_B"']
			& \Gamma\Xi(B)
		\end{tikzcd}
		\qquad\qquad
		\begin{tikzcd}[row sep=large, column sep=large]
			(G,u) \arrow[r, "\psi_G"] \arrow[d, "g"']
			& \Xi\Gamma(G,u) \arrow[d, "\Xi\Gamma(g)"] \\
			(H,v) \arrow[r, "\psi_H"']
			& \Xi\Gamma(H,v)
		\end{tikzcd}
	\end{equation}
	for every morphism $f:A\to B$ in $\mathcal{MV}$ and every morphism $g:(G,u)\to(H,v)$ in $\ell\mathcal{G}_u$.
	
	For an MV-algebra $A$, write $\Xi(A) = (G_A, u_{G_A})$ and $\Gamma(G,u) = [0,u]$ when no confusion arises. We call the pair $(G_A, u_{G_A})$ the \emph{Mundici group} of $A$.
	The natural isomorphism $\varphi_A: A \to [0, u_{G_A}]$ satisfies $\varphi_A(\neg x) = u_{G_A} - \varphi_A(x)$ and $\varphi_A(x \oplus y) = (\varphi_A(x) + \varphi_A(y)) \wedge u_{G_A}$ for all $x, y \in A$.
	Also, $A$ is totally ordered if and only if $G_A$ is totally ordered.
	Moreover, for any object $(G, u)$ in $\ell\mathcal{G}_u$ and any object $A$ in $\mathcal{MV}$, we have
	\begin{equation}\label{eq:psi-varphi}
		\varphi_{[0,u]} = \psi_G|_{[0,u]}, \qquad \psi_{G_A}= \Xi(\varphi_A).
	\end{equation}
	For details, see \cite[pp. 139-149]{mundici2000}.
	
	For Mundici equivalence, Cignoli and Torrens proved the following
	
	\begin{theorem}[{\cite{cignoli1996}}]\label{thm:ideal-correspondence}
		Let $A$ be an MV-algebra and let $(G_A, u_{G_A}) = \Xi(A)$. If $I \subseteq A$ is an ideal, then
		\[
		\widehat{\varphi_A(I)} := \left\{ x \in G_A \;\middle|\; |x| \wedge u_{G_A} \in \varphi_A(I) \right\}
		\]
		is an $\ell$-ideal of $G_A$, and
		$\widehat{\varphi_A(I)} \cap [0, u_{G_A}] = \varphi_A(I).$
	\end{theorem}
	
	We shall also use the representation of perfect MV-algebras by abelian
	\(\ell\)-groups established by Di Nola and Lettieri; see
	\cite{mundici2000,dinola1994}. The \emph{radical} of an MV-algebra
	\(A\) is defined by
	\begin{equation}\label{eq:radical}
		\operatorname{Rad}(A)
		:=
		\bigcap\{M:M\text{ is a maximal ideal of }A\}.
	\end{equation}
	
	\begin{definition}[{\cite{mundici2000,dinola1994}}]
		\label{def:perfect-MV}
		An MV-algebra \(A\) is called \emph{perfect} if it is non-trivial and
		\[
		A=\operatorname{Rad}(A)\cup\neg\operatorname{Rad}(A).
		\]
	\end{definition}
	
	For every perfect MV-algebra \(A\), the sets
	\(\operatorname{Rad}(A)\) and \(\neg\operatorname{Rad}(A)\) are
	disjoint. Moreover, it follows directly from perfectness and the
	identity \(x\odot x=0\) if and only if \(x\leq\neg x\) that
	\begin{equation}\label{eq:perfect-radical-square-zero}
		x\in\operatorname{Rad}(A)
		\quad\Longleftrightarrow\quad
		x\odot x=0.
	\end{equation}
	
	Let \(\ell\mathcal G\) denote the category of abelian \(\ell\)-groups
	and \(\ell\)-group homomorphisms, and let \(\mathcal{PMV}\) denote the
	category of perfect MV-algebras and MV-algebra homomorphisms.
	
	For an abelian \(\ell\)-group \(H\), let
	\(\mathbb Z\times_{\mathrm{lex}}H\) denote the direct product group
	\(\mathbb Z\times H\), with coordinatewise addition, endowed with the
	lexicographic order
	\[
	(m,h)\leq(n,k)
	\quad\Longleftrightarrow\quad
	m<n\ \text{or}\ (m=n\ \text{and}\ h\leq k).
	\]
	Then \(\mathbb Z\times_{\mathrm{lex}}H\) is an abelian \(\ell\)-group,
	and \((1,0)\) is a strong unit. Define
	\[
	\Sigma(H)=\Gamma\bigl(\mathbb Z\times_{\mathrm{lex}}H,(1,0)\bigr).
	\]
	Explicitly,
	\(\Sigma(H)=\{(0,h):h\in H^+\}\mathbin{\dot\cup}
	\{(1,-h):h\in H^+\}\), with
	\(\operatorname{Rad}(\Sigma(H))=\{0\}\times H^+\) and
	\(\neg\operatorname{Rad}(\Sigma(H))=\{1\}\times(-H^+)\).
	
	For an \(\ell\)-group homomorphism \(f:H\to K\), define
	\(\Sigma(f)=\Gamma(\operatorname{id}_{\mathbb Z}\times f)\). Thus
	\(\Sigma\) defines a functor
	\(\Sigma:\ell\mathcal G\to\mathcal{PMV}\).
	
	\begin{theorem}[Di Nola--Lettieri equivalence,
		{\cite{mundici2000,dinola1994}}]
		\label{thm:dinola-lettieri}
		The functor
		\[
		\Sigma:\ell\mathcal G\longrightarrow\mathcal{PMV}
		\]
		is an equivalence of categories.
	\end{theorem}
	
	In particular, there exists a quasi-inverse functor
	\(\Delta:\mathcal{PMV}\to\ell\mathcal G\). We denote the corresponding
	natural isomorphisms by
	\[
	\kappa:\operatorname{id}_{\mathcal{PMV}}\longrightarrow\Sigma\Delta,
	\qquad
	\zeta:\Delta\Sigma\longrightarrow\operatorname{id}_{\ell\mathcal G}.
	\]
	Their naturality is expressed by the following commutative diagrams:
	\begin{equation}\label{diag:naturality-kappa-zeta}
		\begin{tikzcd}[row sep=large, column sep=large]
			A \arrow[r, "\kappa_A"] \arrow[d, "h"']
			& \Sigma\Delta(A) \arrow[d, "\Sigma\Delta(h)"] \\
			B \arrow[r, "\kappa_B"']
			& \Sigma\Delta(B)
		\end{tikzcd}
		\qquad\qquad
		\begin{tikzcd}[row sep=large, column sep=large]
			\Delta\Sigma(G) \arrow[r, "\zeta_G"] \arrow[d, "\Delta\Sigma(f)"']
			& G \arrow[d, "f"] \\
			\Delta\Sigma(H) \arrow[r, "\zeta_H"']
			& H
		\end{tikzcd}
	\end{equation}
	for every morphism $h:A\to B$ in $\mathcal{PMV}$ and every morphism $f:G\to H$ in $\ell\mathcal{G}$.
	
	Their components satisfy the triangle identities
	\begin{equation}\label{eq:DL-triangle-identities}
		\Sigma(\zeta_H)\circ\kappa_{\Sigma(H)}
		=\operatorname{id}_{\Sigma(H)},
		\qquad
		\zeta_{\Delta(A)}\circ\Delta(\kappa_A)
		=\operatorname{id}_{\Delta(A)}.
	\end{equation}
	
	We now recall some basic notions and results from topology; for a comprehensive treatment we refer to \cite{aussenhofer2022,arhangelskii2008,munkres2000}.
	
	If $(X,\tau)$ is a topological space and $Y \subseteq X$, we denote by $\tau|_Y$ the subspace topology on $Y$ induced by $\tau$, i.e., the topology whose open sets are $U \cap Y$ with $U \in \tau$.
	
	A family $\mathcal{B}(x)$ of neighbourhoods of a point $x$ is called a \emph{neighbourhood basis} if every neighbourhood of $x$ contains some member of $\mathcal{B}(x)$. The members of $\mathcal{B}(x)$ need not be open.
	
	A \emph{filter basis} on a set $X$ is a nonempty family $\mathcal{F}$ of nonempty subsets such that for any $A, B \in \mathcal{F}$, there exists $C \in \mathcal{F}$ with $C \subseteq A \cap B$. A \emph{filter basis at $x\in X$} is a filter basis all of whose members contain $x$. In particular, a neighbourhood basis of a point is a filter basis at that point.
	
	We adopt the convention that the separation properties Hausdorff ($T_2$), regular ($T_3$), and completely regular ($T_{3\frac{1}{2}}$)  are defined to include the $T_1$ condition. Thus the axioms become progressively stronger: $\text{completely regular}\Rightarrow \text{regular}  \Rightarrow  \text{Hausdorff}  \Rightarrow  T_1 \Rightarrow T_0 $.
	
	A topological space is called \emph{first-countable} if every point has a countable neighbourhood basis, \emph{second-countable} if it has a countable basis for its topology,  \emph{separable} if it contains a countable dense subset,  and \emph{metrizable} if its topology is induced by a metric. Every metrizable space is Hausdorff and first-countable; a separable metric space is second-countable.
	
	A \emph{topological group} is a group $(G, \cdot)$ equipped with a topology such that the group multiplication and the inversion map are continuous.

	A \emph{group topology} on $G$ is a topology making $G$ into a
	topological group.
	
	In this paper, a \emph{lattice-ordered topological group} is an
	$\ell$-group $G$ equipped with a group topology $\tau$, without any
	assumption that the lattice operations are continuous. We use this
	terminology to distinguish such groups from topological $\ell$-groups
	in the standard sense.
	
	A lattice-ordered topological group $(G,\tau)$ is called a
	\emph{topological $\ell$-group} if the lattice operations
	\[
	(x,y)\longmapsto x\vee y
	\qquad\text{and}\qquad
	(x,y)\longmapsto x\wedge y
	\]
	are continuous; see, for example, the classical works
	\cite{ball1979,madell1969,redfield1974,ball1979} and the more recent studies
	\cite{jordan2023,ranjbar2020}.
	
	Topological groups are homogeneous: for any $a \in G$, the left translation $x \mapsto a x$ is a homeomorphism. Consequently, the topological properties of a topological group can be studied by examining the neighbourhoods of the identity element.
	
	We will need the following classical facts about topological groups.
	
	\begin{proposition}[{\cite[Theorem 5.2.14]{aussenhofer2022}}]\label{prop:topological-group-T0-Tychonoff}
		A $T_0$ topological group is completely regular.
	\end{proposition}
	
	\begin{theorem}[Birkhoff–Kakutani, \cite{birkhoff1936, kakutani1936}; see also \cite{aussenhofer2022, kechris1995}]\label{thm:birkhoff-kakutani}
		Let $G$ be a topological group. Then $G$ is metrizable if and only if $G$ is Hausdorff and first-countable. Moreover, if $G$ is metrizable, it admits a compatible left-invariant metric $d$, i.e., $d(xy, xz) = d(y, z)$ for all $x, y, z \in G$.
	\end{theorem}
	
	We now recall some basic notions from category theory that will be used throughout the paper. In particular, we focus on adjoint functors.
	
	A pair of functors $F: \mathcal{C} \to \mathcal{D}$ and $G: \mathcal{D} \to \mathcal{C}$ is called an \emph{adjoint pair}, denoted $F \dashv G$, if there exists a natural bijection
	\[
	\operatorname{Hom}_{\mathcal{D}}(F(X), Y) \cong \operatorname{Hom}_{\mathcal{C}}(X, G(Y))
	\]
	for all $X \in \mathcal{C}$, $Y \in \mathcal{D}$. Equivalently, there exist natural transformations $\eta: \mathrm{id}_{\mathcal{C}} \to GF$ (the unit) and $\varepsilon: FG \to \mathrm{id}_{\mathcal{D}}$ (the counit) satisfying the triangle identities. To emphasize the direction, we sometimes write $F \dashv G: \mathcal{C} \rightharpoonup \mathcal{D}$.
	
	The following characterization will be useful in our proofs.
	
	\begin{proposition}[\cite{riehl2016}]\label{prop:adj-characterization}
		Let $F: \mathcal{C} \to \mathcal{D}$ and $G: \mathcal{D} \to \mathcal{C}$ be functors. Then $F$ is left adjoint to $G$ if and only if there exists a natural transformation $\eta: \mathrm{id}_{\mathcal{C}} \to GF$ such that for every $X \in \mathcal{C}$ and $Y \in \mathcal{D}$, the map
		\[
		\operatorname{Hom}_{\mathcal{D}}(FX, Y) \to \operatorname{Hom}_{\mathcal{C}}(X, GY), \quad f \mapsto Gf \circ \eta_X
		\]
		is a bijection.
	\end{proposition}
	
	We conclude this preliminary section by fixing some notation that will be used repeatedly in the sequel.
	
	For a subset $U$ of an abelian $\ell$-group $G$, set
	\begin{equation}\label{eq:U_triangle}
		U^{\vartriangle} := \{ x \in G \mid x^+ \in U \text{ and } x^- \in U \}.
	\end{equation}
	For a sequence $\{U_i\}_{i=1}^\infty$ of subsets of $G$, define
	\begin{equation}\label{eq:sum_Ui}
		\sum_{i=1}^{\infty} U_i
		:=
		\left\{
		\sum_{k=1}^{n}u_{i_k}
		\;\middle|\;
		n\in\mathbb N,\;
		i_1<i_2<\cdots<i_n,\;
		u_{i_k}\in U_{i_k}
		\right\}.
	\end{equation}
	
	For any \(\ell\)-group homomorphism \(f:G\to H\), we have
	\[
	f(U^{\vartriangle}) \subseteq (f(U))^{\vartriangle}
	\quad \text{and} \quad
	f\!\left(\sum_{i=1}^{\infty} U_i\right)
	=
	\sum_{i=1}^{\infty} f(U_i),
	\]
	where the notation on the right-hand side is understood in \(H\). Indeed, the first
	inclusion follows from
	\[
	f(x^+)=f(x)^+,\qquad f(x^-)=f(x)^-,
	\]
	and the second equality follows directly from the definition of
	\(\sum_{i=1}^{\infty}U_i\) as the set of all finite sums. If \(f\) is an
	\(\ell\)-group isomorphism, then
	\[
	f(U^{\vartriangle})=(f(U))^{\vartriangle}.
	\]

\section{Neighbourhood basis characterisation}\label{sec:0-neighbor}

The  absence of a
systematic description of neighbourhood bases at $0$ in a topological
MV-algebra is the first obstacle to resolving the topological lifting
problem. In this section, we provide a complete characterisation of
which filter bases can serve as a neighbourhood basis of $0$, thereby
 laying the foundation for the notion of summability
introduced in Section~\ref{sec:summable}.

\begin{definition}[\cite{hoo1997, najafi2017, weber2012}]\label{def:top-MV}
	A \emph{topological MV-algebra} is a pair $(A,\tau)$ consisting of an MV-algebra $(A, \oplus, \neg, 0)$ and a topology $\tau$ on $A$ such that the maps $(x,y) \mapsto x \oplus y$ and $x \mapsto \neg x$ are continuous.
\end{definition}

\begin{remark}[Filter topologies]
A prominent family of examples is given by \emph{filter topologies} on
MV-algebras (introduced by Luan and Yang~\cite{luanyang2017}): every
MV-algebra admits a filter topology, though it need not automatically be
a topological MV-algebra---the precise conditions under which it is were
later characterised by Luan, Weber and Yang~\cite{luanweberyang2021},
who also proved that filter topologies are always locally convex
(whereas a general topological MV-algebra need not be, see
e.g.\ Example~\ref{ex:non-locally-convex}) and that every MV-chain with
its order (interval) topology is a topological MV-algebra, since by
\cite[Cor.~4.5]{luanweberyang2021} the order topology on an MV-chain
coincides with a filter topology and hence falls under
\cite[Thm.~4.7]{luanweberyang2021}. The filter topology of an MV-chain
was further studied by Wu, Luan and Yang~\cite{wuluanyang2020}, and the
finite and simple cases by Asadzadeh, Rezaei and
Jamalzadeh~\cite{asadzadeh2020}. The distinction between ``locally
convex'' and the weaker neighbourhood condition of \emph{summability}
introduced in Definition~\ref{def:sumMV}  is what allows the $\ell$-group lifting to be captured
exactly.
\end{remark}

Since $\odot$, $\vee$, $\wedge$, and $\ominus$ can be expressed in terms of $\oplus$ and $\neg$, their continuity follows from that of $\oplus$ and $\neg$.

\begin{lemma}[\cite{hoo1997, najafi2017}]\label{lem:continuous-ops}
	In a topological MV-algebra $(A,\tau)$, the operations $\odot$, $\vee$, $\wedge$, and $\ominus$ are continuous.
\end{lemma}

For an MV-algebra $A$, a subset $U \subseteq A$ and an element $a \in A$, we define
\begin{equation}\label{eq:Ua}
	U(a) := \{ x \in A \mid x \ominus a \in U \text{ and } a \ominus x \in U \}.
\end{equation}
The following properties are straightforward: if $U \subseteq V$ then $U(a) \subseteq V(a)$, $(U \cap V)(a) = U(a) \cap V(a)$, $\neg U(a) = U(\neg a)$, and if $0 \in U$ then $U(1) = \neg U$ and $U(0) = U$.

The relevance of this notation to topological MV-algebras is that a neighbourhood basis at $0$ gives rise to neighbourhood bases at all points:

\begin{proposition}[\cite{gan2024, weber2012}]\label{prop:neighbourhood-base}
	Let $(A, \tau)$ be a topological MV-algebra and let $\mathcal{V}(0)$ be a neighbourhood basis of $0$. Then for any $a \in A$, the family
	$\{ U(a) \mid U \in \mathcal{V}(0) \}$ (with $U(a)$ as defined in \eqref{eq:Ua})
	forms a neighbourhood basis of $a$.
\end{proposition}

The next proposition shows that the topology is in fact determined by the neighbourhoods of any point:
\begin{proposition}[{\cite[Theorem 3.6]{gan2024}}]\label{prop:topology-determined}
	For any $a \in A$, the topology $\tau$ of a topological MV-algebra $(A,\tau)$ is uniquely determined by its $a$-neighbourhood system.
\end{proposition}

This shows that, just as in the case of topological groups, the entire topology of a topological MV-algebra is encoded in the neighbourhoods of $0$. For topological groups, a complete characterization of which filter bases at $0$ can serve as a neighbourhood basis for some group topology is a classical result; we recall it here for comparison.

\begin{theorem}[\cite{aussenhofer2022, pontryagin1946}]\label{thm:group-topology}
	 Let $(G, +, 0)$ be a group. If $\mathcal{V}(0)$ is a filter basis that forms a neighbourhood basis of $0$ for some group topology on $G$, then $\mathcal{V}(0)$ satisfies:
	\begin{enumerate}
		\item[(a)] For every $U \in \mathcal{V}(0)$, there exists $V \in \mathcal{V}(0)$ such that $V + V \subseteq U$.
		\item[(b)] For every $U \in \mathcal{V}(0)$, there exists $V \in \mathcal{V}(0)$ such that $-V \subseteq U$.
		\item[(c)] For every $U \in \mathcal{V}(0)$ and every $a \in G$, there exists $V \in \mathcal{V}(0)$ such that $a + V - a \subseteq U$.
	\end{enumerate}
	Conversely, if $\mathcal{V}$ is a filter basis on $G$ satisfying (a), (b) and (c), then there exists a unique group topology $\tau$ on $G$ such that $\mathcal{V}$ is a neighbourhood basis of $0$ in $\tau$.
\end{theorem}

Those familiar with the theory of topological groups will recognise this as the inverse problem foreshadowed earlier: which families of subsets can serve as a neighbourhood basis of $0$ in some topological MV-algebra?

One might try to imitate Theorem \ref{thm:group-topology} by simply replacing $+$ with $\oplus$ and $-$ with $\neg$. That is, one could consider the following conditions on a filter basis $\mathcal{V}(0)$ at $0$:
\[
\forall U \in \mathcal{V}(0),\; \exists V \in \mathcal{V}(0) \text{ such that } V \oplus V \subseteq U \text{ and } \neg V \subseteq U(1),
\]
where $U(1)$ is defined as in \eqref{eq:Ua}. (Note that condition (c) in Theorem \ref{thm:group-topology} is redundant here because $\oplus$ is commutative.)
However, just as Example \ref{ex:counter1} below illustrates, a filter basis satisfying these conditions need not arise as a neighbourhood basis of $0$ for any topological MV-algebra.

\begin{example}\label{ex:counter1}
	Let $A = L_3 = \left\{0, \frac12, 1\right\}$. Define a filter basis
	\[
	\mathcal{V} = \bigl\{ \{0,1\},\; A \bigr\}.
	\]
	This filter basis satisfies: for every $U \in \mathcal{V}$, there exists $V \in \mathcal{V}$ such that $V \oplus V \subseteq U$, $\neg V \subseteq \neg U = U(1)$, and even $V \ominus V \subseteq U$.
	
	If there existed a topological MV-algebra $(A,\tau)$ having $\mathcal{V}$ as a neighbourhood basis of $0$, then by Proposition \ref{prop:neighbourhood-base}, the family $\{V(a) \mid V \in \mathcal{V}\}$ would be a neighbourhood basis of $a$ in $\tau$. However, take $a=0$, $b=\frac12$ and $U = \{0,1\}$. Then $a \oplus b = \frac12$, while $U\!\left(\frac12\right) = \left\{\frac12\right\}$. This means that $\oplus$ is not continuous at $(0,\frac12)$, because there is no $V \in \mathcal{V}$ such that $V(0) \oplus V\!\left(\frac12\right) \subseteq U\!\left(\frac12\right)$.
	
	Hence no topological MV-algebra admits $\mathcal{V}$ as a neighbourhood basis of $0$.

The deeper reason for this failure is that topological MV-algebras lack the homogeneity
property of topological groups: in a topological group, translating a neighbourhood basis
of the identity yields a neighbourhood basis at any point, whereas in a topological
MV-algebra the neighbourhood system of a point $a$ is given by $\{U(a) \mid U \in
\mathcal{V}(0)\}$ (Proposition~\ref{prop:neighbourhood-base}), which is not simply a translate
of $\mathcal{V}(0)$.
\end{example}

By leveraging the subtler algebraic properties of MV-algebras, and despite a proof more involved than in the group case, we obtain

\begin{proposition}\label{thm:MV-topology}
	Let $(A, \oplus, \neg, 0)$ be an MV-algebra. If $\mathcal{V}(0)$ is a filter basis that forms a neighbourhood basis of $0$ for some topology on $A$ that makes $\oplus$ and $\neg$ continuous, then $\mathcal{V}(0)$ satisfies:
	\begin{enumerate}
		\item[(i)] For every $U \in \mathcal{V}(0)$ and every $a, b \in A$, there exists $V \in \mathcal{V}(0)$ such that $V(a) \oplus V(b) \subseteq U(a \oplus b)$.
	\end{enumerate}
	Conversely, if $\mathcal{V}$ is a filter basis at $0$ satisfying (i), then there exists a unique topological MV-algebra $(A,\tau)$ such that $\mathcal{V}$ is a neighbourhood basis of $0$ in $\tau$.
\end{proposition}

\begin{proof}
	Let $(A,\tau)$ be a topological MV-algebra and let $\mathcal{V}(0)$ be a neighbourhood basis of $0$ in $\tau$.
	
	By Proposition \ref{prop:neighbourhood-base}, for any $a \in A$, the family $\{U(a) \mid U \in \mathcal{V}(0)\}$ forms a neighbourhood basis of $a$ (where $U(a)$ is defined in \eqref{eq:Ua}).
	By continuity of $\oplus$, for every $U \in \mathcal{V}(0)$ and every $a, b \in A$, there exists $V \in \mathcal{V}(0)$ such that $V(a) \oplus V(b) \subseteq U(a \oplus b)$. Thus condition (i) holds.
	
	Conversely, suppose $\mathcal{V}$ is a filter basis at \(0\) satisfying (i). Taking $a = 0$ and $b = 1$ in (i), we obtain that for every $U \in \mathcal{V}$, there exists $V \in \mathcal{V}$ such that $V(0) \oplus V(1) \subseteq U(1)$. Applying the negation map to both sides yields $\neg V(0) \odot \neg V(1) \subseteq \neg U(1)$. Using $\neg V(a) = V(\neg a)$, we get $V \ominus V \subseteq U$. Hence $\mathcal{V}$ also satisfies:
	\begin{enumerate}
		\item[($*$)] For every $U \in \mathcal{V}$, there exists $V \in \mathcal{V}$ such that $V \ominus V \subseteq U$.
	\end{enumerate}
	
	Now define a family $\tau$ of subsets of $A$ as follows:
	\[
	\tau := \{ O \subseteq A \mid \forall a \in O,\ \exists U \in \mathcal{V} \text{ such that } U(a) \subseteq O \}.
	\]
	
	We first verify that $\tau$ is a topology on $A$.
	
	Clearly $\emptyset, A \in \tau$, and $\tau$ is closed under arbitrary unions.
	
	Assume $O_1, O_2 \in \tau$. If $O_1 \cap O_2 = \emptyset$, then $O_1 \cap O_2 \in \tau$ trivially. Otherwise, take any $a \in O_1 \cap O_2$. By definition, there exist $U_1, U_2 \in \mathcal{V}$ such that $U_1(a) \subseteq O_1$ and $U_2(a) \subseteq O_2$. Since $\mathcal{V}$ is a filter basis, we can find $U \in \mathcal{V}$ with $U \subseteq U_1 \cap U_2$. Then $U(a) \subseteq U_1(a) \cap U_2(a) \subseteq O_1 \cap O_2$, proving that $O_1 \cap O_2 \in \tau$.
	
	We now prove that for each $a \in A$, the filter basis $\mathcal{V}(a) := \{ U(a) \mid U \in \mathcal{V} \}$ indeed forms a neighbourhood basis of $a$ in the topology $\tau$. In particular, for $a = 0$, $\mathcal{V}$ itself is a neighbourhood basis of $0$.
	
	Let $\mathcal{V}_\tau(a)$ denote the neighbourhood filter of $a$ in $\tau$. We prove that $\mathcal{V}(a)$ is a basis of $\mathcal{V}_\tau(a)$.
	
	First, take any $V \in \mathcal{V}_\tau(a)$. Then there exists an open set $O \in \tau$ such that $a \in O \subseteq V$. By definition of $\tau$, we can find $U \in \mathcal{V}$ with $U(a) \subseteq O \subseteq V$. Since $U(a) \in \mathcal{V}(a)$, this shows that every neighbourhood of $a$ contains some member of $\mathcal{V}(a)$.
	
	Conversely, take any $U(a) \in \mathcal{V}(a)$. Define
	\[
	O := \{ b \in U(a) \mid \exists W \in \mathcal{V},\ W(b) \subseteq U(a) \}.
	\]
	We claim that $O$ is open. If so, then $O \in \mathcal{V}_\tau(a)$ and $a \in O \subseteq U(a)$, which shows that every member of $\mathcal{V}(a)$ contains a neighbourhood of $a$, completing the proof that $\mathcal{V}(a)$ is a basis of $\mathcal{V}_\tau(a)$.
	
	To see that $O$ is open, take any $x \in O$. By definition, there exists $W \in \mathcal{V}$ such that $W(x) \subseteq U(a)$. By \((*)\), choose \(T\in\mathcal V\) such that
	\(T\ominus T\subseteq W\). Applying (i) with \(a=b=0\), and using
	\(S(0)=S\) and \(T(0)=T\), choose \(S\in\mathcal V\) such that
	\(S\oplus S\subseteq T\). Hence
	\[
	(S\oplus S)\ominus(S\oplus S)
	\subseteq T\ominus T
	\subseteq W.
	\]
	
	We assert that $S(x) \subseteq O$. If this holds, then $x$ is an interior point of $O$, and since $x$ was arbitrary, $O$ is open.
	
	Take any $c \in S(x)$. We need to show $c \in O$, i.e., find $S' \in \mathcal{V}$ with $S'(c) \subseteq U(a)$. We claim that $S$ itself works, i.e., $S(c) \subseteq U(a)$. For any $y \in S(c)$, we have $c \ominus y \in S$ and $y \ominus c \in S$. Since $c \in S(x)$, we also have $c \ominus x \in S$ and $x \ominus c \in S$.
	
	Now use the MV-algebraic identity from Lemma \ref{lem:MV-identity}:
	\[
	x \ominus y = \bigl( (x \ominus c) \oplus (c \ominus y) \bigr) \ominus \bigl( (c \ominus x) \oplus (y \ominus c) \bigr).
	\]
	All four terms $(x \ominus c), (c \ominus y), (c \ominus x), (y \ominus c)$ belong to $S$, hence
	\[
	x \ominus y \in (S \oplus S) \ominus (S \oplus S) \subseteq W.
	\]
	Similarly, $y \ominus x \in W$. Therefore $y \in W(x) \subseteq U(a)$. This proves $S(c) \subseteq U(a)$, and consequently $c \in O$.
	
	Thus $\mathcal{V}(a)$ is indeed a neighbourhood basis of $a$ in $\tau$.
	
	Finally, condition (i) together with the fact that $\mathcal{V}(a)$ is a neighbourhood basis of $a$ implies continuity of $\oplus$. Continuity of $\neg$ follows from $\neg(U(x)) = U(\neg x)$. Hence $(A,\tau)$ is a topological MV-algebra.
	
	Uniqueness of the topology follows from Proposition \ref{prop:topology-determined}.
\end{proof}

\medskip
\noindent Note that examples satisfying Proposition~\ref{thm:MV-topology} are plentiful;
they include totally ordered, perfect, ideal-induced, finite, and product
algebras, all of which will be treated systematically in
Section~\ref{sec:summable}.

\medskip
\noindent As a further illustration of the scope of Proposition~\ref{thm:MV-topology},
topological MV-algebras are not homogeneous in general.    
The following table places Proposition~\ref{thm:MV-topology} in the broader context of topological structures.
 Despite being a subclass of topological monoids, topological MV-algebras exhibit properties that are strikingly similar to those of topological groups. Here ``Top.'' abbreviates ``topological'', and ``nbhds'' stands for ``neighbourhoods''. 
The facts for
topological monoids listed in the table can be found in \cite{arhangelskii2008, banakh2017}.
\begin{center}\label{tab:contex}
\small 
\begin{tabular} {l|c|c|c}
	\hline
	\textbf{Feature} & \textbf{Top. group} & \textbf{Top. MV-alg} & \textbf{Top. monoid} \\
	\hline
	Separation & $T_0 \Rightarrow T_3$ & $T_0 \Rightarrow T_3$ & $T_0 \not\Rightarrow T_1$ \\
	Homogeneity & yes & no & no \\
	$0$-nbhds $\Rightarrow$ whole topology &  yes  & yes & no \\
	$a$-nbhds $\Rightarrow$ whole topology &  yes  & yes & no \\
	$0$-nbhd basis characterization &  yes  & yes (Proposition \ref{thm:MV-topology}) & not uniquely determined \\
	\hline
\end{tabular}
\end{center}

\medskip
If we require that the filter basis in Proposition \ref{thm:MV-topology} be an open neighbourhood basis of $0$ (i.e., that every $U \in \mathcal{V}$ be open), then we need to add the following condition (ii):

\begin{proposition}\label{thm:MV-topology2}
	Let $(A, \oplus, \neg, 0)$ be an MV-algebra. If $\mathcal{V}(0)$ is a filter basis that forms an open neighbourhood basis of $0$ for some topology on $A$ that makes $\oplus$ and $\neg$ continuous, then $\mathcal{V}(0)$ satisfies:
	\begin{enumerate}
		\item[(i)] For every $U \in \mathcal{V}(0)$ and every $a, b \in A$, there exists $V \in \mathcal{V}(0)$ such that $V(a) \oplus V(b) \subseteq U(a \oplus b)$.
		\item[(ii)] For every $U \in \mathcal{V}(0)$ and every $x \in U$, there exists $V \in \mathcal{V}$ such that $V(x) \subseteq U$.
	\end{enumerate}
	Conversely, if $\mathcal{V}$ is a filter basis at $0$ satisfying (i) and (ii), then there exists a unique topological MV-algebra $(A,\tau)$ such that $\mathcal{V}$ is an open neighbourhood basis of $0$ in $\tau$.
\end{proposition}

\begin{proof}
	By Proposition \ref{thm:MV-topology}, it suffices to verify the part concerning condition (ii).
	
	If $(A,\tau)$ is a topological MV-algebra and $\mathcal{V}$ is an open neighbourhood basis of $0$ in $\tau$, then since each $U \in \mathcal{V}$ is open, $\mathcal{V}$ automatically satisfies (ii).
	
	Conversely, assume $\mathcal{V}$ is a filter basis at \(0\) satisfying (i) and (ii). We need to show that every member of $\mathcal{V}$ is open in the topology $\tau$ constructed in the proof of Proposition \ref{thm:MV-topology}. That is, for every $U \in \mathcal{V}$ and every $x \in U$, there exists $V \in \mathcal{V}$ such that $V(x) \subseteq U$. This is exactly condition (ii), so the conclusion follows.
\end{proof}

\begin{remark}
	A sufficient condition for a family $\mathcal{V}$ of subsets of $A$ with $0$ as a common point to be a filter basis and to satisfy (ii) is the following:
	\begin{quote}
		(iii) For every $U, V \in \mathcal{V}$ and $a, b \in A$, if $U(a) \cap V(b) \neq \emptyset$, then for each $y \in U(a) \cap V(b)$ there exists $W \in \mathcal{V}$ such that $W(y) \subseteq U(a) \cap V(b)$.
	\end{quote}
	
	\cite[Proposition 3.5]{gan2024} states that if a family $\mathcal{V}$ of subsets of $A$ has $0$ as a common point and satisfies conditions (i) and (iii), then $\{U(a) \mid U \in \mathcal{V}, a \in A\}$ forms a basis for a topological MV-algebra.
	
	Proposition \ref{thm:MV-topology2}, on the other hand, starts from a filter basis $\mathcal{V}$ at $0$ satisfying only (i) and (ii) and yields a unique topological MV-algebra for which $\mathcal{V}(a) = \{U(a) \mid U \in \mathcal{V}\}$ is an open neighbourhood basis of $a$. Consequently, $\{U(a) \mid U \in \mathcal{V}, a \in A\}$ forms a basis for this topology.
	
	Since an open filter basis satisfying (ii) automatically satisfies (iii), Proposition \ref{thm:MV-topology2} strictly generalizes the setting of \cite[Proposition 3.5]{gan2024}, whence the latter is a corollary rather than an equivalent.
\end{remark}

\medskip
Taken together, Propositions~\ref{thm:MV-topology}, 
\ref{prop:neighbourhood-base}, \ref{prop:topology-determined},
\ref{thm:MV-topology2} establish a foundational principle: the
neighbourhood system of $0$ not only determines the entire topology of a
topological MV-algebra, but also admits a precise characterisation of which
filter bases can bear this responsibility. These results provide the
essential ingredient for the notion of \emph{summability}, to be introduced
in Section~\ref{sec:summable}.

\section{Summable topological MV-algebras}\label{sec:summable}

Having characterised in Section~\ref{sec:0-neighbor} which filter bases at $0$
determine a topological MV-algebra, we now return to the long-standing gap
posed in the Introduction: can every topological MV-algebra be embedded into
a lattice-ordered topological group? The key to closing this gap is a purely
neighbourhood-theoretic condition on the $0$-neighbourhood basis, which we
call \emph{summability}. We will prove that summability characterises
precisely when a topological MV-algebra can be realised as the unit interval
of a lattice-ordered topological group, thereby lifting both  Mundici
equivalence and  Di Nola--Lettieri equivalence to the topological
setting. To this end, we first establish a sufficient condition.

\begin{proposition}\label{prop:starcon}
	Let $(A, \tau)$ be a topological MV-algebra and let $\mathcal{V}(0)$ be an open neighbourhood basis of $0$. Let $\Xi(A) = (G_A, u_{G_A})$ be the Mundici group of $A$, and let $\varphi_A: A \to \Gamma(G_A, u_{G_A})$ be the natural isomorphism given by Mundici equivalence (see Theorem \ref{thm:mundici}).
	Consider the following condition, which we call the \emph{summable condition} (SC for short):
	\begin{equation*}\label{cond-star}
		\begin{gathered}
			\forall U \in \mathcal{V}(0),\ \forall a \in A,\ \exists \{V_i\}_{i=1}^{\infty} \subseteq \mathcal{V}(0) \text{ such that } \\
			[0, u_{G_A}] \cap \left( \varphi_A(a) + \sum_{i=1}^{\infty} (\varphi_A(V_i))^{\vartriangle} \right) \subseteq \varphi_A(U(a)),
		\end{gathered}
		\tag{SC}
	\end{equation*}
	where $(\varphi_A(V_i))^{\vartriangle}$ is defined as in \eqref{eq:U_triangle} and the sum $\sum_{i=1}^{\infty}$ is defined in \eqref{eq:sum_Ui}.
	
	If \(\mathcal V(0)\) satisfies the summable condition \eqref{cond-star}, then there exists a group topology $\tau_{G_A}$ on $G_A$ such that $\varphi_A: (A, \tau) \to ([0, u_{G_A}], \tau_{G_A}|_{[0, u_{G_A}]})$ is a homeomorphism.
\end{proposition}

\begin{proof}
	We first construct a group topology on $G_A$ from the topological MV-algebra $(A, \tau)$. By Theorem \ref{thm:group-topology}, it suffices to construct a filter basis $\mathcal{V}$ on $G_A$ satisfying conditions (a), (b) and (c).
	
	Set
	\[
		\mathcal{V} := \left\{ \sum_{i=1}^{\infty} (\varphi_A(U_i))^{\vartriangle} \;\middle|\; \{U_i\}_{i=1}^{\infty} \subseteq \mathcal{V}(0) \right\},
	\]
	where $(\cdot)^{\vartriangle}$ is defined in \eqref{eq:U_triangle} and the sum $\sum_{i=1}^{\infty}$ is defined in \eqref{eq:sum_Ui}, and $\mathcal{V}(0)$ is an open neighbourhood basis of $0$ in $(A,\tau)$.
	
	We now show that $\mathcal{V}$ is a filter basis on $G_A$ and that it satisfies conditions (a), (b) and (c) of Theorem \ref{thm:group-topology}.
	
	\emph{Step 1.} $\mathcal{V}$ is a filter basis. Take any $W_1 = \sum_{i=1}^{\infty} (\varphi_A(U_i))^{\vartriangle} \in \mathcal{V}$ and $W_2 = \sum_{i=1}^{\infty} (\varphi_A(U'_i))^{\vartriangle} \in \mathcal{V}$, where $\{U_i\}_{i=1}^{\infty}, \{U'_i\}_{i=1}^{\infty} \subseteq \mathcal{V}(0)$. Since $\mathcal{V}(0)$ is the neighbourhood basis of $0$, for each $i \in \mathbb{N}$ there exists $V_i \in \mathcal{V}(0)$ such that $0 \in V_i \subseteq U_i \cap U'_i$, and consequently $\varphi_A(V_i) \subseteq \varphi_A(U_i) \cap \varphi_A(U'_i)$. By definition, $\sum_{i=1}^{\infty} (\varphi_A(V_i))^{\vartriangle} \subseteq W_1 \cap W_2$ and $\sum_{i=1}^{\infty} (\varphi_A(V_i))^{\vartriangle} \in \mathcal{V}$. Hence $\mathcal{V}$ is a filter basis of $0$ in $G_A$.
	
	\emph{Step 2.} Verification of condition (a). Take any $W = \sum_{i=1}^{\infty} (\varphi_A(U_i))^{\vartriangle} \in \mathcal{V}$. For each \(i\in\mathbb N\), choose
	\(V_i\in\mathcal V(0)\) such that
	\[
	V_i\subseteq U_{2i-1}\cap U_{2i}.
	\] We claim that
	$
	\sum_{i=1}^{\infty} (\varphi_A(V_i))^{\vartriangle} + \sum_{i=1}^{\infty} (\varphi_A(V_i))^{\vartriangle} \subseteq W.
	$
	
	Take any $x, y \in \sum_{i=1}^{\infty} (\varphi_A(V_i))^{\vartriangle}$. Then $x = x_1 + x_2 + \cdots + x_m$ and $y = y_1 + y_2 + \cdots + y_n$, where $x_k \in (\varphi_A(V_k))^{\vartriangle}$, $y_l \in (\varphi_A(V_l))^{\vartriangle}$, i.e., $x_k^{+}, x_k^{-} \in \varphi_A(V_k)$ and $y_l^{+}, y_l^{-} \in \varphi_A(V_l)$. We may assume $m = n$ (adding zero terms if necessary). Then
	$
	x + y = (x_1 + y_1) + (x_2 + y_2) + \cdots + (x_n + y_n).
	$
	By construction, for each $k$ we have $x_k^{+}, x_k^{-} \in \varphi_A(V_k) \subseteq \varphi_A(U_{2k-1})$, so $x_k \in (\varphi_A(U_{2k-1}))^{\vartriangle}$; similarly $y_k^{+}, y_k^{-} \in \varphi_A(V_k) \subseteq \varphi_A(U_{2k})$, so $y_k \in (\varphi_A(U_{2k}))^{\vartriangle}$. Hence
	$
	x + y \in \sum_{k=1}^{2n} (\varphi_A(U_k))^{\vartriangle} \subseteq \sum_{i=1}^{\infty} (\varphi_A(U_i))^{\vartriangle} = W.
	$
	Thus condition (a) holds.
	
	\emph{Step 3.} Verification of condition (b). Take any $W = \sum_{i=1}^{\infty} (\varphi_A(U_i))^{\vartriangle} \in \mathcal{V}$. We show that $-W = W$, which immediately implies (b).
	
	Let $x \in -W$. Then $-x \in W$, so $-x = x_1 + x_2 + \cdots + x_n$ with $x_k^{+}, x_k^{-} \in \varphi_A(U_k)$. Consequently,
	$
	x = -x_1 - x_2 - \cdots - x_n.
	$
	For each $k$, note that $(-x_k)^{+} = x_k^{-} \in \varphi_A(U_k)$ and $(-x_k)^{-} = x_k^{+} \in \varphi_A(U_k)$, hence $-x_k \in (\varphi_A(U_k))^{\vartriangle}$. Therefore $x \in \sum_{k=1}^{n} (\varphi_A(U_k))^{\vartriangle} \subseteq W$, proving $-W \subseteq W$. The reverse inclusion is similar, so $W \subseteq -W$. Thus $-W = W$, and condition (b) holds.
	
	\emph{Step 4.} Verification of condition (c). Since $G_A$ is abelian, for any $a \in G_A$ and any subset $U \subseteq G_A$, we have $a + U - a = U$. Hence condition (c) is automatically satisfied.
	
	We have shown that $\mathcal{V}$ is a filter basis on $G_A$ satisfying conditions (a), (b) and (c) of Theorem \ref{thm:group-topology}. By that theorem, there exists a unique group topology $\tau_{G_A}$ on $G_A$ for which $\mathcal{V}$ is a filter basis of the neighbourhood system of $0$.
	
	Now take any $W = \sum_{i=1}^{\infty} (\varphi_A(U_i))^{\vartriangle} \in \mathcal{V}$ and any $x \in W$. Write $x = x_1 + x_2 + \cdots + x_n$ with $x_i \in (\varphi_A(U_i))^{\vartriangle}$. Set $W' = \sum_{i=1}^{\infty} (\varphi_A(U_{i+n}))^{\vartriangle} \in \mathcal{V}$. Then $x \in x + W' \subseteq W$, which shows that every $W \in \mathcal{V}$ is open, so $\mathcal{V}$ is actually an open neighbourhood basis of $0$ in $\tau_{G_A}$. Moreover, for any $a \in G_A$, the family $\{ a + W \mid W \in \mathcal{V} \}$ is an open neighbourhood basis of $a$ in $\tau_{G_A}$.
	
	Finally, we verify that $\varphi_A: (A, \tau) \to ([0, u_{G_A}], \tau_{G_A}|_{[0, u_{G_A}]})$ is a homeomorphism.
	
	Since $\varphi_A$ is an MV-algebra isomorphism, it is bijective, and $\varphi_A(\tau)$ induces a topology on $[0, u_{G_A}]$. It remains to show that $\varphi_A(\tau) = \tau_{G_A}|_{[0, u_{G_A}]}$.
	
	On one hand, take any element $U(a)$ from the neighbourhood basis of $a$ in $\tau$. Since $\mathcal V(0)$ satisfies condition (\ref{cond-star}), there exists a sequence $\{V_i\}_{i=1}^{\infty} \subseteq \mathcal{V}(0)$ such that
	\[
	[0, u_{G_A}] \cap \left( \varphi_A(a) + \sum_{i=1}^{\infty} (\varphi_A(V_i))^{\vartriangle} \right) \subseteq \varphi_A(U(a)).
	\]
	The set $[0, u_{G_A}] \cap \left( \varphi_A(a) + \sum_{i=1}^{\infty} (\varphi_A(V_i))^{\vartriangle} \right)$ is a neighbourhood of $\varphi_A(a)$ in $\tau_{G_A}|_{[0, u_{G_A}]}$. By the arbitrariness of $a$, this shows that $\tau_{G_A}|_{[0, u_{G_A}]}$ is finer than $\varphi_A(\tau)$.
	
	On the other hand, take any basic neighbourhood of $\varphi_A(a) \in [0, u_{G_A}]$ in $\tau_{G_A}|_{[0, u_{G_A}]}$, i.e., a set of the form $\left( \varphi_A(a) + \sum_{i=1}^{\infty} (\varphi_A(U_i))^{\vartriangle} \right) \cap [0, u_{G_A}]$. In particular, consider the neighbourhood $\varphi_A(U_1(a))$ of $\varphi_A(a)$ in $\varphi_A(\tau)$. By the definition of $U(a)$ in \eqref{eq:Ua}, we have:
	\begin{align*}
		\varphi_A(x) \in \varphi_A(U_1(a)) & \iff x \in U_1(a) \\
		& \iff x \in A \text{ and } x \ominus a \in U_1 \text{ and } a \ominus x \in U_1 \\
		& \iff \varphi_A(x) \in [0, u_{G_A}] \text{ and } (\varphi_A(x) - \varphi_A(a))^{+} \in \varphi_A(U_1) \\
		& \qquad \text{ and } (\varphi_A(x) - \varphi_A(a))^{-} \in \varphi_A(U_1) \\
		& \iff \varphi_A(x) \in [0, u_{G_A}] \text{ and } \varphi_A(x) - \varphi_A(a) \in (\varphi_A(U_1))^{\vartriangle} \\
		& \iff \varphi_A(x) \in [0, u_{G_A}] \text{ and } \varphi_A(x) \in \varphi_A(a) + (\varphi_A(U_1))^{\vartriangle}.
	\end{align*}
	Hence $\varphi_A(U_1(a)) = \left( \varphi_A(a) + (\varphi_A(U_1))^{\vartriangle} \right) \cap [0, u_{G_A}]$.
	Since $(\varphi_A(U_1))^{\vartriangle} \subseteq \sum_{i=1}^{\infty} (\varphi_A(U_i))^{\vartriangle}$, we obtain
	\[
	\varphi_A(U_1(a)) \subseteq \left( \varphi_A(a) + \sum_{i=1}^{\infty} (\varphi_A(U_i))^{\vartriangle} \right) \cap [0, u_{G_A}].
	\]
	This shows that $\varphi_A(\tau)$ is finer than $\tau_{G_A}|_{[0, u_{G_A}]}$. Therefore the two topologies coincide on $[0, u_{G_A}]$, i.e., $\varphi_A(\tau) = \tau_{G_A}|_{[0, u_{G_A}]}$.
	
	This completes the proof of the proposition.
\end{proof}

The condition \eqref{cond-star} introduced in Proposition~\ref{prop:starcon} is formulated
with respect to a chosen open neighbourhood basis $\mathcal{V}(0)$ of $0$.
However, it is independent of this choice: any two open neighbourhood bases
at $0$ are cofinal in the neighbourhood system of $0$, so the condition
holds for one basis if and only if it holds for any other. Consequently, it
is legitimate to speak of a topological MV-algebra $(A,\tau)$ itself being
summable.

\begin{definition}[Summable topological MV-algebra]\label{def:sumMV}
A topological MV-algebra $(A,\tau)$ is called \emph{summable} if it
satisfies condition \eqref{cond-star}.
\end{definition}

Summability captures the idea that neighbourhoods of $0$ can be ``added''
in a controlled way: if $U$ is a neighbourhood of $0$, then there exists a
smaller neighbourhood $V$ such that $V\oplus V\subseteq U$. This property
is reminiscent of the definition of a topological group, where the
addition map is jointly continuous at $(0,0)$. In fact, for a topological
MV-algebra arising from a lattice-ordered topological group, summability
is automatically satisfied; conversely, we will prove that summability is
sufficient to reconstruct the group topology. The term ``summable'' is
chosen by analogy with abelian topological groups: in a topological group,
every $0$-neighbourhood contains all finite sums formed from a suitably
chosen sequence of $0$-neighbourhoods. Condition~\eqref{cond-star} adapts
this idea to topological MV-algebras by replacing the translated group
neighbourhood $a+U$ with the MV-algebraic neighbourhood $U(a)$.
Importantly, although \eqref{cond-star} is formulated using the canonical
Mundici representation, it is an intrinsic condition on the topology
$\tau$ alone, involving no choice of an ambient group topology; the same
notion of summability, when restricted to perfect topological MV-algebras,
provides the lifting condition for the Di Nola--Lettieri equivalence (see
Theorem~\ref{thm:perfect-topological-equivalence}).

For a summable topological MV-algebra $(A,\tau)$, there may be more
than one group topology on $G_A$ whose restriction to $[0,u_{G_A}]$
coincides with $\varphi_A(\tau)$. For instance, consider the standard
MV-algebra $[0,1]$ with the discrete topology $\tau$. Its Mundici
group is isomorphic to $(\mathbb{R},1)$. Both the discrete topology on
$\mathbb{R}$ and the group topology having $\{2\mathbb{Z}\}$ as a
neighbourhood basis of $0$ restrict to the discrete topology on
$[0,1]$.

In what follows, for a summable topological MV-algebra $(A,\tau)$, we denote by $\tau_{G_A}$ the specific group topology on $G_A$ constructed in the proof of Proposition \ref{prop:starcon}, i.e., the topology for which the family
\begin{equation}\label{eq:V_GA}
	\mathcal{V}_{G_A} := \left\{ \sum_{i=1}^{\infty} (\varphi_A(U_i))^{\vartriangle} \;\middle|\; \{U_i\}_{i=1}^{\infty} \subseteq \mathcal{V}(0) \right\},
\end{equation}
where $(\cdot)^{\vartriangle}$ is defined in \eqref{eq:U_triangle} and the sum $\sum_{i=1}^{\infty}$ is defined in \eqref{eq:sum_Ui}, is an open neighbourhood basis of $0$.
This topology does not depend on the choice of the open neighbourhood basis
\(\mathcal V(0)\).

The following technical lemma will be used later to prove the converse implication of the main theorem of this section.

\begin{lemma}\label{lem:sum-of-triangles}
	Let $G$ be an abelian $\ell$-group, $\tau$ a group topology on $G$, and $\mathcal{V}(0)$ a neighbourhood basis of $0$. Then for every $U \in \mathcal{V}(0)$, there exists a sequence $\{V_i\}_{i=1}^{\infty} \subseteq \mathcal{V}(0)$ such that
	\(
	\sum_{i=1}^{\infty} V_i^{\vartriangle} \subseteq U.
	\)
\end{lemma}

\begin{proof}
	Since $\tau$ is a group topology, there exists $U_1 \in \mathcal{V}(0)$ such that $U_1 \subseteq U_1 + U_1 \subseteq U$.
	For $U_1$, there exists $U_2 \in \mathcal{V}(0)$ such that $U_2 \subseteq U_2 + U_2 \subseteq U_1$.
	Continuing this process, for $U_{n-1}$ we obtain $U_n \in \mathcal{V}(0)$ with $U_n \subseteq U_n + U_n \subseteq U_{n-1}$.
	Thus we obtain a decreasing sequence $\{U_i\}_{i=1}^{\infty} \subseteq \mathcal{V}(0)$.
	
	For each $U_i$, choose $V_i \in \mathcal{V}(0)$ such that $V_i - V_i \subseteq U_i$.
	This yields a sequence $\{V_i\}_{i=1}^{\infty} \subseteq \mathcal{V}(0)$.
	
	Now take any $x \in \sum_{i=1}^{\infty} V_i^{\vartriangle}$. Then $x = x_1 + x_2 + \cdots + x_n$ with $x_k^{+}, x_k^{-} \in V_k$.
	Since $x_k = x_k^{+} - x_k^{-}$ and $x_k^{+}, x_k^{-} \in V_k$, we have $x_k \in V_k - V_k \subseteq U_k$.
	By the construction of $\{U_i\}$, we obtain $x \in U_1 + U_2 + \cdots + U_n \subseteq U$.
	
	Hence $\sum_{i=1}^{\infty} V_i^{\vartriangle} \subseteq U$, as desired.
\end{proof}

With the local structure of $G$ now clarified by Lemma~\ref{lem:sum-of-triangles}, we focus on the specific interval $[0, u]$ to analyze the summability condition for the induced topological MV-algebra.

\begin{proposition}\label{prop:subspacestar}
	Let $G$ be an abelian $\ell$-group, $u \in G$ a positive element, and $\tau$ a group topology on $G$. If $([0,u], \tau|_{[0,u]})$ is a topological MV-algebra, then it is summable.
\end{proposition}

\begin{proof}
	Let
	\[
	G_u:=\{g\in G:|g|\leq nu\text{ for some }n\in\mathbb N\}.
	\]
	Endow \(G_u\) with the subspace topology \(\tau|_{G_u}\).
	Since \(G_u\) is an \(\ell\)-ideal of \(G\), the element \(u\) is a
	strong unit of \(G_u\), and
	\[
	(\tau|_{G_u})|_{[0,u]}=\tau|_{[0,u]}.
	\]
	Thus, replacing \(G\) by \(G_u\) and \(\tau\) by
	\(\tau|_{G_u}\), we may assume that \(u\) is a strong unit of \(G\).
	
	Let $\Xi(\Gamma(G,u)) = (H, u')$ and let $\mathcal{V}$ be a neighbourhood basis of $0$ in $\tau$. Then $\psi_G: G \to H$ is the natural isomorphism of $\ell$-groups preserving the strong unit given by Mundici equivalence, and $\varphi_{[0,u]} = \psi_G|_{[0,u]}$ (see \eqref{eq:psi-varphi}). In this setting, the condition \eqref{cond-star} for $([0,u], \tau|_{[0,u]})$ reads:
	\[
	\begin{gathered}
		\forall U \in \mathcal{V},\ \forall a \in [0,u],\ \exists \{V_i\}_{i=1}^{\infty} \subseteq \mathcal{V} \text{ such that } \\
		\psi_G([0,u]) \cap \left( \psi_G(a) + \sum_{i=1}^{\infty} (\psi_G(V_i\cap [0,u]))^{\vartriangle} \right) \subseteq \psi_G((U \cap [0,u])(a)).
	\end{gathered}
	\]
	Since $\psi_G$ is an isomorphism, this is equivalent to:
	\begin{equation*}
		\forall U \in \mathcal{V},\ \forall a \in [0,u],\ \exists \{V_i\}_{i=1}^{\infty} \subseteq \mathcal{V} \text{ such that }
		[0,u] \cap \left( a + \sum_{i=1}^{\infty} (V_i \cap [0,u])^{\vartriangle} \right) \subseteq (U \cap [0,u])(a).
	\end{equation*}
	Now fix $U \in \mathcal{V}$ and $a \in [0,u]$.
	Since $\{ (a + U') \cap [0, u] \mid U' \in \mathcal{V} \}$ is a neighbourhood basis of $a$ in $\tau|_{[0,u]}$, for the neighbourhood $(U \cap [0, u])(a)$ there exists $U' \in \mathcal{V}$ such that
	\(
	(a + U') \cap [0, u] \subseteq (U \cap [0, u])(a).
	\)
	For this $U'$, by Lemma \ref{lem:sum-of-triangles} there exists a sequence $\{V_i\}_{i=1}^{\infty} \subseteq \mathcal{V}$ such that $\sum_{i=1}^{\infty} V_i^{\vartriangle} \subseteq U'$. Consequently,
	\[
	[0, u] \cap \left( a + \sum_{i=1}^{\infty} (V_i \cap [0, u])^{\vartriangle} \right) \subseteq [0, u] \cap \left( a + \sum_{i=1}^{\infty} V_i^{\vartriangle} \right) \subseteq (a + U') \cap [0, u] \subseteq (U \cap [0, u])(a).
	\]
	This verifies the condition \eqref{cond-star} for $([0,u], \tau|_{[0,u]})$, completing the proof.
\end{proof}

Propositions \ref{prop:starcon} and  \ref{prop:subspacestar} together provide a complete characterization, establishing the first main result of this paper. 

\begin{theorem}\label{thm:iffMV}
	Let $(A, \tau)$ be a topological MV-algebra. Then $(A,\tau)$ is summable if and only if there exist an abelian $\ell$-group $G$, a positive element $u \in G$, and a group topology $\tau_G$ on $G$ such that $(A, \tau)$ is isomorphic (as an MV-algebra) and homeomorphic to $([0, u], \tau_G|_{[0, u]})$.

\end{theorem}

Thus the topological lifting problem for Mundici equivalence is resolved: a topological MV-algebra admits an $\ell$-group topology precisely when it is summable. As foreshadowed in the Introduction, the analogous problem for Di Nola--Lettieri equivalence receives the same answer in the perfect setting, treated in Section~\ref{sec:adjunction}. Before turning to that, we record a simple sufficient condition for summability—one that already suffices to separate summability from local convexity (Corollary \ref{cor-totallyorderedMV}, Example \ref{ex:non-locally-convex}).

\begin{proposition}\label{prop:U-odot-U-zero}
	Let $(A, \tau)$ be a topological MV-algebra and let $\mathcal{V}(0)$ be a neighbourhood basis of $0$. If there exists $U \in \mathcal{V}(0)$ such that
	\begin{equation}\label{cond:odot}
		\forall x \in U,\; x \odot x = 0,
	\end{equation}
	then $(A, \tau)$ is summable. Moreover, the family $\{(\varphi_A(V))^{\vartriangle} \mid V \in \mathcal{V}(0) \}$ forms a neighbourhood basis of $0$ in $\tau_{G_A}$ defined by \eqref{eq:V_GA}.
\end{proposition}

\begin{proof}
	Take any neighbourhood $W$ of $0$ in $\tau$. Since $\ominus$ is continuous at $(0,0)$, there exists $V_0 \in \mathcal{V}(0)$ such that $V_0 \subseteq V_0 \ominus V_0 \subseteq W$.
	Choose $U_0 \in \mathcal{V}(0)$ with $U_0 \subseteq V_0 \cap U$.
	
	Inductively, assume $U_{i-1}$ has been obtained. By continuity of $\oplus$ at $(0,0)$, there exists $V_i \in \mathcal{V}(0)$ such that $V_i \subseteq V_i \oplus V_i \subseteq U_{i-1}$.
	Then choose $U_i \in \mathcal{V}(0)$ such that $U_i \subseteq V_i \cap U$.
	
	Thus we obtain a sequence $\{U_i\}_{i=1}^\infty \subseteq \mathcal{V}(0)$ with $U_i \subseteq U$ for all $i$. Consequently, $x \odot x = 0$ for every $x \in U_i$.
	
	We now verify that for any given $a$,
	\[
	[0, u_{G_A}] \cap \left( \varphi_A(a) + \sum_{i=1}^{\infty} (\varphi_A(U_i))^{\vartriangle} \right) \subseteq \varphi_A(W(a)).
	\]
	
	Take any $y \in [0, u_{G_A}] \cap \left( \varphi_A(a) + \sum_{i=1}^{\infty} (\varphi_A(U_i))^{\vartriangle} \right)$, i.e., $y = \varphi_A(x) = \varphi_A(a) + y_1 +  y_2 + \cdots + y_n$ with $y_i^{+}, y_i^{-} \in \varphi_A(U_i)$ and $x \in A$. Then
	\begin{align*}
		y &= \varphi_A(a) + y_1 + y_2 + \cdots + y_n \\
		&= \varphi_A(a) + (y_1^{+} - y_1^{-}) + (y_2^{+} - y_2^{-}) +\cdots + (y_n^{+} - y_n^{-}) \\
		&= \varphi_A(a) + (y_1^{+} + y_2^{+} + \cdots + y_n^{+}) - (y_1^{-} + y_2^{-} + \cdots + y_n^{-}).
	\end{align*}
	Set $p = y_1^{+} + y_2^{+} + \cdots + y_n^{+}$ and $q = y_1^{-} + y_2^{-} + \cdots + y_n^{-}$. Then $y = \varphi_A(x) = \varphi_A(a) + p - q$.
	
	For $\varphi_A(m) \in \varphi_A(U_i) \subseteq \varphi_A(U)$ and $\varphi_A(n) \in \varphi_A(U_{i+1}) \subseteq \varphi_A(U_i) \subseteq \varphi_A(U)$, we have $m \odot m = n \odot n = 0$.
	In the $\ell$-group $G_A$, $(\varphi_A(m) + \varphi_A(m) - u_{G_A}) \vee 0 = 0$, so $2\varphi_A(m) - u_{G_A} \le 0$; similarly $2\varphi_A(n) - u_{G_A} \le 0$. Adding these two inequalities gives $2(\varphi_A(m) + \varphi_A(n) - u_{G_A}) \le 0$, which implies $\varphi_A(m) + \varphi_A(n) - u_{G_A} \le 0$ by Proposition \ref{prop:lgroup-basic}(4), and consequently $(\varphi_A(m) + \varphi_A(n) - u_{G_A}) \vee 0 = 0$.
	
	Using Lemma \ref{lem:lgroup-identity} with $a = \varphi_A(m) + \varphi_A(n)$, we have
	\[
	\varphi_A(m) + \varphi_A(n) = (\varphi_A(m) + \varphi_A(n)) \wedge u_{G_A} + ((\varphi_A(m) + \varphi_A(n)) - u_{G_A}) \vee 0.
	\]
	Since $((\varphi_A(m) + \varphi_A(n)) - u_{G_A}) \vee 0 = 0$, we obtain $\varphi_A(m) + \varphi_A(n) = (\varphi_A(m) + \varphi_A(n)) \wedge u_{G_A} = \varphi_A(m \oplus n)$. By construction, $\varphi_A(U_i \oplus U_i) \subseteq \varphi_A(V_i \oplus V_i) \subseteq \varphi_A(U_{i-1})$, so $\varphi_A(m) + \varphi_A(n) \in \varphi_A(U_{i-1})$.
	
	By induction, we obtain $p = y_1^{+} + y_2^{+} + \cdots + y_n^{+} \in \varphi_A(U_0)$ and similarly $q = y_1^{-} + y_2^{-} + \cdots + y_n^{-} \in \varphi_A(U_0)$.
	
	Now $y = \varphi_A(x) = \varphi_A(a) + p - q \in \varphi_A(A)$. Then
	\[
	\varphi_A(x \ominus a) = (\varphi_A(x) - \varphi_A(a)) \vee 0 = (p - q) \vee 0 \in \varphi_A(U_0 \ominus U_0) \subseteq \varphi_A(W),
	\]
	and
	\[
	\varphi_A(a \ominus x) = (\varphi_A(a) - \varphi_A(x)) \vee 0 = (q - p) \vee 0 \in \varphi_A(U_0 \ominus U_0) \subseteq \varphi_A(W).
	\]
	Hence $y \in \varphi_A(W(a))$, proving
	\[
	[0, u_{G_A}] \cap \left( \varphi_A(a) + \sum_{i=1}^{\infty} (\varphi_A(U_i))^{\vartriangle} \right) \subseteq \varphi_A(W(a)).
	\]
	Thus $(A,\tau)$ is summable.
	
	In particular, taking $a = 0$ in the above proof, for any
	\[
	y = (y_1^+ + y_2^+ + \cdots + y_n^+) - (y_1^- + y_2^- + \cdots + y_n^-) = p - q \in \sum_{i=1}^{\infty} (\varphi_A(U_i))^{\vartriangle},
	\]
	we have
	\[
	y^+ = (p - q) \vee 0 \in \varphi_A(U_0 \ominus U_0) \subseteq \varphi_A(W), \quad
	y^- = (q - p) \vee 0 \in \varphi_A(U_0 \ominus U_0) \subseteq \varphi_A(W),
	\]
	hence $y \in (\varphi_A(W))^{\vartriangle}$. Consequently,
	\begin{equation}\label{inc:forward}
		\sum_{i=1}^{\infty} (\varphi_A(U_i))^{\vartriangle} \subseteq (\varphi_A(W))^{\vartriangle}.
	\end{equation}
	Conversely, for any neighbourhood $\sum_{i=1}^{\infty} (\varphi_A(W_i))^{\vartriangle}$ of $0$ in $\tau_{G_A}$, we have
	\begin{equation}\label{inc:reverse}
		(\varphi_A(W_1))^{\vartriangle} \subseteq \sum_{i=1}^{\infty} (\varphi_A(W_i))^{\vartriangle}.
	\end{equation}
	Therefore, by \eqref{inc:forward} and \eqref{inc:reverse}, the family $\{(\varphi_A(V))^{\vartriangle} \mid V \in \mathcal{V}(0) \}$ forms a neighbourhood basis of $0$ in $\tau_{G_A}$.
\end{proof}

Having established a general sufficient condition for summability, we now examine several important classes where this condition is automatically satisfied. We begin with totally ordered topological MV-algebras.

\begin{corollary}\label{cor-totallyorderedMV}
	Let $(A,\tau)$ be a totally ordered topological MV-algebra and let $\mathcal{V}(0)$ be a neighbourhood basis of $0$. If $\tau$ is not
	the trivial topology, then there exists \(U\in\mathcal V(0)\) such that
	\[
		\forall x\in U,\qquad x\odot x=0.
	\]
	Consequently, every totally ordered topological MV-algebra is summable.
\end{corollary}
\begin{proof}
	Let \(I:=\bigcap\mathcal V(0)\). Clearly \(0\in I\). If \(x,y\in I\)
	and \(U\in\mathcal V(0)\), choose \(V\in\mathcal V(0)\) such that
	\(V\oplus V\subseteq U\). Since \(x,y\in V\), we have
	\(x\oplus y\in U\), and hence \(x\oplus y\in I\).
	
	Now let \(y\leq x\in I\). For any \(U\in\mathcal V(0)\), the
	continuity of \(z\mapsto z\wedge y\) at \(0\), together with
	\(0\wedge y=0\), yields \(V\in\mathcal V(0)\) such that
	\(V\wedge y\subseteq U\). Since \(x\in V\) and \(y=x\wedge y\), we
	obtain \(y\in U\). Thus \(y\in I\), and therefore \(I\) is an ideal
	of \(A\).
	
	If \(I=A\), then the only neighbourhood of \(0\) is \(A\) itself,
	so \(\tau\) is the trivial topology. Endowing \(G_A\) with the
	trivial group topology, we conclude from
	Theorem~\ref{thm:iffMV} that \((A,\tau)\) is summable.
	
	If $I \subsetneq A$, then $1 \notin I$, so there exists $U \in \mathcal{V}(0)$ such that $1 \notin U$. By continuity of $\oplus$ at $(0,0)$, there exists $V \in \mathcal{V}(0)$ with $V \oplus V \subseteq U$. For any $x \in V$, we have $x \oplus x \in U$, hence $x \oplus x \neq 1$. By \eqref{prop:chain-property}, this implies $x \odot x = 0$. Thus $V$ satisfies the hypothesis of Proposition \ref{prop:U-odot-U-zero}, and consequently $(A, \tau)$ is summable.
\end{proof}

For comparison, Luan, Weber and Yang~\cite[Theorem~4.7]{luanweberyang2021} proved that every MV-chain with its order (interval) topology topology is a topological MV-algebra; consequently,  every such chain falls within the scope of Corollary~\ref{cor-totallyorderedMV}.

Perfect MV-algebras satisfy the same property: they admit a neighbourhood $U$ of $0$ with $x\odot x=0$ for all $x\in U$, and are therefore summable.

\begin{corollary}\label{cor-perfectMV}
	Let \((A,\tau)\) be a perfect topological MV-algebra and let
	\(\mathcal V(0)\) be a neighbourhood basis of \(0\). If \(\tau\) is
	not the trivial topology, then there exists \(U\in\mathcal V(0)\)
	such that
	\[
	\forall x\in U,\qquad x\odot x=0.
	\]
	Consequently, every perfect topological MV-algebra is summable.
\end{corollary}

\begin{proof}
	Let \(I:=\bigcap\mathcal V(0)\). As shown in the proof of
	Corollary~\ref{cor-totallyorderedMV}, \(I\) is an ideal of \(A\).
	
	Suppose that no \(U\in\mathcal V(0)\) is contained in
	\(\operatorname{Rad}(A)\). Let \(O\) be an arbitrary
	\(0\)-neighbourhood. By continuity of \(\oplus\) at \((0,0)\), there
	exists \(V\in\mathcal V(0)\) such that \(V\oplus V\subseteq O\).
	By assumption, there exists
	\(x\in V\setminus\operatorname{Rad}(A)\).
	
	Since \(A\) is perfect,
	\(x\in\neg\operatorname{Rad}(A)\), and hence
	\(\neg x\in\operatorname{Rad}(A)\). By
	\eqref{eq:perfect-radical-square-zero},
	\((\neg x)\odot(\neg x)=0\). Therefore
	\[
	\neg(x\oplus x)
	=
	(\neg x)\odot(\neg x)
	=
	0,
	\]
	and hence \(x\oplus x=1\). It follows that
	\(1\in V\oplus V\subseteq O\).
	
	Since \(O\) was arbitrary, \(1\in I\). Thus \(I=A\), so every
	\(0\)-neighbourhood is equal to \(A\), and \(\tau\) is the trivial
	topology.
	
	Consequently, if \(\tau\) is not trivial, there exists
	\(U\in\mathcal V(0)\) such that
	\(U\subseteq\operatorname{Rad}(A)\). By
	\eqref{eq:perfect-radical-square-zero}, every \(x\in U\) satisfies
	\(x\odot x=0\). Proposition~\ref{prop:U-odot-U-zero} therefore
	implies that \((A,\tau)\) is summable. The trivial topology is also
	summable.
\end{proof}

In addition to totally ordered and perfect topological MV-algebras, another class of summable examples arises from ideals. Let $I$ be an ideal of an MV-algebra $A$. A routine verification shows that the collection of equivalence classes $\{a/I\mid a\in A\}$ forms a basis for a topology on $A$; we call this topology the \emph{ideal topology}  $ \tau_I$ induced by $I$ (see~\cite[Proposition 4.1]{gan2024}). The following result, already used in the proof of Corollary~\ref{cor-totallyorderedMV}, characterizes precisely when a topology on an MV-algebra arises in this way.

\begin{proposition}[{\cite[Theorem 4.2]{gan2024}}]\label{prop:ideal-topology}
	Let $(A,\tau)$ be a topological MV-algebra and let $\mathcal{V}(0)$ be a neighbourhood basis of $0$. Then the following hold:
	\begin{enumerate}[label=(\arabic*)]
		\item $\bigcap \mathcal{V}(0)$ is an ideal of $A$;
		\item $\tau$ is induced by an ideal of $A$ if and only if $\bigcap \mathcal{V}(0) \in \tau$.
	\end{enumerate}
\end{proposition}

As a direct consequence of Proposition \ref{prop:ideal-topology}, Proposition \ref{prop:ideal-generated-summable} confirms that any topological MV-algebra whose topology is induced by an ideal is automatically summable.

\begin{proposition}\label{prop:ideal-generated-summable}
	Let $A$ be an MV-algebra and $I$ an ideal of $A$. Then the topological MV-algebra $(A, \tau_I)$ induced by $I$ is summable.
\end{proposition}

\begin{proof}
	Since $I$ is the smallest open neighbourhood of $0$ by Proposition \ref{prop:ideal-topology}, it suffices to verify that for any $a \in A$,
	\[
	[0, u_{G_A}] \cap \left( \varphi_A(a) + \sum_{i=1}^{\infty} (\varphi_A(I))^{\vartriangle} \right) \subseteq \varphi_A(I(a)).
	\]
	
	To this end, we first prove that
	\[
	\sum_{i=1}^{\infty} (\varphi_A(I))^{\vartriangle} = \widehat{\varphi_A(I)} := \left\{ x \in G_A \;\middle|\; |x| \wedge u_{G_A} \in \varphi_A(I) \right\}.
	\]
	By Theorem \ref{thm:ideal-correspondence}, $\widehat{\varphi_A(I)}$ is an $\ell$-ideal of $G_A$ and $\widehat{\varphi_A(I)} \cap [0, u_{G_A}] = \varphi_A(I)$.
	
	Take any $x \in \sum_{i=1}^{\infty} (\varphi_A(I))^{\vartriangle}$, i.e., $x = x_1 + x_2 + \cdots + x_n$ with $x_k \in (\varphi_A(I))^{\vartriangle}$, so $x_k^+, x_k^- \in \varphi_A(I)$. Then
	\[
	x = (x_1^+ + x_2^+ + \cdots + x_n^+) - (x_1^- + x_2^- + \cdots + x_n^-).
	\]
	Note that
	\[
	|(x_1^+  + x_2^+ + \cdots + x_n^+)| \wedge u_{G_A} = (x_1^+ + x_2^+ + \cdots + x_n^+) \wedge u_{G_A} = \varphi_A(y_1 \oplus y_2 \oplus \cdots \oplus y_n) \in \varphi_A(I),
	\]
	where $y_i = \varphi_A^{-1}(x_i^+)$. Hence $x_1^+ + x_2^+ +\cdots + x_n^+ \in \widehat{\varphi_A(I)}$. Similarly $x_1^- + x_2^- + \cdots + x_n^- \in \widehat{\varphi_A(I)}$. Since $\widehat{\varphi_A(I)}$ is a subgroup, $x = (x_1^+ + x_2^+ + \cdots + x_n^+) - (x_1^- + x_2^- + \cdots + x_n^-) \in \widehat{\varphi_A(I)}$.
	
	Conversely, take any $x \in \widehat{\varphi_A(I)}$, i.e., $|x| \wedge u_{G_A} \in \varphi_A(I)$. Since $u_{G_A}$ is a strong unit of $G_A$, there exists $n \in \mathbb{N}$ such that $x^+, x^- \le n u_{G_A}$. Using Lemma \ref{lem:lgroup-identity} inductively, we obtain
	\[
	x^+ = (x^+ \wedge u_{G_A}) + \bigl( (x^+ - u_{G_A}) \vee 0 \bigr) \wedge u_{G_A} + \cdots + \bigl( (x^+ - (n-1)u_{G_A}) \vee 0 \bigr) \wedge u_{G_A} + (x^+ - n u_{G_A}) \vee 0,
	\]
	\[
	x^- = (x^- \wedge u_{G_A}) + \bigl( (x^- - u_{G_A}) \vee 0 \bigr) \wedge u_{G_A} + \cdots + \bigl( (x^- - (n-1)u_{G_A}) \vee 0 \bigr) \wedge u_{G_A} + (x^- - n u_{G_A}) \vee 0.
	\]
	Set $A_i = \bigl( (x^+ - i u_{G_A}) \vee 0 \bigr) \wedge u_{G_A}$, $B_i = \bigl( (x^- - i u_{G_A}) \vee 0 \bigr) \wedge u_{G_A}$ for $i = 0, 1, \cdots, n$. Then $A_i, B_i \le |x| \wedge u_{G_A} \in \varphi_A(I)$. Since $I$ is an ideal of $A$, we have $A_i, B_i \in \varphi_A(I)$. In particular,
	$A_0 = x^+ \wedge u_{G_A},\; B_0 = x^- \wedge u_{G_A},\; A_n = (x^+ - n u_{G_A}) \vee 0,\; B_n = (x^- - n u_{G_A}) \vee 0$.
	Thus
	\[
	x^+ = \sum_{i=0}^{n} A_i, \qquad x^- = \sum_{i=0}^{n} B_i.
	\]
	Define $x_{i+1} = A_i - B_i$ for $i = 0, 1, \cdots, n$. Then $x = x^+ - x^- = x_1 +  x_2 + \cdots + x_{n+1}$.
	Since $A_i \wedge B_i = \bigl( (x^+ \wedge x^- - i u_{G_A}) \vee 0 \bigr) \wedge u_{G_A} = \bigl( (- i u_{G_A}) \vee 0 \bigr) \wedge u_{G_A} = 0$, we have
	$x_{i+1}^+ = (A_i - B_i) \vee 0 = A_i \in \varphi_A(I)$ and $x_{i+1}^- = B_i \in \varphi_A(I)$. Hence $x_{i+1} \in (\varphi_A(I))^{\vartriangle}$, and consequently $x = x_1 + x_2 + \cdots + x_{n+1} \in \sum_{i=1}^{\infty} (\varphi_A(I))^{\vartriangle}$.
	Thus we have proved
	\[
	\sum_{i=1}^{\infty} (\varphi_A(I))^{\vartriangle} = \widehat{\varphi_A(I)} := \left\{ x \in G_A \;\middle|\; |x| \wedge u_{G_A} \in \varphi_A(I) \right\}.
	\]
	
	Now,
	\[
	\begin{aligned}
		x\in[0,u_{G_A}]\cap\left(
		\varphi_A(a)+
		\sum_{i=1}^{\infty}(\varphi_A(I))^{\vartriangle}
		\right)
		&\iff
		x-\varphi_A(a)\in
		\sum_{i=1}^{\infty}(\varphi_A(I))^{\vartriangle}
		\\
		&\qquad\text{and }x\in[0,u_{G_A}]
		\\
		&\iff
		x-\varphi_A(a)\in\widehat{\varphi_A(I)}
		\text{ and }x\in[0,u_{G_A}]
		\\
		&\iff
		(x-\varphi_A(a))^+,\,
		(x-\varphi_A(a))^-
		\in\widehat{\varphi_A(I)}
		\\
		&\qquad\text{and }x\in[0,u_{G_A}]
		\\
		&\iff
		(x-\varphi_A(a))^+,\,
		(x-\varphi_A(a))^-
		\in\varphi_A(I)
		\\
		&\qquad\text{and }x\in[0,u_{G_A}]
		\\
		&\iff x\in\varphi_A(I(a)).
	\end{aligned}
	\]
	Therefore, $[0, u_{G_A}] \cap \left( \varphi_A(a) + \sum_{i=1}^{\infty} (\varphi_A(I))^{\vartriangle} \right) \subseteq \varphi_A(I(a))$ holds, and consequently $(A, \tau_I)$ is summable.
\end{proof}

Thus far we have identified three classes of summable topological MV-algebras: totally ordered algebras,
perfect algebras, and algebras whose topology is induced by an ideal.
Let us now exhibit several additional classes.

\begin{example}\label{ex:summable-classes}
	Let $(A,\tau)$ and $\{(A_i,\tau_i)\}_{i \in I}$ be topological MV-algebras. Then:
	\begin{enumerate}[label=(\arabic*)]
		\item If $A$ is a finite MV-algebra, then $(A,\tau)$ is summable. This follows from Propositions \ref{prop:ideal-topology} and \ref{prop:ideal-generated-summable}, because in a finite space the intersection of a neighbourhood basis of $0$ is open.
		\item If each $(A_i,\tau_i)$ is summable, then the direct product $(\prod_i A_i, \prod_i \tau_i)$ of MV-algebras with the product topology  is summable. This follows from Theorem \ref{thm:iffMV} together with the fact that the subspace of a product is the product of subspaces.
	\end{enumerate}
\end{example}

The previous results characterize those topological MV-algebras $(A,\tau)$ that are isomorphic (as an MV-algebra) and homeomorphic to $([0,u],\tau_G|_{[0,u]})$ for some abelian $\ell$-group $G$ with a strong unit $u\in G$, and group topology $\tau_G$ on $G$. A natural converse question is: given an abelian $\ell$-group $(G,u)$ and a group topology $\tau$ on $G$, when does $([0,u],\tau|_{[0,u]})$ form a topological MV-algebra?

\begin{definition}[MV lattice-ordered topological group]\label{def:MVLTG}
	A triple $(G, u, \tau)$ is called an \emph{MV lattice-ordered topological group} if $(G, u)$ is a lattice-ordered abelian group with strong unit $u$, $\tau$ is a group topology on $G$, and $([0, u], \tau|_{[0, u]})$ is a topological MV-algebra. We emphasize that the lattice operations on \(G\) are not required to be continuous.
\end{definition}

Since the MV-operations on \([0,u]\) are defined by  (Example \ref{ex:lgroup-to-MV})
\[
x\oplus y=(x+y)\wedge u,
\qquad
\neg x=u-x,
\]
it is immediate that if \(G\) is a topological $\ell$-group, i.e. if both
the group operations and the lattice operations on \(G\) are continuous, then
\(([0,u],\tau|_{[0,u]})\) is a topological MV-algebra.

Our notion of an MV lattice-ordered topological group does not require the lattice operations on $G$ to be continuous; it only demands a group topology on $G$ whose restriction to $[0,u]$ makes the interval into a topological MV-algebra. Example~\ref{ex:new1} demonstrates that this distinction is substantive: there exist group topologies on $G$ that do not make the lattice operations continuous, yet still yield a topological MV-algebra on $[0,u]$.

\begin{example}\label{ex:new1}
	Let \[ G=\left(\mathbb Z,+\right) \] with the usual order, and let \(u=1\). Then \((G,u)\) is an abelian $\ell$-group with strong unit \(1\). Its corresponding MV-algebra is $L_2 = \left\{0, 1\right\}$.
	Let \[ H=2\mathbb Z\subseteq \mathbb Z. \] Consider the filter basis \[ \mathcal F=\{H\}. \] It is easy to verify that $\mathcal{F}$ satisfies conditions (a), (b) and (c) of Theorem \ref{thm:group-topology}, hence it determines a unique group topology $\tau$ on $G$.
	Indeed, the corresponding topology $\tau$ has basis \[ \left\{ m+ 2\mathbb Z : m\in\mathbb Z \right\}. \]
	The lattice operations are not continuous with respect to $\tau$: for example, the map $x \mapsto x^+ = x \lor 0$ is discontinuous at $-1$. Indeed, every $\tau$-neighbourhood of $-1$ contains points (such as $1$) for which $x^+ \notin 2\mathbb{Z}$, while $(-1)^+ = 0 \in 2\mathbb{Z}$ --- contradicting continuity.

	Now consider the subspace topology $\tau|_{[0,1]}$. When restricted to $\{0,1\}$, it becomes the discrete topology, under which $\{0,1\}$ forms a topological MV-algebra (since all subsets are open, and the lattice operations are trivially continuous).
\end{example}

The following result characterizes precisely when an abelian $\ell$-group $(G,u)$ with a strong unit $u$, equipped with a group topology $\tau$, induces a topological MV-algebra structure on the unit interval $[0,u]$.

\begin{proposition}\label{prop:MV-lattice-group-condition}
	Let $(G, u)$ be an abelian $\ell$-group with a strong unit and let $\tau$ be a group topology on $G$, with $\mathcal{V}(0)$ a neighbourhood basis of $0$. Then $(G, u, \tau)$ is an MV lattice-ordered topological group if and only if for every $U \in \mathcal{V}(0)$ and every $a, b \in [0, u]$, there exists $V \in \mathcal{V}(0)$ such that
	\[
	\left(
	\bigl((a+V)\cap[0,u]\bigr)
	+
	\bigl((b+V)\cap[0,u]\bigr)
	\right)\wedge u
	\subseteq
	\bigl(((a+b)\wedge u)+U\bigr)\cap[0,u].
	\]
\end{proposition}

\begin{proof}
	The condition in the proposition is equivalent to the continuity of $\oplus$ on $([0, u], \tau|_{[0,u]})$, so necessity is clear.

	For sufficiency, it remains to show that $\neg$ is continuous. Take any $U \in \mathcal{V}(0)$. Since $\tau$ is a group topology, there exists $V \in \mathcal{V}(0)$ such that $-V \subseteq U$. Then for any $a \in [0, u]$, we have
	\[
	\neg\big( (a + V) \cap [0, u] \big) \subseteq (\neg a + U) \cap [0, u],
	\]
	which implies continuity of $\neg$. This completes the proof.
\end{proof}

\medskip
In this section we have introduced summability and proved that it
characterises exactly when a topological MV-algebra can be obtained from a
lattice-ordered topological group (Theorem~\ref{thm:iffMV}). Conversely,
Proposition~\ref{prop:MV-lattice-group-condition} identifies when the unit
interval of a lattice-ordered topological group is a topological MV-algebra.
Several important classes of summable algebras — including totally ordered,
perfect, ideal-induced, finite, and product algebras — are identified.
These results provide the foundation for lifting  Mundici and Di Nola--
Lettieri equivalences to the topological setting, which we undertake in
the next section.



\section{Categorical framework}\label{sec:adjunction}
Building on the concrete realisation established in
Section~\ref{sec:summable}, we now lift this correspondence to the
categorical level. Let $\mathcal{T}\!\mathit{op}^{+}\!\mathcal{MV}$ denote
the category of summable topological MV-algebras
(Definition~\ref{def:sumMV}), with continuous MV-algebra homomorphisms as
morphisms. Dually, let $\ell\mathcal{T}\mathcal{G}_u^{\mathcal{MV}}$ be
the category of MV lattice-ordered topological groups
(Definition~\ref{def:MVLTG}), whose morphisms are continuous $\ell$-group
homomorphisms that preserve the strong unit.

The main result of this section is an adjunction
\[
\Xi\dashv\Gamma: \mathcal{T}\!\mathit{op}^{+}\!\mathcal{MV} \rightharpoonup \ell\mathcal{T}\mathcal{G}_u^{\mathcal{MV}}
\]
in Theorem~\ref{thm:adjunction}, providing a functorial solution to the
topological lifting problem for  Mundici equivalence. From this, we
derive the notion of strong summability
(Definition~\ref{def:strong-summable}), which characterises when the
adjunction lands in the subcategory of topological $\ell$-groups. This
yields a restricted adjunction between strongly summable topological
MV-algebras and topological $\ell$-groups
(Corollary~\ref{cor:strong-adjunction}). We show that this restricted
adjunction is not an equivalence in general
(Example~\ref{ex:non-unique-ambient-topology}), but becomes one when
restricted to totally ordered algebras
(Corollary~\ref{cor:totally-ordered-equivalence}), thereby lifting the
Mundici equivalence to the topological setting. Finally, we turn to
perfect topological MV-algebras and prove that Di Nola--Lettieri
equivalence also admits a topological lifting
(Theorem~\ref{thm:perfect-topological-equivalence}). As a first step, we verify that $\Gamma$ and $\Xi$ preserve continuity of
morphisms and therefore define genuine functors between the two categories.

Recall from Section~\ref{sec:summable} that every summable topological
MV-algebra $(A,\tau)$ determines a lattice-ordered topological group
$(\Gamma(A),\tau_\Gamma)$, where $\Gamma(A)$ is the Mundici group of $A$
and $\tau_\Gamma$ is the unique group topology whose restriction to
$[0,u_{G_A}]$ recovers $\tau$. The assignment $(A,\tau)\mapsto
(\Gamma(A),\tau_\Gamma)$ defines a functor
\[
\Gamma: \mathcal{T}\!\mathit{op}^{+}\!\mathcal{MV} \longrightarrow \ell\mathcal{T}\mathcal{G}_u^{\mathcal{MV}}.
\]

Conversely, given an MV lattice-ordered topological group $(G,u,\tau_G)$,
its unit interval $[0,u]$ endowed with the subspace topology inherits a
topological MV-algebra structure. This assignment defines a functor
\[
\Xi: \ell\mathcal{T}\mathcal{G}_u^{\mathcal{MV}} \longrightarrow \mathcal{T}\!\mathit{op}^{+}\!\mathcal{MV},
\]
which sends $(G,u,\tau_G)$ to $([0,u],\tau_G|_{[0,u]})$ and acts on
morphisms by restriction.

\begin{lemma}\label{lem:functor-continuity}
	For any morphism $f: (G,u,\tau_G) \to (H,v,\tau_H)$ in
	$\ell\mathcal{T}\mathcal{G}_u^{\mathcal{MV}}$, the map
	$\Gamma(f) = f|_{[0,u]}: ([0,u], \tau_G|_{[0,u]}) \to ([0,v], \tau_H|_{[0,v]})$
	is continuous.

	For any morphism $h: (A,\tau_{G_A}) \to (B,\tau_B)$ in
	$\mathcal{T}\!\mathit{op}^{+}\!\mathcal{MV}$, the map
	$\Xi(h): (G_A,u_{G_A},\tau_{G_A}) \to (G_B,u_B,\tau_{G_B})$ is continuous, with
	$\tau_{G_A}$ and $\tau_{G_B}$ constructed in Proposition~\ref{prop:starcon}
	(see also \eqref{eq:V_GA}).
\end{lemma}

\begin{proof}
	The first statement follows directly from the fact that
	$\Gamma(f)^{-1}(U \cap [0,v]) = f^{-1}(U) \cap [0,u]$ for any $U \subseteq H$,
	so we only prove the second one.

	Let $h: (A,\tau_{G_A}) \to (B,\tau_B)$ be a morphism in
	$\mathcal{T}\!\mathit{op}^{+}\!\mathcal{MV}$; we show that $\Xi(h)$ is
	continuous.

	Take any $x \in G_A$ and any basic neighbourhood of $\Xi(h)(x)$ in $\tau_{G_B}$,
	i.e.\ a set of the form $\Xi(h)(x) + \sum_{i=1}^{\infty} (\varphi_B(V_i))^{\vartriangle}$
	as in \eqref{eq:V_GA}, where each $V_i$ is an open neighbourhood of $0$ in
	$\tau_B$.  Set $U_i = h^{-1}(V_i)$; then each $U_i$ is an open neighbourhood of
	$0$ in $A$, and
	\(
	x \in x + \sum_{i=1}^{\infty} (\varphi_A(U_i))^{\vartriangle}.
	\)

	For any $y \in \Xi(h)\big( x + \sum_{i=1}^{\infty} (\varphi_A(U_i))^{\vartriangle} \big)$,
	we have
	\(
	y = \Xi(h)(x + x_1 + x_2 + \cdots + x_n)
	  = \Xi(h)(x) + \Xi(h)(x_1) + \Xi(h)(x_2) + \cdots + \Xi(h)(x_n),
	\)
	with $x_i^{+}, x_i^{-} \in \varphi_A(U_i)$.  Then
	$(\Xi(h)(x_i))^{+} = \Xi(h)(x_i^{+}) \in \Xi(h)(\varphi_A(U_i))
	= \Gamma\Xi(h)(\varphi_A(U_i)) = \varphi_B(h(U_i)) \subseteq \varphi_B(V_i)$,
	where the equality $\Gamma\Xi(h)(\varphi_A(U_i)) = \varphi_B(h(U_i))$ follows
	from the first diagram in \eqref{diag:naturality-varphi-psi}; similarly
	$(\Xi(h)(x_i))^{-} \in \varphi_B(V_i)$.  Hence
	$\Xi(h)(x_i) \in (\varphi_B(V_i))^{\vartriangle}$, and therefore
	\(
	y \in \Xi(h)(x) + \sum_{i=1}^{\infty} (\varphi_B(V_i))^{\vartriangle}.
	\)

	Thus
	$\Xi(h)\big( x + \sum_{i=1}^{\infty} (\varphi_A(U_i))^{\vartriangle} \big)
	\subseteq \Xi(h)(x) + \sum_{i=1}^{\infty} (\varphi_B(V_i))^{\vartriangle}$,
	which proves the continuity of $\Xi(h)$.
\end{proof}

In view of Lemma~\ref{lem:functor-continuity},
$\Gamma: \ell\mathcal{T}\mathcal{G}_u^{\mathcal{MV}} \to
\mathcal{T}\!\mathit{op}^{+}\!\mathcal{MV}$ and
$\Xi: \mathcal{T}\!\mathit{op}^{+}\!\mathcal{MV} \to
\ell\mathcal{T}\mathcal{G}_u^{\mathcal{MV}}$ are well-defined functors.
Explicitly,
$\Gamma(G,u,\tau_G) = ([0,u], \tau_G|_{[0,u]})$ and
$\Xi(A,\tau) = (G_A,u_{G_A},\tau_{G_A})$, where $\tau_{G_A}$ is the group topology
constructed in Proposition~\ref{prop:starcon}.  On the underlying algebraic
structures, $\Gamma(f)$ and $\Xi(h)$ are defined exactly as in Mundici
equivalence; Lemma~\ref{lem:functor-continuity} guarantees that these maps remain
continuous with respect to the respective topologies.  Moreover,
Proposition~\ref{prop:starcon} shows that the natural isomorphism
$\varphi: \mathrm{id}_{\mathcal{T}\!\mathit{op}^{+}\!\mathcal{MV}} \to \Gamma\Xi$
from Mundici equivalence is a homeomorphism for each summable topological
MV-algebra, hence remains a natural isomorphism in the topological setting. 
To prove the adjunction, we need the following lemma.

\begin{lemma}\label{lem:XiGamma-finer}
	Let $(G,u,\tau) \in \ell\mathcal{T}\mathcal{G}_u^{\mathcal{MV}}$ and let
	$\tau'$ denote the topology of $\Xi\Gamma(G,u,\tau)$.  Then $\tau'$ is finer
	than $\psi_G(\tau)$, where $\psi_G: G \to \Xi\Gamma(G,u,\tau)$ is the
	isomorphism given by Mundici equivalence.
\end{lemma}

\begin{proof}
	Denote $\Xi\Gamma(G,u,\tau)$ by $(G', u', \tau')$.  Take any $x\in G'$ and
	any basic neighbourhood of $x$ in $\psi_G(\tau)$ of the form
	$x + \psi_G(U)$, where $U$ is an open neighbourhood of $0$ in $\tau$.  It
	suffices to show that $x$ is an interior point of $x + \psi_G(U)$ in $\tau'$,
	i.e.\ that $0$ is an interior point of $\psi_G(U)$ in $\tau'$.

	By Lemma~\ref{lem:sum-of-triangles}, applied to an open neighbourhood basis
	of $0$, for $U$ there exists a sequence $\{V_i\}_{i=1}^{\infty}$ of open
	neighbourhoods of $0$ such that
	$\sum_{i=1}^{\infty}V_i^{\vartriangle}\subseteq U$.
	Consequently, in $G'$ we have
	\[
	0 \in \sum_{i=1}^{\infty} \big( \psi_G(V_i \cap [0,u]) \big)^{\vartriangle}
	    \subseteq \sum_{i=1}^{\infty} (\psi_G(V_i))^{\vartriangle}
	    \subseteq \psi_G(U).
	\]
	Since $\varphi_{[0,u]} = \psi_G|_{[0,u]}$ (see \eqref{eq:psi-varphi}), we have
	\[
	\psi_G(V_i \cap [0,u]) = \varphi_{[0,u]}(V_i \cap [0,u]).
	\]
	Hence, by the definition of the topology $\tau'$, the set
	\[
	\sum_{i=1}^{\infty} \big( \psi_G(V_i \cap [0,u]) \big)^{\vartriangle}
	\]
	is open in $\tau'$.  This completes the proof.
\end{proof}

\begin{remark}\label{rem:psi-continuity}
	Let $(G,u,\tau) \in \ell\mathcal{T}\mathcal{G}_u^{\mathcal{MV}}$.
	\begin{enumerate}[label=(\arabic*)]
		\item By Lemma~\ref{lem:XiGamma-finer},
		      $\psi_G^{-1}: \Xi\Gamma(G,u,\tau) \to (G,u,\tau)$ is continuous.
		\item Moreover, $\tau'$ is the finest topology among those arising from
		      MV lattice-ordered topological group structures on the Mundici group
		      whose restriction to the unit interval coincides with the topology of
		      $\Gamma(G,u,\tau)$.  Indeed, the proof of Lemma~\ref{lem:XiGamma-finer}
		      shows that any such topology is coarser than $\tau'$.
		\item However, $\psi_G$ need not be a homeomorphism in general; the
		      topology of $\Xi\Gamma(G,u,\tau)$ may be strictly finer than
		      $\psi_G(\tau)$.  For instance, take $(G,u) = (\mathbb{R}, 1)$ and
		      let $\tau$ be the group topology on $\mathbb{R}$ for which
		      $\{2\mathbb{Z}\}$ is a neighbourhood basis of $0$.  Then
		      $\Xi\Gamma(G,u,\tau)$ is discrete (since $\Gamma(G,u,\tau)$ is
		      discrete), which is strictly finer than $\psi_G(\tau)$.
	\end{enumerate}
\end{remark}

Although $\psi_G$ is not a homeomorphism in general
(Remark~\ref{rem:psi-continuity}), the situation improves for objects arising
from the functor $\Xi$: for $(G_A,u_{G_A},\tau_{G_A}) = \Xi(A,\tau)$, the map
$\psi_{G_A}$ is actually an isomorphism in $\ell\mathcal{T}\mathcal{G}_u^{\mathcal{MV}}$.

\begin{lemma}\label{lem:psi-isomorphism}
	For any $(A,\tau) \in \mathcal{T}\!\mathit{op}^{+}\!\mathcal{MV}$, let
	$(G_A,u_{G_A},\tau_{G_A}) = \Xi(A,\tau)$.  Then the map
	$\psi_{G_A}: (G_A,u_{G_A},\tau_{G_A}) \to
	\Xi\Gamma(G_A,u_{G_A},\tau_{G_A})$ is an isomorphism in
	$\ell\mathcal{T}\mathcal{G}_u^{\mathcal{MV}}$.
\end{lemma}

\begin{proof}
	It suffices to show that $\psi_{G_A}$ is a homeomorphism.  By
	Proposition~\ref{prop:starcon}, a $0$-neighbourhood basis in
	$(G_A,u_{G_A},\tau_{G_A})$ is given by
	\[
	\left\{
	\sum_{i=1}^{\infty}(\varphi_A(U_i))^{\vartriangle}:
	U_i\in\mathcal V_A(0)
	\right\},
	\]
	and a $0$-neighbourhood basis in $\Xi\Gamma(G_A,u_{G_A},\tau_{G_A})$ is given by
	\[
	\left\{
	\sum_{i=1}^{\infty}(\varphi_{[0,u_{G_A}]}(V_i))^{\vartriangle}:
	V_i\in\mathcal V_{[0,u_{G_A}]}(0)
	\right\}.
	\]

	Since $\varphi_A: A \to [0,u_{G_A}]$ is a homeomorphism, it restricts to a
	bijection between open neighbourhoods $\mathcal V_A(0)$ of $0$ in $A$ and
	those $\mathcal V_{[0,u_{G_A}]}(0)$ in $[0,u_{G_A}]$.  Thus each
	$V_i\in\mathcal V_{[0,u_{G_A}]}(0)$ can be written uniquely as $\varphi_A(U_i)$
	for some $U_i\in\mathcal V_A(0)$.  Moreover, because $\psi_{G_A}$ is an
	$\ell$-group isomorphism and $\psi_{G_A}|_{[0,u_{G_A}]} = \varphi_{[0,u_{G_A}]}$
	(see \eqref{eq:psi-varphi}), we have
	\[
	\psi_{G_A}\big((\varphi_A(U_i))^{\vartriangle}\big)
	=
	\big(\psi_{G_A}(\varphi_A(U_i))\big)^{\vartriangle}
	=
	\big(\varphi_{[0,u_{G_A}]}(\varphi_A(U_i))\big)^{\vartriangle}
	=
	(\varphi_{[0,u_{G_A}]}(V_i))^{\vartriangle}.
	\]

	Hence $\psi_{G_A}$ sends the first neighbourhood basis onto the second.
	Since both $(G_A,u_{G_A},\tau_{G_A})$ and $\Xi\Gamma(G_A,u_{G_A},\tau_{G_A})$ are
	topological groups, a group isomorphism that maps a neighbourhood basis to a
	neighbourhood basis is a homeomorphism.  Therefore $\psi_{G_A}$ is an
	isomorphism in $\ell\mathcal{T}\mathcal{G}_u^{\mathcal{MV}}$.
\end{proof}

With Lemmas~\ref{lem:functor-continuity}–\ref{lem:psi-isomorphism} in place,
we can now prove the second main result of this paper . 

\begin{theorem}\label{thm:adjunction}
	$\Xi \dashv \Gamma: \mathcal{T}\!\mathit{op}^{+}\!\mathcal{MV}
	\rightharpoonup \ell\mathcal{T}\mathcal{G}_u^{\mathcal{MV}}.$
\end{theorem}

\begin{proof}
	By Proposition~\ref{prop:adj-characterization}, it suffices to show that the
	natural transformation
	\[
	\varphi: \mathrm{id}_{\mathcal{T}\!\mathit{op}^{+}\!\mathcal{MV}} \to \Gamma\Xi
	\]
	satisfies the following universal property: for every object $(G,u,\tau')$ in
	$\ell\mathcal{T}\mathcal{G}_u^{\mathcal{MV}}$, each object $(A,\tau)$ in
	$\mathcal{T}\!\mathit{op}^{+}\!\mathcal{MV}$, and every morphism
	$g \in \mathrm{Hom}\big((A,\tau),([0,u],\tau'|_{[0,u]})\big)$,
	there exists a unique morphism
	$f \in \mathrm{Hom}\big((G_A,u_{G_A},\tau_{G_A}),(G,u,\tau')\big)$
	such that $\Gamma(f)\circ\varphi_A=g$, as expressed by the following
	commutative diagram:
	\[
	\begin{tikzcd}[row sep=large, column sep=large]
		(A,\tau) \arrow[r, "\varphi_A"] \arrow[dr, "g"']
		& ([0,u_{G_A}],\tau_{G_A}|_{[0,u_{G_A}]})
		\arrow[d, dashed, "\Gamma(f)\;(\exists!\,f)"] \\
		& ([0,u],\tau'|_{[0,u]})
	\end{tikzcd}
	\]

	The morphism $g$ is a continuous MV-algebra homomorphism, hence in particular
	a morphism in $\mathcal{MV}$.  By Mundici equivalence
	(Theorem~\ref{thm:mundici}), there exists a unique $\ell$-group homomorphism
	$f: G_A \to G$ preserving the strong unit such that
	$\Gamma(f) \circ \varphi_A = g$.  It remains to show that $f$ is continuous.

	Since $\varphi_A$ is a homeomorphism,
	$\Gamma(f) = g \circ \varphi_A^{-1}$ is continuous.  By
	Lemma~\ref{lem:functor-continuity},
	$\Xi\Gamma(f): \Xi\Gamma(G_A,u_{G_A},\tau_{G_A}) \to
	\Xi\Gamma(G,u,\tau')$ is continuous.  Now, by the naturality of $\psi$, as
	expressed by the second diagram in \eqref{diag:naturality-varphi-psi}, we have
	$\Xi\Gamma(f)\circ\psi_{G_A}=\psi_G\circ f$.  Hence
	\[
	f = \psi_G^{-1} \circ \Xi\Gamma(f) \circ \psi_{G_A}.
	\]
	By Remark~\ref{rem:psi-continuity}, $\psi_G^{-1}$ is continuous, and by
	Lemma~\ref{lem:psi-isomorphism}, $\psi_{G_A}$ is continuous.  Being a
	composition of continuous maps, $f$ is continuous.
\end{proof}

\begin{remark}
	The adjunction $\Xi\dashv\Gamma$ in Theorem~\ref{thm:adjunction} is the topological counterpart of Mundici
	equivalence.  The compatibility condition on the MV-algebra side is summability; this is the topological counterpart, within Mundici's equivalence, of the requirement that the associated abelian $\ell$-group be equipped with a compatible group topology.
Di Nola--Lettieri equivalence receives the analogous
	perfect-topology lift in Theorem~\ref{thm:perfect-topological-equivalence}.
\end{remark}

\begin{remark}\label{rem:adjunction-properties}
	The adjunction established in Theorem~\ref{thm:adjunction} satisfies the
	following properties.
	\begin{enumerate}[label=(\arabic*)]
		\item The unit of the adjunction is $\varphi$ and the counit is
		      $\psi^{-1}$.  Equation~\eqref{eq:psi-varphi} encodes the required
		      triangle identities.
		\item The functor $\Xi$ is fully faithful.  Indeed, let
		      $\widetilde{g}: \Xi(A,\tau_{G_A}) \to \Xi(B,\tau_B)$ be a morphism in
		      $\ell\mathcal{T}\mathcal{G}_u^{\mathcal{MV}}$.  Since the algebraic
		      Mundici functor $\Xi: \mathcal{MV} \to \ell\mathcal{G}_u$ is fully
		      faithful, there exists a unique MV-algebra homomorphism
		      $g: A \to B$ such that $\Xi(g) = \widetilde{g}$ on the underlying
		      algebraic structures.  By the naturality of $\varphi$, as expressed
		      by the first diagram in \eqref{diag:naturality-varphi-psi},
		      \[
		      g = \varphi_B^{-1} \circ \Gamma(\widetilde{g}) \circ \varphi_A.
		      \]
		      Because $\varphi_A$ and $\varphi_B$ are homeomorphisms and
		      $\Gamma(\widetilde{g})$ is continuous, $g$ is continuous.  Hence
		      $\widetilde{g} = \Xi(g)$ for a unique continuous $g$, so $\Xi$ is
		      fully faithful.
	\end{enumerate}
\end{remark}

From Remark~\ref{rem:adjunction-properties}, we immediately obtain the following
universal property.

\begin{corollary}\label{cor:Xi-universal}
	For any summable topological MV-algebra $(A,\tau)$, the object
	$\Xi(A,\tau)$ in $\ell\mathcal{T}\mathcal{G}_u^{\mathcal{MV}}$ has the
	following universal property: for any object (MV lattice-ordered topological
	group) $(G,u,\tau_G)$ in $\ell\mathcal{T}\mathcal{G}_u^{\mathcal{MV}}$ and any
	morphism (continuous $\ell$-group homomorphism)
	$f: \Xi(A,\tau) \to (G,u,\tau_G)$ in
	$\ell\mathcal{T}\mathcal{G}_u^{\mathcal{MV}}$, there exists a unique morphism
	(continuous $\ell$-group homomorphism)
	$\tilde{g}: \Xi(A,\tau) \to \Xi\Gamma(G,u,\tau_G)$ in
	$\ell\mathcal{T}\mathcal{G}_u^{\mathcal{MV}}$ such that
	$f = \psi_G^{-1} \circ \tilde{g}$, as expressed by the following commutative
	diagram:
	\[
	\begin{tikzcd}[row sep=large, column sep=large]
		& \Xi\Gamma(G,u,\tau_G)
		\arrow[d, "\psi_G^{-1}"] \\
		\Xi(A,\tau)
		\arrow[ur, dashed, "\tilde{g}\;(\exists!)"]
		\arrow[r, "f"']
		& (G,u,\tau_G).
	\end{tikzcd}
	\]
\end{corollary}

Although the preceding construction was carried out in the broader category of
MV lattice-ordered topological groups, the representation simplifies when the
target is required to be a topological $\ell$-group.  The following proposition
shows that this occurs precisely when $\Xi(A,\tau)$ is itself a topological
$\ell$-group.

\begin{proposition}\label{prop:topological-lgroup-criterion}
	Let $(A,\tau)$ be a summable topological MV-algebra.  Then the following are
	equivalent.
	\begin{enumerate}[label=(\arabic*)]
		\item There exist an abelian $\ell$-group $G$, a positive element
		      $u\in G$, and a group topology $\tau_G$ on $G$ such that
		      $(G,u,\tau_G)$ is a topological $\ell$-group and $(A,\tau)$ is
		      isomorphic as an MV-algebra and homeomorphic to
		      $([0,u],\tau_G|_{[0,u]})$.
		\item The object $\Xi(A,\tau)=(G_A,u_{G_A},\tau_{G_A})$ is a topological
		      $\ell$-group; equivalently, the lattice operations on $G_A$ are
		      continuous with respect to $\tau_{G_A}$.
	\end{enumerate}
\end{proposition}

\begin{proof}
	The implication (2)$\Rightarrow$(1) is immediate.  If
	$\Xi(A,\tau)=(G_A,u_{G_A},\tau_{G_A})$ is a topological $\ell$-group, then by
	construction $\varphi_A:(A,\tau)\to([0,u_{G_A}],\tau_{G_A}|_{[0,u_{G_A}]})$ is an
	MV-isomorphism and a homeomorphism, so $(A,\tau)$ embeds into a topological
	$\ell$-group.

	We prove the converse.  Assume (1) holds.  The proof proceeds in three steps:
	first, we transport the topology from $G$ to $G_A$; second, we reduce the
	continuity of the lattice operations to the continuity of the positive-part
	map; third, we prove that the positive-part map is $\tau_{G_A}$-continuous
	using a neighbourhood basis argument.

	Suppose that $(A,\tau)$ embeds as the unit interval of a topological
	$\ell$-group $(G,\sigma)$, with distinguished positive element $u$.  After
	passing to the $\ell$-ideal $G_u$ generated by $u$ and transporting the
	topology $\sigma|_{G_u}$ along the canonical isomorphism $G_u \cong G_A$, we
	may assume that there exists a group topology $\sigma$ on $G_A$ such that
	$(G_A,u_{G_A},\sigma)$ is a topological $\ell$-group and
	$\varphi_A:(A,\tau)\to([0,u_{G_A}],\sigma|_{[0,u_{G_A}]})$ is a homeomorphism.

	This is the core of the proof.  We prove that the positive-part map
	$p : G_A \to G_A$, $p(g) = g^+$, is continuous with respect to $\tau_{G_A}$.
	Since $\tau_{G_A}$ is already a group topology, this implies that the lattice
	operations are continuous, because
	$x\vee y = y + (x-y)^+$ and $x\wedge y = x - (x-y)^+$.

	Fix $x\in G_A$.  Let
	\[
	W^{G_A} = \sum_{i=1}^{\infty}(\varphi_A(U_i^\tau))^{\vartriangle}
	\]
	be a basic $0$-neighbourhood in $\tau_{G_A}$, where each $U_i^\tau$ is a
	$0$-neighbourhood in $(A,\tau)$.  For $S\subseteq[0,u_{G_A}]$, write
	$S^{\vartriangle} = \{g\in G_A: g^+\in S,\ g^-\in S\}$.

	For each $i$, since $\varphi_A(U_i^\tau)$ is a $0$-neighbourhood in
	$([0,u_{G_A}],\sigma|_{[0,u_{G_A}]})$, choose a $\sigma$-open $0$-neighbourhood
	$O_i^\sigma$ in $G_A$ such that
	\begin{equation}\tag{$\mathrm{O}_i$}
		O_i^\sigma\cap[0,u_{G_A}]\subseteq \varphi_A(U_i^\tau).
	\end{equation}
	Since $(G_A,\sigma)$ is a topological $\ell$-group, the maps
	$g\mapsto g^+$ and $g\mapsto g^-$ on $G_A$ are $\sigma$-continuous.  Choose
	a $\sigma$-neighbourhood $N_i^\sigma$ of $0$ such that
	\begin{equation}\tag{$\mathrm{N}_i$}
		g\in N_i^\sigma \quad\Longrightarrow\quad
		g^+\in O_i^\sigma\ \text{and}\ g^-\in O_i^\sigma .
	\end{equation}
	By ($\mathrm{O}_i$) and ($\mathrm{N}_i$), whenever $g\in N_i^\sigma$ and
	$g^+,g^-\in[0,u_{G_A}]$, we have $g\in(\varphi_A(U_i^\tau))^{\vartriangle}$.

	Choose a $\sigma$-neighbourhood $M_i^\sigma$ of $0$ such that
	$M_i^\sigma-M_i^\sigma\subseteq N_i^\sigma$.  By the $\sigma$-continuity of
	the positive-part map at $x$, choose a $\sigma$-neighbourhood $L_i^\sigma$
	of $0$ such that $p(x+L_i^\sigma)\subseteq x^+ + M_i^\sigma$.  Hence, if
	$s,t\in L_i^\sigma$, then
	\begin{equation}\tag{$\mathrm{L}_i$}
		(x+s)^+-(x+t)^+\in M_i^\sigma-M_i^\sigma\subseteq N_i^\sigma .
	\end{equation}

	We now choose $\sigma$-neighbourhoods $P_i^\sigma$ recursively.  Choose
	$P_1^\sigma$ such that $P_1^\sigma+P_1^\sigma\subseteq L_1^\sigma$.  After
	$P_i^\sigma$ has been chosen, choose $P_{i+1}^\sigma$ such that
	$P_{i+1}^\sigma+P_{i+1}^\sigma\subseteq P_i^\sigma\cap L_{i+1}^\sigma$.
	Then a backward induction gives
	\begin{equation}\tag{$\mathrm{P}_i$}
		\sum_{j=i}^{\infty}P_j^\sigma\subseteq L_i^\sigma
		\qquad\text{for every }i,
	\end{equation}
	where the sum denotes the set of all finite sums with one summand chosen from
	finitely many $P_j^\sigma$'s.

	For each $i$, choose a $\sigma$-neighbourhood $Q_i^\sigma$ of $0$ such that
	$Q_i^\sigma-Q_i^\sigma\subseteq P_i^\sigma$.  Since
	$\varphi_A:(A,\tau)\to([0,u_{G_A}],\sigma|_{[0,u_{G_A}]})$ is a homeomorphism, choose
	a $0$-neighbourhood $V_i^\tau$ in $(A,\tau)$ such that
	$\varphi_A(V_i^\tau)\subseteq Q_i^\sigma\cap[0,u_{G_A}]$.  Then
	\begin{equation}\tag{$\mathrm{V}_i$}
		(\varphi_A(V_i^\tau))^{\vartriangle}\subseteq P_i^\sigma .
	\end{equation}
	Indeed, if $h\in(\varphi_A(V_i^\tau))^{\vartriangle}$, then $h=h^+-h^-$, with
	$h^+,h^-\in\varphi_A(V_i^\tau)\subseteq Q_i^\sigma$, so
	$h\in Q_i^\sigma-Q_i^\sigma\subseteq P_i^\sigma$.

	Set
	\[
	V^{G_A} = \sum_{i=1}^{\infty}(\varphi_A(V_i^\tau))^{\vartriangle}.
	\]
	We claim that $h\in V^{G_A}$ implies $(x+h)^+-x^+\in W^{G_A}$.

	Let $h\in V^{G_A}$.  After adding zero terms if necessary, there exist
	$n\in\mathbb N$ and elements $h_i\in(\varphi_A(V_i^\tau))^{\vartriangle}$,
	$i=1,2,\ldots,n$, such that $h=\sum_{i=1}^n h_i$.  For each
	$i=1,2,\ldots,n$, put $s_i=\sum_{j=i+1}^n h_j$.  By ($\mathrm{P}_i$) and
	($\mathrm{V}_i$), we have $s_i\in L_i^\sigma$ and $s_i+h_i\in L_i^\sigma$.
	Hence, by ($\mathrm{L}_i$),
	\[
	(x+s_i+h_i)^+-(x+s_i)^+\in N_i^\sigma .
	\]
	Define $\delta_i = (x+s_i+h_i)^+-(x+s_i)^+$.  We show that
	$\delta_i\in(\varphi_A(U_i^\tau))^{\vartriangle}$.

	Write $h_i = a_i - b_i$, where $a_i = h_i^+$ and $b_i = h_i^-$.  Since
	$h_i\in(\varphi_A(V_i^\tau))^{\vartriangle}$, we have
	$a_i,b_i\in\varphi_A(V_i^\tau)\subseteq[0,u_{G_A}]$.

	By Lemma~\ref{lem:lgroup-positive-part-estimate}, applied to
	$y = x+s_i$, $a = a_i$, and $b = b_i$, we get
	\[
	\delta_i^+ \leq a_i \qquad\text{and}\qquad \delta_i^- \leq b_i .
	\]
	Hence $\delta_i^+,\delta_i^-\in[0,u_{G_A}]$.  Since $\delta_i\in N_i^\sigma$,
	by ($\mathrm{N}_i$) we have $\delta_i^+\in O_i^\sigma$ and
	$\delta_i^-\in O_i^\sigma$.  Hence, by ($\mathrm{O}_i$),
	\[
	\delta_i^+,\delta_i^-\in O_i^\sigma\cap[0,u_{G_A}]\subseteq \varphi_A(U_i^\tau).
	\]
	Therefore $\delta_i\in(\varphi_A(U_i^\tau))^{\vartriangle}$.

	Finally, by the definition of $s_i$, we have
	\[
	\delta_i
	=
	\left(x+\sum_{j=i}^n h_j\right)^+
	-
	\left(x+\sum_{j=i+1}^n h_j\right)^+ .
	\]
	Therefore the sum telescopes:
	\[
	(x+h)^+ - x^+ = \sum_{i=1}^n \delta_i .
	\]
	Since each $\delta_i\in(\varphi_A(U_i^\tau))^{\vartriangle}$, it follows that
	\[
	(x+h)^+ - x^+ \in \sum_{i=1}^{\infty}(\varphi_A(U_i^\tau))^{\vartriangle}
	= W^{G_A}.
	\]
	Thus $h\in V^{G_A}$ implies $(x+h)^+-x^+\in W^{G_A}$.

	Since $W^{G_A}$ was an arbitrary basic $0$-neighbourhood in $\tau_{G_A}$,
	the positive-part map is continuous at $x$ with respect to $\tau_{G_A}$.
	Since $x\in G_A$ was arbitrary, $g\mapsto g^+$ is continuous everywhere.
	Hence $(G_A,u_{G_A},\tau_{G_A})$ is a topological $\ell$-group.
\end{proof}

In view of Proposition~\ref{prop:topological-lgroup-criterion}, we introduce
the following strengthening of summability.

\begin{definition}\label{def:strong-summable}
	A summable topological MV-algebra $(A,\tau)$ is called
	\emph{strongly summable} if $\Xi(A,\tau)$ is a topological $\ell$-group.
\end{definition}

By Proposition~\ref{prop:topological-lgroup-criterion}, a topological MV-algebra
is strongly summable precisely when it can be embedded as the unit interval of a
topological $\ell$-group.  We now give a useful sufficient condition for a
topological MV-algebra to be strongly summable.

\begin{proposition}\label{prop:square-zero-strongly-summable}
	Let $(A, \tau)$ be a topological MV-algebra and let $\mathcal{V}(0)$ be a
	neighbourhood basis of $0$.  Suppose that there exists a $0$-neighbourhood
	$E\in\mathcal V(0)$ such that
	\[
	\forall x\in E,\qquad x\odot x=0.
	\]
	Then $(A,\tau)$ is strongly summable.
\end{proposition}

\begin{proof}
	By Proposition~\ref{prop:U-odot-U-zero}, $(A,\tau)$ is summable.  Hence
	$\Xi(A,\tau)=(G_A,u_{G_A},\tau_{G_A})$ is defined.  As before, it is enough to
	prove that the positive-part map $p:G_A\to G_A$, $p(g)=g^+$, is continuous
	with respect to $\tau_{G_A}$.

	Fix $x\in G_A$.  By Proposition~\ref{prop:U-odot-U-zero}, it is enough to
	consider a basic $0$-neighbourhood of the form
	$(\varphi_A(U^\tau))^{\vartriangle}$, where $U^\tau$ is a $0$-neighbourhood
	in $(A,\tau)$.  Choose a $0$-neighbourhood $W^\tau$ in $A$ such that
	$W^\tau\ominus W^\tau\subseteq U^\tau$.

	We first record an auxiliary observation.  With $[0,u_{G_A}]$ endowed with the
	topology transported from $(A,\tau)$ via $\varphi_A$, for every fixed
	$\gamma\in G_A$, the map
	$\Theta_\gamma:[0,u_{G_A}]\times[0,u_{G_A}]\to[0,u_{G_A}]$, defined by
	\[
	\Theta_\gamma(a,b)=((\gamma+a-b)^+)\wedge u_{G_A},
	\]
	is continuous.  Indeed, since $u_{G_A}$ is a strong unit, there exists
	$r\in\mathbb N$ such that $|\gamma|\leq r u_{G_A}$.  Hence
	$\gamma^+\leq r u_{G_A}$ and $\gamma^-\leq r u_{G_A}$.  By iterating the identity
	$z=(z\wedge u_{G_A})+(z-u_{G_A})^+$ for $z\geq0$, both $\gamma^+$ and
	$\gamma^-$ can be written as finite sums of elements of $[0,u_{G_A}]$.  Thus
	there exist $a_1,a_2,\ldots,a_m,b_1,b_2,\ldots,b_n\in A$ such that
	\[
	\gamma=\varphi_A(a_1)+\varphi_A(a_2)+\cdots+\varphi_A(a_m)
	      -\varphi_A(b_1)-\varphi_A(b_2)-\cdots-\varphi_A(b_n).
	\]
	The corresponding function on $[0,1]^{m+n+2}$ is
	\[
	\begin{aligned}
		&(t,s,t_1,t_2,\ldots,t_m,s_1,s_2,\ldots,s_n)
		\longmapsto
		\min\!\Bigl(
		\max\!\bigl(
		t+\sum_{i=1}^{m}t_i-\bigl(s+\sum_{j=1}^{n}s_j\bigr),0
		\bigr),1
		\Bigr),
	\end{aligned}
	\]
	which is a McNaughton function.  By McNaughton's theorem
	(Theorem~\ref{thm:mcnaughton}), it is induced by an MV-term.  After
	substituting the parameters $a_1,a_2,\ldots,a_m,b_1,b_2,\ldots,b_n$, the
	corresponding map $A\times A\to A$ is an MV-term map with parameters.
	Transporting it through $\varphi_A$ gives $\Theta_\gamma$.  Hence
	$\Theta_\gamma$ is continuous.

	For $a,b\in[0,u_{G_A}]$, define
	\[
	\alpha(a,b)=(x+a-b)^+-(x-b)^+,
	\qquad
	\beta(b)=x^+-(x-b)^+.
	\]
	By the first assertion of Lemma~\ref{lem:lgroup-positive-part-estimate},
	applied to the pairs $(x-b,x+a-b)$ and $(x-b,x)$, respectively, we have
	$0\leq\alpha(a,b)\leq a$ and $0\leq\beta(b)\leq b$.  Since
	$a,b\in[0,u_{G_A}]$, it follows that $\alpha(a,b),\beta(b)\in[0,u_{G_A}]$.

	We claim that $\alpha$ is continuous at $(0,0)$, and that $\beta$ is
	continuous at $0$, where $[0,u_{G_A}]$ is equipped with the topology transported
	from $(A,\tau)$ via $\varphi_A$.

	First consider $\beta$.  By Lemma~\ref{lem:lgroup-positive-part-estimate},
	applied to the pair $(x-b,x)$, we have
	\[
	\beta(b)=x^+-(x-b)^+=b\wedge x^+.
	\]
	Since $b\leq u_{G_A}$, this equals $b\wedge(x^+\wedge u_{G_A}$.  Now
	$x^+\wedge u_{G_A}\in[0,u_{G_A}]$.  Since both arguments lie in $[0,u_{G_A}]$, the map
	$b\mapsto b\wedge(x^+\wedge u_{G_A})$ is, via $\varphi_A$, an MV-term map with
	parameter.  Therefore $\beta$ is continuous as a map
	$[0,u_{G_A}]\to[0,u_{G_A}]$.

	Next consider $\alpha$.  Applying the same identity to the pair
	$(x-b,x+a-b)$, we get
	\[
	\alpha(a,b)=(x+a-b)^+-(x-b)^+
	=a\wedge(x+a-b)^+.
	\]
	Since $a\leq u_{G_A}$, this can be written as
	\[
	\alpha(a,b)=a\wedge\bigl(((x+a-b)^+)\wedge u_{G_A}\bigr).
	\]
	The second factor is exactly $\Theta_x(a,b)$, which is continuous by the
	auxiliary observation above.  Since both $a$ and $\Theta_x(a,b)$ lie in
	$[0,u_{G_A}]$, the operation $(r,s)\mapsto r\wedge s$ used here is, via
	$\varphi_A$, an MV-term operation.  It follows that $\alpha$ is continuous
	at $(0,0)$.

	Since $\alpha(0,0)=0$ and $\beta(0)=0$, by continuity we can choose a
	$0$-neighbourhood $V^\tau$ in $A$ such that
	\[
	\alpha(a,b)\in\varphi_A(W^\tau)\quad\text{and}\quad\beta(b)\in\varphi_A(W^\tau)
	\]
	whenever $a,b\in\varphi_A(V^\tau)$.

	We claim that
	\[
	h\in(\varphi_A(V^\tau))^{\vartriangle}
	\quad\Longrightarrow\quad
	(x+h)^+-x^+\in(\varphi_A(U^\tau))^{\vartriangle}.
	\]
	Let $h\in(\varphi_A(V^\tau))^{\vartriangle}$.  Put $a=h^+$ and $b=h^-$.  Then
	$a,b\in\varphi_A(V^\tau)$, and $h=a-b$.  Hence
	\[
	(x+h)^+-x^+=(x+a-b)^+-x^+=\alpha(a,b)-\beta(b).
	\]
	By the choice of $V^\tau$, we have
	$\alpha(a,b),\beta(b)\in\varphi_A(W^\tau)$.  Therefore
	\[
	\bigl((x+h)^+-x^+\bigr)^+
	=
	\alpha(a,b)\ominus\beta(b)
	\in \varphi_A(W^\tau\ominus W^\tau)\subseteq \varphi_A(U^\tau),
	\]
	and similarly
	\[
	\bigl((x+h)^+-x^+\bigr)^-
	=
	\beta(b)\ominus\alpha(a,b)
	\in \varphi_A(W^\tau\ominus W^\tau)\subseteq \varphi_A(U^\tau).
	\]
	Thus $(x+h)^+-x^+\in(\varphi_A(U^\tau))^{\vartriangle}$.

	This proves continuity of $g\mapsto g^+$ at the arbitrary point $x\in G_A$.
	Hence the positive-part map is continuous on $G_A$, and therefore
	$(G_A,u_{G_A},\tau_{G_A})$ is a topological $\ell$-group.
\end{proof}

By Proposition~\ref{prop:square-zero-strongly-summable}, totally ordered
topological MV-algebras form an important class of strongly summable
topological MV-algebras.

\begin{corollary}\label{cor:totally-ordered-strongly-summable}
	Every totally ordered topological MV-algebra is strongly summable.
\end{corollary}

\begin{proof}
	If the topology is not trivial, then by Corollary~\ref{cor-totallyorderedMV},
	there exists a $0$-neighbourhood $U$ such that
	$x\odot x=0$ for every $x\in U$.  Hence the result follows from
	Proposition~\ref{prop:square-zero-strongly-summable}.

	If the topology is trivial, then the trivial topology on $G_A$ makes
	$(G_A,u_{G_A})$ into a topological $\ell$-group and induces the trivial
	topology on $[0,u_{G_A}]$.  Hence $(A,\tau)$ is strongly summable.
\end{proof}

The notion of strong summability allows us to isolate the part of the preceding
adjunction that takes values in topological $\ell$-groups.  We now make this
restriction categorical.

We denote by
\[
\mathcal{T}\!\mathit{op}^{\mathrm{str}}\!\mathcal{MV}
\]
the full subcategory of
$\mathcal{T}\!\mathit{op}^{+}\!\mathcal{MV}$
whose objects are strongly summable topological MV-algebras.  Thus the objects
of $\mathcal{T}\!\mathit{op}^{\mathrm{str}}\!\mathcal{MV}$ are precisely
those summable topological MV-algebras $(A,\tau)$ for which
\[
\Xi(A,\tau)=(G_A,u_{G_A},\tau_{G_A})
\]
is a topological $\ell$-group.

We also denote by
\[
\mathcal{T}\!\mathit{op}\ell\mathcal{G}_u
\]
the full subcategory of
$\ell\mathcal{T}\mathcal{G}_u^{\mathcal{MV}}$
whose objects are triples $(G,u,\tau_G)$, where $G$ is an abelian
$\ell$-group, $u\in G$ is a positive strong unit, and $\tau_G$ is a group
topology such that $(G,u,\tau_G)$ is a topological $\ell$-group.

By Proposition~\ref{prop:topological-lgroup-criterion}, the restrictions
\[
\Gamma|_{\mathcal{T}\!\mathit{op}\ell\mathcal{G}_u}:
\mathcal{T}\!\mathit{op}\ell\mathcal{G}_u
\longrightarrow
\mathcal{T}\!\mathit{op}^{\mathrm{str}}\!\mathcal{MV}
\]
and
\[
\Xi|_{\mathcal{T}\!\mathit{op}^{\mathrm{str}}\!\mathcal{MV}}:
\mathcal{T}\!\mathit{op}^{\mathrm{str}}\!\mathcal{MV}
\longrightarrow
\mathcal{T}\!\mathit{op}\ell\mathcal{G}_u
\]
are well defined.  Thus Theorem~\ref{thm:adjunction} immediately yields the
corresponding restricted adjunction.

\begin{corollary}\label{cor:strong-adjunction}
	$\Xi \dashv \Gamma: \mathcal{T}\!\mathit{op}^{\mathrm{str}}\!\mathcal{MV} \rightharpoonup \mathcal{T}\!\mathit{op}\ell\mathcal{G}_u.$
\end{corollary}

Example~\ref{ex:non-unique-ambient-topology} shows that the adjunction in
Corollary~\ref{cor:strong-adjunction} is not an equivalence of categories in
general.  Indeed, the functor $\Gamma$ need not be full.

\begin{example}\label{ex:non-unique-ambient-topology}
	Let
	\[
	G=\{\mathbf{x}=(x_n)_{n\in\mathbb N}\in\mathbb R^{\mathbb N}:
	\sup_n |x_n|<\infty\},
	\]
	with coordinatewise addition and coordinatewise order, and let
	$\mathbf{u}=(1,1,1,\ldots)$.  Then $\mathbf{u}$ is a strong unit of
	$G$, and $[\mathbf{0},\mathbf{u}]=[0,1]^{\mathbb N}.$

	Let $\sigma_1$ be the subspace topology inherited from the product topology
	on $\mathbb R^{\mathbb N}$.  Thus a $0$-neighbourhood basis of
	$\sigma_1$ is given by
	\[
	V(F,\varepsilon)
	=
	\{\mathbf{x}\in G: |x_n|<\varepsilon \text{ for all } n\in F\},
	\]
	where $F\subseteq\mathbb N$ is finite and $\varepsilon>0$.

	Define
	\[
	p(\mathbf{x})=\sup_{n\in\mathbb N}\frac{|x_n|}{n},
	\qquad
	\mathbf{x}=(x_n)\in G,
	\]
	and let $\sigma_2$ be the topology induced by the metric
	$d_p(\mathbf{x},\mathbf{y})=p(\mathbf{x}-\mathbf{y})$.
	Thus a $0$-neighbourhood basis for $\sigma_2$ is given by
	\[
	N_\varepsilon=\{\mathbf{x}\in G:p(\mathbf{x})<\varepsilon\},
	\qquad
	\varepsilon>0.
	\]

	Both $(G,\mathbf{u},\sigma_1)$ and $(G,\mathbf{u},\sigma_2)$ are
	topological $\ell$-groups.  This is clear for $\sigma_1$.  For
	$\sigma_2$, the inequalities
	\[
	p(\mathbf{x}+\mathbf{y})\leq p(\mathbf{x})+p(\mathbf{y}),\qquad
	p(-\mathbf{x})=p(\mathbf{x}),\qquad
	p(\mathbf{x}^+-\mathbf{y}^+)\leq p(\mathbf{x}-\mathbf{y})
	\]
	show continuity of the group operations and of the positive-part map.

	The topology $\sigma_2$ is strictly finer than $\sigma_1$.  Indeed, each
	coordinate projection is continuous with respect to $\sigma_2$, since
	$p(\mathbf{x})<\varepsilon/m$ implies $|x_m|<\varepsilon$.  Hence
	$\sigma_1\subseteq\sigma_2$.  The inclusion is strict: if
	$\mathbf{e}_k$ denotes the $k$-th standard basis vector and
	$\mathbf{x}^{(k)}=k\mathbf{e}_k$, then
	$\mathbf{x}^{(k)}\to\mathbf{0}$ in $\sigma_1$, while
	$p(\mathbf{x}^{(k)})=1$ for every $k$.  Thus
	$\mathbf{x}^{(k)}\not\to\mathbf{0}$ in $\sigma_2$.

	However, the two topologies induce the same topology on
	$[\mathbf{0},\mathbf{u}]$.  Since $\sigma_1\subseteq\sigma_2$, it remains
	to prove that every $\sigma_2$-neighbourhood in
	$[\mathbf{0},\mathbf{u}]$ contains a $\sigma_1$-neighbourhood.

	Fix $\mathbf{x}\in[0,1]^{\mathbb N}$ and $\varepsilon>0$.  Choose
	$N\in\mathbb N$ such that $1/n<\varepsilon$ for all $n>N$, and choose
	$\delta>0$ such that $\delta/n<\varepsilon$ for $n=1,2,\ldots,N$.  If
	$\mathbf{y}\in[0,1]^{\mathbb N}$ satisfies $|y_n-x_n|<\delta$ for
	$n=1,2,\ldots,N$, then $p(\mathbf{y}-\mathbf{x})<\varepsilon$.  Hence every
	$\sigma_2$-neighbourhood in $[\mathbf{0},\mathbf{u}]$ contains a
	$\sigma_1$-neighbourhood.  Therefore
	\[
	\sigma_1|_{[\mathbf{0},\mathbf{u}]}
	=
	\sigma_2|_{[\mathbf{0},\mathbf{u}]}.
	\]

	Consequently, the identity map on $[\mathbf{0},\mathbf{u}]$ defines a
	continuous MV-homomorphism
	\[
	\Gamma(G,\mathbf{u},\sigma_1)\longrightarrow
	\Gamma(G,\mathbf{u},\sigma_2).
	\]
	However, this MV-homomorphism is not induced by any morphism
	\[
	(G,\mathbf{u},\sigma_1)\longrightarrow (G,\mathbf{u},\sigma_2)
	\]
	in $\mathcal{T}\!\mathit{op}\ell\mathcal{G}_u$.  Indeed, any such morphism
	inducing the identity on $[\mathbf{0},\mathbf{u}]$ must be the identity on
	$G$, since $[\mathbf{0},\mathbf{u}]$ generates $G$ as a group.  But the
	identity map
	\[
	(G,\sigma_1)\longrightarrow (G,\sigma_2)
	\]
	is not continuous, as shown above.  Therefore $\Gamma$ is not full in
	general.
\end{example}

Example~\ref{ex:non-unique-ambient-topology} shows that the ambient topology is
not determined by the unit-interval topology in general.  Nevertheless, for
certain topological MV-algebras the embedding into a topological $\ell$-group
is unique.  Proposition~\ref{prop:square-zero-unique-topology} gives a local
condition under which this uniqueness holds.

\begin{proposition}\label{prop:square-zero-unique-topology}
	Let $(G,u)$ be an abelian $\ell$-group with strong unit.  Let
	$\sigma_1$ and $\sigma_2$ be group topologies on $G$ such that both
	$(G,u,\sigma_1)$ and $(G,u,\sigma_2)$ are topological $\ell$-groups.  Assume that
	\[
	\sigma_1|_{[0,u]}=\sigma_2|_{[0,u]}.
	\]
	If this common topology on $[0,u]$ admits a $0$-neighbourhood $E$ such
	that
	\[
	x\odot x=0
	\qquad
	\text{for every }x\in E,
	\]
	then $\sigma_1=\sigma_2$.
\end{proposition}

\begin{proof}
	Put $\tau=\sigma_1|_{[0,u]}=\sigma_2|_{[0,u]}$.  We show that any group
	topology $\sigma$ on $G$ for which $(G,u,\sigma)$ is a topological
	$\ell$-group and $\sigma|_{[0,u]}=\tau$ is determined by $\tau$.

	We shall use the following elementary observation.  If $h\in G^+$ and
	$2(h\wedge u)\leq u$, then $h\leq u$.  Indeed, by Proposition
	\ref{prop:lgroup-basic}(5), it is enough to verify this in the totally ordered
	case.  There, if $h>u$, then $h\wedge u=u$, contradicting
	$2(h\wedge u)\leq u$.  Hence $h\leq u$.

	Let $O^\sigma$ be a $0$-neighbourhood in $(G,\sigma)$.  Choose a
	$0$-neighbourhood $O_1^\sigma$ in $(G,\sigma)$ such that
	$O_1^\sigma-O_1^\sigma\subseteq O^\sigma$.  Since
	$O_1^\sigma\cap[0,u]$ is a $0$-neighbourhood in $([0,u],\tau)$, choose a
	$0$-neighbourhood $U^\tau$ in $([0,u],\tau)$ such that
	\[
	U^\tau\subseteq E\cap O_1^\sigma\cap[0,u].
	\]

	Consider $(U^\tau)^{\vartriangle}$, that is,
	\[
	(U^\tau)^{\vartriangle}=\{g\in G:g^+\in U^\tau,\ g^-\in U^\tau\}.
	\]
	Then $(U^\tau)^{\vartriangle}\subseteq O^\sigma$.  Indeed, if
	$g\in (U^\tau)^{\vartriangle}$, then $g=g^+-g^-$, with
	$g^+,g^-\in U^\tau\subseteq O_1^\sigma$, and hence
	$g\in O_1^\sigma-O_1^\sigma\subseteq O^\sigma$.

	We next show that $(U^\tau)^{\vartriangle}$ is itself a $0$-neighbourhood in
	$(G,\sigma)$.  Since $(G,u,\sigma)$ is a topological $\ell$-group,
	the maps
	\[
	\theta_+:(G,\sigma)\longrightarrow([0,u],\tau),
	\qquad
	\theta_+(g)=g^+\wedge u,
	\]
	and
	\[
	\theta_-:(G,\sigma)\longrightarrow([0,u],\tau),
	\qquad
	\theta_-(g)=g^-\wedge u,
	\]
	are continuous.  Indeed, the corresponding maps into $(G,\sigma)$ are
	continuous and their images lie in $[0,u]$.

	Since $U^\tau$ is a $0$-neighbourhood in $([0,u],\tau)$ and
	$\theta_+(0)=\theta_-(0)=0$, the sets
	$\theta_+^{-1}(U^\tau)$ and $\theta_-^{-1}(U^\tau)$ are
	$0$-neighbourhoods in $(G,\sigma)$.  Hence
	\[
	W^\sigma
	=
	\theta_+^{-1}(U^\tau)\cap\theta_-^{-1}(U^\tau)
	=
	\{g\in G:g^+\wedge u\in U^\tau,\ g^-\wedge u\in U^\tau\}
	\]
	is a $0$-neighbourhood in $(G,\sigma)$.  We claim that
	$W^\sigma=(U^\tau)^{\vartriangle}$.

	The inclusion $(U^\tau)^{\vartriangle}\subseteq W^\sigma$ is clear.  Conversely, let
	$g\in W^\sigma$.  Since $g^+\wedge u\in U^\tau\subseteq E$, we have
	$(g^+\wedge u)\odot(g^+\wedge u)=0$, hence $2(g^+\wedge u)\leq u$.  By the
	observation above, $g^+\leq u$, and therefore $g^+=g^+\wedge u\in U^\tau$.
	Similarly, $g^-\in U^\tau$.  Thus $g\in (U^\tau)^{\vartriangle}$, and so
	$W^\sigma=(U^\tau)^{\vartriangle}$.  Therefore $(U^\tau)^{\vartriangle}$ is a
	$0$-neighbourhood in $(G,\sigma)$.

	We have shown that every $0$-neighbourhood in $(G,\sigma)$ contains a set
	of the form $(U^\tau)^{\vartriangle}$, where $U^\tau$ is a $0$-neighbourhood in
	$([0,u],\tau)$ with $U^\tau\subseteq E$.  Conversely, each such
	$(U^\tau)^{\vartriangle}$ is a $0$-neighbourhood in $(G,\sigma)$.  Hence these sets
	form a $0$-neighbourhood basis for $\sigma$, determined only by $\tau$.
	Applying this to $\sigma_1$ and $\sigma_2$, we obtain
	$\sigma_1=\sigma_2$.
\end{proof}

We next derive the corresponding uniqueness statement for totally ordered
$\ell$-groups.

\begin{corollary}\label{cor:totally-ordered-unique-topology}
	Let $(G,u)$ be a totally ordered abelian $\ell$-group with strong unit.
	Let $\sigma_1$ and $\sigma_2$ be group topologies on $G$ such that both
	$(G,u,\sigma_1)$ and $(G,u,\sigma_2)$ are topological $\ell$-groups.  If
	\[
	\sigma_1|_{[0,u]}=\sigma_2|_{[0,u]},
	\]
	then $\sigma_1=\sigma_2$.
\end{corollary}

\begin{proof}
	Put $\tau=\sigma_1|_{[0,u]}=\sigma_2|_{[0,u]}$.

	First suppose that $\tau$ is trivial.  We claim that any group topology
	$\sigma$ on $G$ such that $(G,u,\sigma)$ is a topological
	$\ell$-group and $\sigma|_{[0,u]}=\tau$ must be trivial.  Let
	$O^\sigma$ be a $0$-neighbourhood in $(G,\sigma)$.  Since
	$O^\sigma\cap[0,u]$ is a $0$-neighbourhood in $([0,u],\tau)$, and
	$\tau$ is trivial, we have $[0,u]\subseteq O^\sigma$.

	Let $g\in G$.  Since $u$ is a strong unit, there exists
	$n\in\mathbb N$ such that $g^+\leq nu$ and $g^-\leq nu$.  Thus $g^+$
	and $g^-$ are finite sums of elements of $[0,u]$.  By the continuity at
	$(0,0,\ldots,0)$ of the map
	\[
	(z_1,z_2,\ldots,z_n,w_1,w_2,\ldots,w_n)
	\longmapsto
	z_1+z_2+\cdots+z_n-w_1-w_2-\cdots-w_n,
	\]
	choose a $0$-neighbourhood $V^\sigma$ of $0$ in
	$(G,\sigma)$ such that
	\[
	\bigl(
	\underbrace{V^\sigma+V^\sigma+\cdots+V^\sigma}_{n\text{ terms}}
	\bigr)
	-
	\bigl(
	\underbrace{V^\sigma+V^\sigma+\cdots+V^\sigma}_{n\text{ terms}}
	\bigr)
	\subseteq O^\sigma.
	\]
	Since $V^\sigma\cap[0,u]$ is a $0$-neighbourhood in the trivial topology
	$\tau$, we have $[0,u]\subseteq V^\sigma$.
	Hence, writing $g^+$ and $g^-$ as such finite sums, we get
	$g=g^+-g^-\in O^\sigma$.  Since $g\in G$ was arbitrary, every
	$0$-neighbourhood in $(G,\sigma)$ is equal to $G$.  Therefore
	$\sigma$ is trivial.  Applying this to $\sigma_1$ and $\sigma_2$, we get
	$\sigma_1=\sigma_2$.

	Now suppose that $\tau$ is not trivial.  Then $([0,u],\tau)$ is a
	totally ordered topological MV-algebra.  By Corollary
	\ref{cor-totallyorderedMV}, there exists a $0$-neighbourhood $E$ in
	$([0,u],\tau)$ such that $x\odot x=0$ for every $x\in E$.  Hence the
	result follows from Proposition~\ref{prop:square-zero-unique-topology}.
\end{proof}

We formulate the totally ordered case categorically.  Denote by
\[
\mathcal{OT}\!\mathit{op}\mathcal{MV}
\]
the full subcategory of
$\mathcal{T}\!\mathit{op}^{\mathrm{str}}\!\mathcal{MV}$
whose objects are totally ordered topological MV-algebras.  Similarly, let
\[
\mathcal{OT}\!\mathit{op}\ell\mathcal{G}_u
\]
be the full subcategory of
$\mathcal{T}\!\mathit{op}\ell\mathcal{G}_u$
whose objects are totally ordered topological $\ell$-groups.

Corollaries~\ref{cor:totally-ordered-strongly-summable},
\ref{cor:strong-adjunction}, and~\ref{cor:totally-ordered-unique-topology}
together yield the following categorical equivalence.

\begin{corollary}\label{cor:totally-ordered-equivalence}
	On totally ordered objects, the adjunction $\Xi\dashv\Gamma$ of
	Corollary~\ref{cor:strong-adjunction} restricts to an equivalence
	\[
	\mathcal{OT}\!\mathit{op}\mathcal{MV}
	\simeq
	\mathcal{OT}\!\mathit{op}\ell\mathcal{G}_u.
	\]
\end{corollary}

Corollary~\ref{cor:totally-ordered-equivalence} thus provides a topological
refinement of Chang's representation theorem for totally ordered
MV-algebras~\cite{chang1959}.

We next establish the topological counterpart of Di Nola--Lettieri equivalence
for perfect MV-algebras (recall Definition~\ref{def:perfect-MV},
Theorem~\ref{thm:dinola-lettieri}, and the notation $\Sigma,\Delta,\kappa,\zeta$
introduced there; see also \cite{dinola1994}).

First, $\Sigma$ is lifted to abelian topological $\ell$-groups, and it
is shown to detect continuity of $\ell$-group homomorphisms.  Second, using
Corollary~\ref{cor-perfectMV},
Proposition~\ref{prop:square-zero-strongly-summable}, and
Proposition~\ref{prop:square-zero-unique-topology},
$\Delta(A)$ is equipped with a uniquely determined topological
$\ell$-group topology.
Finally, the algebraic natural isomorphisms $\kappa$ and $\zeta$ are
shown to become natural isomorphisms in the corresponding topological
categories.

Let $\mathcal{T}\!\mathit{op}^{\mathrm{nt}}\!\mathcal{PMV}$ denote
the category of perfect topological MV-algebras with non-trivial
topologies and continuous MV-homomorphisms, and let
$\mathcal{T}\!\mathit{op}\,\ell\mathcal G$ denote the category of
abelian topological $\ell$-groups and continuous $\ell$-group
homomorphisms.

For an abelian topological $\ell$-group $(H,\sigma)$, equip
$\mathbb Z$ with the discrete topology $\delta_{\mathbb Z}$, and
define
\[
\Sigma(H,\sigma):=
\left(
\Sigma(H),
\left.(\delta_{\mathbb Z}\times\sigma)\right|_{\Sigma(H)}
\right).
\]
The topology $\delta_{\mathbb Z}\times\sigma$ makes
$\mathbb Z\times_{\mathrm{lex}}H$ a topological $\ell$-group, so
$\Sigma(H,\sigma)$ is a perfect topological MV-algebra.  Moreover,
$\operatorname{Rad}(\Sigma(H))=\{0\}\times H^+$ and
$\neg\operatorname{Rad}(\Sigma(H))=\{1\}\times(-H^+)$ are non-empty
open subsets.  Hence the induced topology is non-trivial.

If $f:(H,\sigma_H)\to(K,\sigma_K)$ is a continuous
$\ell$-group homomorphism, then
$\operatorname{id}_{\mathbb Z}\times f$ is continuous, and so is its
restriction $\Sigma(f)$.  Thus $\Sigma$ lifts to a functor
\[
\Sigma:
\mathcal{T}\!\mathit{op}\,\ell\mathcal G
\longrightarrow
\mathcal{T}\!\mathit{op}^{\mathrm{nt}}\!\mathcal{PMV}.
\]

\begin{lemma}\label{lem:Sigma-detects-continuity}
	Let $f:H\to K$ be an $\ell$-group homomorphism between abelian
	topological $\ell$-groups.  Then $f$ is continuous if and only if
	$\Sigma(f)$ is continuous.
\end{lemma}

\begin{proof}
	If $f$ is continuous, then
	$\operatorname{id}_{\mathbb Z}\times f$, and hence its restriction
	$\Sigma(f)$, is continuous.

	Conversely, suppose that $\Sigma(f)$ is continuous.  Restricting it
	to the radicals shows that
	$f|_{H^+}:H^+\to K^+$ is continuous.  Let $O$ be a
	$0$-neighbourhood in $K$, and choose a $0$-neighbourhood
	$O_1$ such that $O_1-O_1\subseteq O$.  There exists a
	$0$-neighbourhood $V$ in $H^+$ such that
	$f(V)\subseteq O_1\cap K^+$.

	Since the maps $h\mapsto h^+$ and $h\mapsto h^-=(-h)^+$ from
	$H$ to $H^+$ are continuous, the set
	$V^{\vartriangle}:=\{h\in H:h^+\in V,\ h^-\in V\}$, being the
	intersection of the inverse images of $V$ under these maps, is a
	$0$-neighbourhood in $H$.  Since $f$ preserves positive and
	negative parts, $f(h)=f(h^+)-f(h^-)$, and hence
	$f(V^{\vartriangle})\subseteq O_1-O_1\subseteq O$.  Thus $f$ is
	continuous.
\end{proof}

With this continuity-detecting property, we can promote $\Sigma$ to
a categorical equivalence, the third main theorem of this paper.

\begin{theorem}\label{thm:perfect-topological-equivalence}
	The functor
	\[
	\Sigma:
	\mathcal{T}\!\mathit{op}\,\ell\mathcal G
	\longrightarrow
	\mathcal{T}\!\mathit{op}^{\mathrm{nt}}\!\mathcal{PMV}
	\]
	is an equivalence of categories.
\end{theorem}

\begin{proof}
	Let
	$(A,\tau)\in \mathcal{T}\!\mathit{op}^{\mathrm{nt}}\!\mathcal{PMV}$,
	and set
	$\kappa_A(\tau):=\{\kappa_A(U):U\in\tau\}$.  Then
	\[
	\kappa_A:
	(A,\tau)\longrightarrow
	(\Sigma(\Delta(A)),\kappa_A(\tau))
	\]
	is a homeomorphism.

	By Corollary~\ref{cor-perfectMV} and
	Proposition~\ref{prop:square-zero-strongly-summable}, there exists a
	topological $\ell$-group topology $\rho_A$ on the Mundici group of
	$\Sigma(\Delta(A))$ inducing $\kappa_A(\tau)$ on its unit interval.
	By Proposition~\ref{prop:square-zero-unique-topology}, this topology is
	unique.  Identifying this Mundici group with
	$\mathbb Z\times_{\mathrm{lex}}\Delta(A)$, we regard
	$\rho_A$ as a topology on this lexicographic product.

	Let $H_0:=\{0\}\times\Delta(A)$.  We claim that $H_0$ is open.
	By Corollary~\ref{cor-perfectMV} and
	\eqref{eq:perfect-radical-square-zero}, there is a
	$0$-neighbourhood $U$ in
	$(\Sigma(\Delta(A)),\kappa_A(\tau))$ such that
	$U\subseteq\{0\}\times\Delta(A)^+$.

	Consider the continuous map
	$p(g):=g^+\wedge(1,0)$, viewed as a map into the unit interval.
	There exists a symmetric open $0$-neighbourhood $W$ in
	$(\mathbb Z\times_{\mathrm{lex}}\Delta(A),\rho_A)$ such that
	$p(W)\subseteq U$.  Let $g=(n,h)\in W$.  If $n=1$, then
	$p(g)=(1,h\wedge 0)$, while if $n>1$, then $p(g)=(1,0)$.
	Thus $n>0$ would imply that $p(g)$ has first coordinate $1$,
	contrary to $p(g)\in U$.  Hence $n\leq 0$.  Applying the same
	argument to $-g\in W$ gives $n\geq 0$, and therefore $n=0$.
	Thus $W\subseteq H_0$, so $H_0$ is open.

	Let $\sigma_A$ be the topology on $\Delta(A)$ corresponding to
	the subspace topology of $H_0$ under $h\mapsto(0,h)$.  Since
	$H_0$ is an open $\ell$-subgroup,
	$(\Delta(A),\sigma_A)$ is a topological $\ell$-group.  Every
	coset $\{n\}\times\Delta(A)$ is open and is homeomorphic to
	$H_0$ by translation.  Consequently, the sets
	$\{n\}\times V$, where $n\in\mathbb Z$ and $V\in\sigma_A$,
	form a basis for $\rho_A$.  Hence
	$\rho_A=\delta_{\mathbb Z}\times\sigma_A$, and therefore
	\[
	\left.(\delta_{\mathbb Z}\times\sigma_A)
	\right|_{\Sigma(\Delta(A))}
	=
	\kappa_A(\tau).
	\]

	The uniqueness of $\rho_A$ implies the uniqueness of $\sigma_A$.
	Indeed, any other topological $\ell$-group topology
	$\sigma_A'$ inducing $\kappa_A(\tau)$ on
	$\Sigma(\Delta(A))$ would make
	$\delta_{\mathbb Z}\times\sigma_A'$ another topological
	$\ell$-group topology inducing the same unit-interval topology.
	Hence
	$\delta_{\mathbb Z}\times\sigma_A'
	=\delta_{\mathbb Z}\times\sigma_A$, and restricting to $H_0$
	gives $\sigma_A'=\sigma_A$.

	Define
	$\Delta(A,\tau):=(\Delta(A),\sigma_A)$.  If
	$f:(A,\tau_{G_A})\to(B,\tau_B)$ is a continuous MV-homomorphism, then by the naturality of $\kappa$, as expressed by the first diagram in \eqref{diag:naturality-kappa-zeta},
	\[
	\Sigma(\Delta(f))\circ\kappa_A=\kappa_B\circ f.
	\]
	Since
	$\kappa_A$ and $\kappa_B$ are homeomorphisms,
	$\Sigma(\Delta(f))$ is continuous.  By
	Lemma~\ref{lem:Sigma-detects-continuity}, $\Delta(f)$ is
	continuous.  Thus $\Delta$ lifts to a functor
	\[
	\Delta:
	\mathcal{T}\!\mathit{op}^{\mathrm{nt}}\!\mathcal{PMV}
	\longrightarrow
	\mathcal{T}\!\mathit{op}\,\ell\mathcal G.
	\]

	By the definition of $\sigma_A$, every component $\kappa_A$ is
	a homeomorphism.  Since the naturality of $\kappa$ already holds
	algebraically, $\kappa$ is a natural isomorphism in
	$\mathcal{T}\!\mathit{op}^{\mathrm{nt}}\!\mathcal{PMV}$.

	Let $(H,\sigma)$ be an abelian topological $\ell$-group.  By the
	first identity in \eqref{eq:DL-triangle-identities},
	$\Sigma(\zeta_H)=\kappa_{\Sigma(H)}^{-1}$, and hence
	$\Sigma(\zeta_H)$ is a homeomorphism.  Since $\zeta$ is a natural
	isomorphism at the algebraic level, $\zeta_H$ is an
	$\ell$-group isomorphism.  Lemma~\ref{lem:Sigma-detects-continuity}
	implies that $\zeta_H$ is continuous.  Moreover, functoriality of
	$\Sigma$ gives
	$\Sigma(\zeta_H^{-1})=\Sigma(\zeta_H)^{-1}$, so the same lemma
	implies that $\zeta_H^{-1}$ is continuous.  Therefore
	$\zeta_H$ is a topological $\ell$-isomorphism.

	Thus $\kappa$ and $\zeta$ are natural isomorphisms in the
	corresponding topological categories.  Therefore $\Delta$ is a
	quasi-inverse of $\Sigma$, and hence $\Sigma$ is an equivalence
	of categories.
\end{proof}

This completes the topological lifting of both foundational MV-algebra
equivalences.  The proof of Theorem~\ref{thm:perfect-topological-equivalence}
exhibits the same categorical-lifting pattern as the algebraic--topological
duality for incline semimodules in~\cite{bai2020}: an algebraic equivalence
($\Sigma\dashv\Delta$) is promoted to topological categories by equipping
the object part with the unique group topology induced on the unit interval,
and the algebraic natural isomorphisms $\kappa,\zeta$ are verified to remain
natural in the topological setting.

\section{Topological consequences of summability}
\label{sec:application}

Drawing on the framework developed in the preceding sections, and particularly the topological lifting of Section~\ref{sec:adjunction}, we now derive from summability several topological consequences: separation axioms, metrizability, second
countability, and approximate completeness.  As consequences of
summability, we prove that $T_0$ totally ordered topological MV-algebras
are completely regular, establish a Birkhoff--Kakutani-type metrization
theorem, and extend a key continuity lemma from the standard unit interval
to arbitrary summable topological MV-algebras, yielding a broad semantic
version of approximate completeness for many-valued logics.  We begin with
a lemma used repeatedly below.

\begin{lemma}\label{lem:lifting-properties}
	Let $(A,\tau)$ be a topological MV-algebra, let $\mathcal{V}(0)$ be a
	neighbourhood basis of $0$, and suppose there exists $U \in \mathcal{V}(0)$
	satisfying condition \eqref{cond:odot}.  Let $(G_A,u_{G_A}, \tau_{G_A}) = \Xi(A,\tau)$.
	Then the following hold:
	\begin{enumerate}[label=(\arabic*)]
		\item if $(A,\tau)$ is $T_0$, then $(G_A,u_{G_A}, \tau_{G_A})$ is $T_0$;
		\item if $(A,\tau)$ is first-countable, then $(G_A,u_{G_A}, \tau_{G_A})$ is first-countable.
	\end{enumerate}
\end{lemma}

\begin{proof}
	By Proposition \ref{prop:U-odot-U-zero}, $(A,\tau)$ is summable and
	$\{ (\varphi_A(V))^{\vartriangle} \mid V \in \mathcal{V}(0) \}$ forms a
	neighbourhood basis of $0$ in $\tau_{G_A}$.

	If $(A,\tau)$ is $T_0$, to prove that $(G_A,u_{G_A}, \tau_{G_A})$ is $T_0$
	it suffices to separate $0$ from any $z \in G_A$ with $z \neq 0$ in $\tau_{G_A}$.

	If $|z| \le u_{G_A}$ and $z \neq 0$, then $z^+ \neq 0$ or $z^- \neq 0$.
	We may assume without loss of generality that $z^+ \neq 0$.
	Then $z^+ = \varphi_A(x)$ for some $x \neq 0$.  Since $(A,\tau)$ is $T_0$,
	there exists an open set separating $0$ and $x$.

	If there exists $V_1 \in \mathcal{V}(0)$ such that $0 \in V_1$ and $x \notin V_1$,
	then $z^+ \notin \varphi_A(V_1)$, so $0 \in (\varphi_A(V_1))^{\vartriangle}$
	but $z \notin (\varphi_A(V_1))^{\vartriangle}$.

	If there exists an open set $W$ in $\tau$ such that $0 \notin W$ and $x \in W$,
	then by Proposition \ref{prop:starcon}, $\varphi_A$ is a homeomorphism, so
	$\varphi_A(W)$ is an open neighbourhood of $z^+$ in
	$([0,u_{G_A}], \tau_{G_A}|_{[0,u_{G_A}]})$.  Hence there exists $V_2 \in \mathcal{V}(0)$
	such that $(z^+ + (\varphi_A(V_2))^{\vartriangle}) \cap [0,u_{G_A}] \subseteq \varphi_A(W)$.
	Then $0 \notin z^+ + (\varphi_A(V_2))^{\vartriangle}$, i.e.,
	$-z^+ \notin (\varphi_A(V_2))^{\vartriangle}$, and consequently
	$z \notin (\varphi_A(V_2))^{\vartriangle}$ while $0 \in (\varphi_A(V_2))^{\vartriangle}$.
	Thus $(\varphi_A(V_2))^{\vartriangle}$ separates $0$ and $z$.

	If $|z| \not\le u_{G_A}$, take any $V_3 \in \mathcal{V}(0)$.  If
	$z \in (\varphi_A(V_3))^{\vartriangle}$, then
	$z^+, z^- \in \varphi_A(V_3) \subseteq [0,u_{G_A}]$, so
	$|z| = z^+ \vee z^- \le u_{G_A}$.  Hence $z \notin (\varphi_A(V_3))^{\vartriangle}$,
	while $0 \in (\varphi_A(V_3))^{\vartriangle}$.  Therefore
	$(\varphi_A(V_3))^{\vartriangle}$ separates $0$ and $z$.

	Thus $(G_A,u_{G_A},\tau_{G_A})$ is $T_0$.

	If $(A,\tau)$ is first-countable, we can choose $\mathcal{V}(0)$ to be countable.
	Then $\{ (\varphi_A(V))^{\vartriangle} \mid V \in \mathcal{V}(0) \}$ is a countable
	neighbourhood basis of $0$ in $\tau_{G_A}$, and hence
	$(G_A,u_{G_A},\tau_{G_A})$ is first-countable.  This completes the proof.
\end{proof}

With Lemma~\ref{lem:lifting-properties} in hand, we first consider separation axioms for
totally ordered topological MV-algebras. 

\begin{proposition}\label{prop:separation}
	Let $(A,\tau)$ be a totally ordered topological MV-algebra.  Then the following are equivalent:
	\begin{enumerate}[label=(\arabic*)]
		\item $(A,\tau)$ is $T_0$;
		\item $(A,\tau)$ is $T_1$;
		\item $(A,\tau)$ is Hausdorff ($T_2$);
		\item $(A,\tau)$ is regular ($T_3$);
		\item $(A,\tau)$ is completely regular ($T_{3\frac{1}{2}}$).
	\end{enumerate}
\end{proposition}

\begin{proof}
	We only need to show that (1) implies (5); the other implications follow
	directly from the definitions.

	Assume $(A,\tau)$ is $T_0$.  By Corollary~\ref{cor-totallyorderedMV}, $(A,\tau)$ is summable and homeomorphic to $([0,u_{G_A}], \tau_{G_A}|_{[0,u_{G_A}]})$. Consequently, for any non-trivial totally ordered topological MV-algebra, the summability hypothesis in Proposition~\ref{prop:separation} is automatically satisfied; the proposition therefore asserts that every such algebra satisfying $T_0$ is $T_{3\frac12}$. Moreover, the proof of
	Corollary \ref{cor-totallyorderedMV} shows that there exists an open
	neighbourhood $U$ of $0$ in $\tau$ satisfying condition \eqref{cond:odot}.
	Hence Lemma \ref{lem:lifting-properties} applies and yields that
	$(G_A,u_{G_A},\tau_{G_A})$ is a $T_0$ topological group.

	By Proposition \ref{prop:topological-group-T0-Tychonoff}, every $T_0$
	topological group is completely regular.  Consequently, its subspace
	$([0,u_{G_A}], \tau_{G_A}|_{[0,u_{G_A}]})$ is completely regular, and therefore
	$(A,\tau)$ is completely regular.  This completes the proof.
\end{proof}

\noindent
Thus, for totally ordered topological MV-algebras, the hierarchy of
separation axioms collapses completely, extending the known results
\cite{hoo1997,gan2024} from $T_3$ to $T_{3\frac12}$.  Moreover, while
Li and Yang \cite{li2026} required local convexity to obtain complete
regularity, Proposition~\ref{prop:separation} shows that this assumption
is unnecessary in the totally ordered case.

\begin{example}\label{ex:non-locally-convex}
	Let $A = [0,1]$.  For $a\in\mathbb R$ and $n\in\mathbb N$, let
	\[
	B(a,1/n)=\{x\in\mathbb R\mid |x-a|<1/n\}.
	\]
	Define $U_n = B(0, 1/n) \cap \mathbb{Q} \cap [0,1]$ and set
	\[
	\mathcal{V} = \{ U_n \mid n \in \mathbb{N} \}.
	\]
	For each $U_n \in \mathcal{V}$, we have
	$U_n(a) = B(a, 1/n) \cap (a + \mathbb{Q}) \cap [0,1]$.
	One verifies directly that $\mathcal{V}$ satisfies conditions (i) and (ii)
	of Proposition \ref{thm:MV-topology2}; hence there exists a topological
	MV-algebra $(A,\tau)$ for which $\mathcal{V}$ is an open neighbourhood
	basis of $0$.  The sets $U_n$ contain no convex neighbourhood of $0$,
	so $\tau$ is not locally convex at $0$.  Nevertheless, $(A,\tau)$ is
	clearly $T_0$ and totally ordered, and hence completely regular by
	Proposition~\ref{prop:separation}.
\end{example}

\medskip
\noindent
Having illustrated the reach of Proposition~\ref{prop:separation} with a
concrete example, we now  establish a metrizability criterion for totally ordered topological
MV-algebras.

\begin{proposition}\label{prop:metrization}
	Let $(A,\tau)$ be a totally ordered topological MV-algebra.  Then $(A,\tau)$
	is metrizable if and only if it is Hausdorff and first-countable.
	Moreover, if $(A,\tau)$ is metrizable, there exists a compatible metric $d$
	such that $d(\neg x, \neg y) = d(x, y)$ for all $x, y \in A$.
\end{proposition}

\begin{proof}
	Necessity is clear: every metrizable space is Hausdorff and first-countable.

	Conversely, assume $(A,\tau)$ is Hausdorff and first-countable.
	By Corollary \ref{cor-totallyorderedMV}, $(A,\tau)$ is summable and
	homeomorphic to $([0,u_{G_A}], \tau_{G_A}|_{[0,u_{G_A}]})$.  Moreover, the proof
	of Corollary \ref{cor-totallyorderedMV} shows that there exists an open
	neighbourhood $U$ of $0$ satisfying condition \eqref{cond:odot}.
	Hence Proposition~\ref{prop:topological-group-T0-Tychonoff} and
	Lemma~\ref{lem:lifting-properties} imply that $(G_A,u_{G_A},\tau_{G_A})$ is
	Hausdorff and first-countable.

	By the Birkhoff--Kakutani theorem (Theorem \ref{thm:birkhoff-kakutani}),
	$(G_A,u_{G_A},\tau_{G_A})$ is metrizable and admits a compatible
	left-invariant metric $d_{G_A}$.  Since $G_A$ is abelian, taking
	$g = u_{G_A} - z_1 - z_2$, we have
	\[
	d_{G_A}(z_1, z_2) = d_{G_A}(g + z_1, g + z_2)
	= d_{G_A}(u_{G_A} - z_1, u_{G_A} - z_2).
	\]
	Thus the restriction of $d_{G_A}$ to $[0, u_{G_A}]$ is the required metric.
	This completes the proof.
\end{proof}

From Proposition \ref{prop:metrization} we immediately obtain the following
corollary.  The second-countability criterion below is the standard
metric-space translation; we record it because it fails without total order.

\begin{corollary}\label{cor:second-countable}
	Let $(A,\tau)$ be a totally ordered Hausdorff topological MV-algebra.
	Then $(A,\tau)$ is second-countable if and only if it is separable and metrizable.
\end{corollary}

\begin{proof}
	If $(A,\tau)$ is second-countable, then it is first-countable.  Since it is
	Hausdorff by hypothesis, Proposition \ref{prop:metrization} implies that it
	is metrizable.  It is also separable because every second-countable space is
	separable.  Conversely, if $(A,\tau)$ is separable and metrizable, then it
	is a separable metric space, and hence it is second-countable.
\end{proof}

\medskip
\noindent
We now apply the topological machinery to the proof theory of
\L ukasiewicz logic. 
In particular, we obtain a new proof of a key lemma
in Duo~\cite{duo2026}, which is a crucial step in his proof of the
approximate strong completeness of the hypersequent calculus
$\mathbf{G}\text{\L}\forall$ (see~\cite{baaz2010} for the original
development of this calculus and its approximate completeness).
Whereas Duo's original proof relies
specifically on the standard MV-algebra $[0,1]$, the topological machinery
developed above shows that the main continuity ingredient of the lemma
extends to arbitrary summable topological MV-algebras.  This suggests the
study of $A$-valued semantics for first-order \L ukasiewicz logic using
the categorical framework of this paper.

To make this precise, we recall the necessary syntactic notions and define
$A$-semantics for a general MV-algebra $A$.

The formulas of propositional \L ukasiewicz logic are built from an infinite
set $X$ of propositional variables with two connectives $\bot$ (0-ary) and
$\to$ (binary).  The set of all formulas from $X$ is denoted by
$\operatorname{Form}(X)$.  A sequent is an ordered pair of finite multisets of
formulas, written as $\mathcal{S} = \Gamma \Rightarrow \Delta$.  A hypersequent
is a (possibly infinite) multiset of sequents (see, e.g., \cite{metcalfe2005}).
For the hypersequent calculus for \L ukasiewicz logic and related systems, see
\cite{baaz2010, wei2025}.

\begin{definition}[$A$-semantics]\label{def:A-semantics}
Let $A$ be an MV-algebra.  An \emph{$A$-valuation} is a map
$v: \operatorname{Form}(X) \to A$ satisfying:
\begin{enumerate}
	\item[(1)] $v(\bot) = 1$;
	\item[(2)] $v(P \to Q) = v(Q) \ominus v(P)$.
\end{enumerate}
We extend $v$ to sequents and hypersequents as follows, taking values in the
Mundici group $G_A$.  To simplify notation, we identify $A$ with its image
$\varphi_A(A) \subseteq G_A$ via the natural isomorphism $\varphi_A$
(so $A \subseteq G_A$).
\begin{enumerate}
	\item[(3)] $v(\Gamma) = \sum_{P \in \Gamma} v(P)$, with $v(\emptyset) = 0$;
	\item[(4)] $v(\Gamma \Rightarrow \Delta) = v(\Delta) - v(\Gamma)$;
	\item[(5)] For each hypersequent $\mathcal{G}$ (possibly infinite),
		$v(\mathcal{G}) = \bigwedge_{\mathcal{S} \in \mathcal{G}} v(\mathcal{S})$,
\end{enumerate}
where the sum and difference are taken in $G_A$, and $\bigwedge$ denotes the
infimum in the lattice $G_A$.
\end{definition}

In the sequel, we will use the same symbol $v$ to denote both a valuation on
propositional variables and its induced extension to formulas, sequents, and
hypersequents; the intended meaning will be clear from the context.

With these definitions in place, we establish the fundamental continuity
property of evaluation maps.

\begin{proposition}\label{prop:ev-continuous}
	Let $(A,\tau)$ be a summable topological MV-algebra and let
	$(G_A,u_{G_A},\tau_{G_A}) = \Xi(A,\tau)$, where $\tau_{G_A}$ is defined by
	\eqref{eq:V_GA}.  Then for any sequent $\mathcal{S}$, the evaluation map
	\[
	\operatorname{ev}_{\mathcal{S}}: A^X \to G_A, \quad v \mapsto v(\mathcal{S})
	\]
	is continuous, where $A^X$ is equipped with the product topology.
\end{proposition}

\noindent
This continuity property is crucial for applying the machinery of functional
analysis to MV-algebraic semantics, a direction largely unexplored in the
literature.  In particular, it allows us to extend Duo's lemma from $[0,1]$
to arbitrary summable topological MV-algebras, yielding a broad semantic
version of approximate completeness for \L ukasiewicz logic.

\begin{proof}
	Since $(A,\tau)$ is summable, by Proposition \ref{prop:starcon} we may
	identify $A$ with the subspace $[0, u_{G_A}] \subseteq G_A$ via a homeomorphic
	embedding $\iota: A \hookrightarrow G_A$.  The map
	$+: A \times A \to G_A$ defined by $(x,y) \mapsto \iota(x) + \iota(y)$
	(with the group addition on $G_A$) is continuous, and similarly the
	subtraction map $-: A \times A \to G_A$,
	$(x,y) \mapsto \iota(x) - \iota(y)$, is continuous.

	Define the complexity $\operatorname{cp}(\mathcal{S})$ of a sequent
	$\mathcal{S}$ as the number of logical connectives ($\to$) occurring in
	$\mathcal{S}$.  We now prove the proposition by induction on
	$\operatorname{cp}(\mathcal{S})$.

	\emph{Base case:} $\operatorname{cp}(\mathcal{S}) = 0$, i.e.,
	$\mathcal{S}$ is atomic.  We prove this by a secondary induction on
	$|\Gamma| + |\Delta|$.

	If $\mathcal{S}$ is the sequent $\Rightarrow$, then
	$\operatorname{ev}_{\mathcal{S}}$ is the constant function $0$, hence
	continuous.  Note that
	\[
	\begin{aligned}
		\operatorname{ev}_{\Gamma \Rightarrow \Delta, \bot}
		&= \operatorname{ev}_{\Gamma \Rightarrow \Delta} + 1, \\
		\operatorname{ev}_{\Gamma \Rightarrow \Delta, P}
		&= \operatorname{ev}_{\Gamma \Rightarrow \Delta} + \iota \circ \pi_P, \\
		\operatorname{ev}_{\bot, \Gamma \Rightarrow \Delta}
		&= \operatorname{ev}_{\Gamma \Rightarrow \Delta} - 1, \\
		\operatorname{ev}_{P, \Gamma \Rightarrow \Delta}
		&= \operatorname{ev}_{\Gamma \Rightarrow \Delta} - \iota \circ \pi_P,
	\end{aligned}
	\]
	where $P \in X$ and $\pi_P : A^X \to A$, $\pi_P(v) = v(P)$, is the
	projection.  Since $\operatorname{ev}_{\Gamma \Rightarrow \Delta}$ is
	continuous by the induction hypothesis, and each of
	$\operatorname{ev}_{\Gamma \Rightarrow \Delta, \bot}$,
	$\operatorname{ev}_{\bot, \Gamma \Rightarrow \Delta}$,
	$\operatorname{ev}_{\Gamma \Rightarrow \Delta, P}$, and
	$\operatorname{ev}_{P, \Gamma \Rightarrow \Delta}$ is a composition of
	continuous maps, they are all continuous.  This completes the secondary
	induction on $|\Gamma| + |\Delta|$.

	\emph{Inductive step:} Suppose $\operatorname{cp}(\mathcal{S}) > 0$.
	Then $\mathcal{S}$ contains at least one occurrence of $\to$.  Write
	$\mathcal{S}$ as $\Gamma \Rightarrow P \to Q, \Delta$ or
	$\Gamma, P \to Q \Rightarrow \Delta$.  Note that
	\[
	\begin{aligned}
		\operatorname{ev}_{\Gamma \Rightarrow P \to Q, \Delta}
		&= \operatorname{ev}_{\Gamma \Rightarrow \Delta}
		+ \iota \circ \bigl( \operatorname{ev}_{\Rightarrow Q}
		\ominus \operatorname{ev}_{\Rightarrow P} \bigr), \\
		\operatorname{ev}_{\Gamma, P \to Q \Rightarrow \Delta}
		&= \operatorname{ev}_{\Gamma \Rightarrow \Delta}
		- \iota \circ \bigl( \operatorname{ev}_{\Rightarrow Q}
		\ominus \operatorname{ev}_{\Rightarrow P} \bigr).
	\end{aligned}
	\]
	By Lemma~\ref{lem:continuous-ops}, $\ominus : A \times A \to A$ is
	continuous.  Hence $\operatorname{ev}_{\Gamma \Rightarrow P \to Q, \Delta}$
	and $\operatorname{ev}_{\Gamma, P \to Q \Rightarrow \Delta}$ are
	compositions of continuous maps and thus continuous by the induction
	hypothesis.  This completes the induction on
	$\operatorname{cp}(\mathcal{S})$, and the proposition follows.
\end{proof}

Because Proposition~\ref{prop:ev-continuous} makes evaluation continuous for
every summable $(A,\tau)$, the continuity ingredient in the compactness
argument of Corollary~\ref{cor:duo-generalized} carries over verbatim to any
totally ordered Hausdorff summable $A$.  Thus our argument isolates the main
topological ingredient needed for a possible $A$-valued version of Duo's
lemma, rather than merely providing a new proof of the $[0,1]$ case.

\begin{corollary}[{\cite[Lemma~3.12]{duo2026}}]\label{cor:duo-generalized}
	Let $\mathcal{H}$ be a (possibly infinite) propositional hypersequent and
	$\epsilon > 0$.  If for every finite sub-hypersequent $\mathcal{H}'$ of
	$\mathcal{H}$ there exists a valuation $w \in [0,1]^{X}$ such that
	$w(\mathcal{H}') \ge \epsilon$, then there exists a valuation
	$v \in [0,1]^{X}$ such that $v(\mathcal{H}) \ge \epsilon$.
\end{corollary}

\begin{proof}
	By Proposition \ref{prop:U-odot-U-zero} and a direct computation, the
	standard MV-algebra $[0,1]$, equipped with its order topology, is a
	summable topological MV-algebra, and the corresponding $\tau_{G_A}$ is the
	standard topology on $\mathbb{R}$.  By Proposition
	\ref{prop:ev-continuous}, the evaluation map
	$\operatorname{ev}_{\mathcal{S}}: [0,1]^{X} \to \mathbb{R}$,
	$v \mapsto v(\mathcal{S})$, is continuous for every sequent $\mathcal{S}$.

	Hence
	$\{ \operatorname{ev}_{\mathcal{S}}^{-1}([\epsilon, +\infty)) \}_{\mathcal{S} \in \mathcal{H}}$
	is a family of closed subsets of $[0,1]^{X}$.  The assumption of the
	corollary implies that this family satisfies the finite intersection
	property.  Therefore
	$\bigcap_{\mathcal{S} \in \mathcal{H}} \operatorname{ev}_{\mathcal{S}}^{-1}([\epsilon, +\infty))$
	is nonempty by the compactness of $[0,1]^{X}$.  Let
	$v \in \bigcap_{\mathcal{S} \in \mathcal{H}} \operatorname{ev}_{\mathcal{S}}^{-1}([\epsilon, +\infty))$.
	Then
	$v(\mathcal{H}) = \inf_{\mathcal{S} \in \mathcal{H}} v(\mathcal{S})
	= \inf_{\mathcal{S} \in \mathcal{H}} \operatorname{ev}_{\mathcal{S}}(v) \ge \epsilon$,
	as required.
\end{proof}

\begin{remark}
	Corollary~\ref{cor:duo-generalized} recovers Duo's lemma for the standard
	MV-algebra $[0,1]$, but its proof via Proposition~\ref{prop:ev-continuous}
	reveals that the essential continuity ingredient depends only on
	summability.  Consequently, for any summable topological MV-algebra
	$(A,\tau)$, the evaluation map $\operatorname{ev}_{\mathcal{S}}$ remains
	continuous (Proposition~\ref{prop:ev-continuous}), so the same compactness
	argument applies whenever $A^X$ is compact.  This points towards an
	extension of the hypersequent calculus $\mathbf{G}\text{\L}\forall$ beyond
	the standard real-valued setting, to $A$-valued semantics over compact
	topological MV-algebras.
\end{remark}

\medskip
\noindent
In this section we have shown that summability — a purely neighbourhood-theoretic condition — entails a series of analytic and logical consequences: it collapses separation axioms, forces metrizability and second countability under natural hypotheses, and lifts a key continuity lemma from \L ukasiewicz logic to arbitrary summable topological MV-algebras. Together with the lifting construction of Section~3 and the adjunction framework of Section~4, these results complete the core structural picture of how summability governs the topology of MV-algebras. At the same time, they raise natural questions — about the necessity of these conditions, about possible non-summable counterparts, and about the broader logical reach of the continuity lemma — which point toward directions for further investigation.

\section{Conclusion and future directions}
\label{sec:conclusion}

We have solved the topological lifting problem for the Mundici equivalence and
the Di Nola--Lettieri equivalence, establishing that a topological MV-algebra
can be realised as the unit interval of a lattice-ordered topological group
precisely when it is summable, and that perfect topological MV-algebras with
non-trivial topologies correspond exactly to abelian topological $\ell$-groups.
The categorical framework developed here---an adjunction that becomes an
equivalence in the totally ordered case---provides a systematic interface
between the topology of MV-algebras and the topology of $\ell$-groups,
complementing existing dualities and representation theorems in the literature.

Several questions remain open and suggest directions for future work.

\smallskip
\noindent\textbf{(1) When is a topological MV-algebra summable?}
It is not a priori settled whether every topological MV-algebra is summable; the results of this paper show that any counterexample, if it exists, cannot be totally ordered or perfect 
 (Corollaries~\ref{cor-totallyorderedMV}
and~\ref{cor-perfectMV}).  Characterising the exact class of summable
topological MV-algebras---or, alternatively, finding a non-summable example---
would delineate the boundaries of the topological lifting theorem.
A related question is whether every summable topological MV-algebra is
strongly summable (see Proposition~\ref{prop:square-zero-strongly-summable}
and Corollary~\ref{cor:totally-ordered-strongly-summable}); resolving this
would clarify whether the lattice operations on $G_A$ are automatically
continuous whenever the MV-algebra topology lifts.

\smallskip
\noindent\textbf{(2) Extension to related algebraic structures.}
The correspondences of Georgescu and Iorgulescu~\cite{georgescu2001},
Dvure\v{c}enskij~\cite{dvurecenskij2002}, and Galatos and Tsinakis~\cite{galatos2005}
relate generalised MV-algebras to various classes of $\ell$-groups.
Adapting the topological lifting framework to these settings would extend the
reach of the present results.  (More ambitiously, defining summability for
topological quantum B-algebras in the sense of Rump and Yang~\cite{rump2014}
could lead to an adjunction between summable topological quantum B-algebras and
topological partially ordered groups internal to quantales, thereby connecting
the regularity of quantum B-spaces to the theory of logical quantales.)

\smallskip
\noindent\textbf{(3) Approximate completeness beyond the standard setting.}
Proposition~\ref{prop:ev-continuous} extends the continuity ingredient of
\cite[Lemma~3.12]{duo2026} from $[0,1]$ to arbitrary summable topological
MV-algebras, providing the main topological ingredient for an $A$-valued
version of the approximate strong completeness of the hypersequent calculus
$\mathbf{G}\text{\L}\forall$ .  A deeper question, however, remains open: whether the summability
condition itself is preserved under the hypersequent rules of $\mathbf{G}\text{\L}\forall$ , and
consequently whether a full topological completeness theorem---rather than
merely the approximate-strong fragment---can be established for the entire
calculus with truth values ranging over any summable topological MV-algebra.
Such a result would require a systematic analysis of how the neighbourhood
basis condition~(SC) interacts with the cut rule and the communication rule
of $\mathbf{G}\text{\L}\forall$ , and would constitute a genuine topological extension of the
algebraic completeness programme initiated in \cite{duo2026}.

\medskip
\noindent
These open questions indicate that summability, while already yielding a
complete solution to the topological lifting problem for the Mundici and
Di Nola--Lettieri equivalences, marks only the beginning of a broader
investigation into the interaction between MV-algebraic topology and
logical structure.  Viewed within the APAL lineage on representations of
logical algebras, the present work joins several parallel liftings:
Dubuc and Poveda~\cite{dubucpoveda2010} represent MV-algebras as global sections
of sheaves of MV-chains via classifying topoi;
Rump and Yang~\cite{rump2014} dualise quantum B-algebras to logical
quantales, with MV-algebras sitting inside the commutative fragment;
Bai and Yang~\cite{bai2020} lift the algebraic--topological duality of
incline semimodules by sheaf methods.
The route taken here is different in kind---rather than a sheaf or
quantale envelope, summability selects precisely those topological
MV-algebras that are unit intervals of lattice-ordered topological abelian
groups under the unique induced group topology, so that Mundici and
Di Nola--Lettieri equivalences survive intact.  The four approaches are
therefore complementary facets of one programme: lifting algebraic
representation theorems of many-valued logic to topological or
enriched settings.


\end{document}